\documentclass[]{amsart}

\usepackage[margin=1in]{geometry}
\usepackage{amsmath, amsthm, amssymb, mathrsfs, enumitem, mathtools, tikz-cd, setspace, hyperref, float}

\usepackage[backend=biber, style=numeric]{biblatex}
\usepackage[T1]{fontenc}
\usepackage{mlmodern}

\setlist[enumerate,1]{label={(\alph*)}}

\allowdisplaybreaks
\newcommand{\ident}{\mathbb{I}}
\newcommand{\Z}{\mathbb{Z}}

\newcommand{\C}{\mathbb{C}}
\newcommand{\Q}{\mathbb{Q}}
\newcommand{\Qp}{{\mathbb{Q}_p}}
\newcommand{\Qptwo}{{\mathbb{Q}_{p^2}}}
\newcommand{\Qpfour}{{\mathbb{Q}_{p^4}}}
\newcommand{\Qpf}{{\mathbb{Q}_{p^f}}}

\newcommand{\Zp}{{\mathbb{Z}_p}}
\newcommand{\Fpf}{{\mathbb{F}_{p^f}}}

\newcommand{\D}{\mathbf{D}}

\newcommand{\V}{\mathbf{V}}

\newcommand{\F}{\mathbb{F}}
\newcommand{\Ms}{\mathcal{M}}
\newcommand{\Mf}{\mathfrak{M}}

\newcommand{\pu}{[\![u]\!]}

\newcommand{\Rep}{\text{\normalfont Rep}}

\newcommand{\ind}{\text{\normalfont Ind}}
\newcommand{\semi}{\text{\normalfont ss}}
\newcommand{\et}{\text{\normalfont \'{e}t}}

\newcommand{\Fil}{\text{\normalfont Fil}}
\newcommand{\I}{\text{\normalfont I}}
\newcommand{\II}{\text{\normalfont II}}

\newcommand{\Diag}{\text{\normalfont Diag}}
\newcommand{\Gal}{\text{\normalfont Gal}}
\newcommand{\GL}{\text{\normalfont GL}}

\newcommand{\Hom}{\text{\normalfont Hom}}

\newcommand{\Mod}{\text{\normalfont Mod}}

\newcommand{\Mat}{\text{\normalfont Mat}}

\newcommand{\cris}{\text{\normalfont cris}}

\newcommand{\oh}{\mathscr O}

\newcommand{\Ts}{\mathcal T}
\newcommand{\Ss}{\mathcal S}

\newcommand{\Sig}{\mathfrak{S}}

\newcommand{\As}{\mathcal{A}}
\newcommand{\Af}{\mathfrak{A}}
\newcommand{\Bs}{\mathcal{B}}
\newcommand{\Ws}{\mathcal{W}}
\newcommand{\Xs}{\mathcal{X}}
\newcommand{\Ys}{\mathcal{Y}}
\newcommand{\Vs}{\mathcal{V}}
\newcommand{\Is}{\mathcal{I}}
\newcommand{\Ps}{\mathcal{P}}

\newcommand{\Bf}{\mathfrak{B}}

\newtheoremstyle{italics}{}{}{\itshape}{}{\bfseries}{:}{ }{}
\theoremstyle{italics}
\newtheorem{thm}[subsubsection]{Theorem}
\newtheorem{lem}[subsubsection]{Lemma}
\newtheorem{prop}[subsubsection]{Proposition}

\newtheoremstyle{noitalics}{}{}{}{}{\bfseries}{:}{ }{}
\theoremstyle{noitalics}
\newtheorem{rem}[subsubsection]{Remark}

\newtheorem{example}[subsubsection]{Example}
\newtheorem{defin}[subsubsection]{Definition}
\newtheorem{assump}[subsubsection]{Assumption}

\theoremstyle{italics}

\title{Reductions of Crystalline Representations for Small Weights}
\author{Anthony Guzman}
\address{Department of Mathematics, University of Arizona, Tucson, AZ 85721 USA}
\email{awguzman@arizona.edu}
\date{\today}

\begin{document}
	
\begin{abstract}
    This article is an extension of \cite{Guz1}. Using ideas found in loc. cit.,  we compute explicit reductions of crystalline representations of the absolute Galois group $\Gal(\overline{\Q}_p/\Qpf)$  with labeled Hodge-Tate weights in the range $p+2\le k_{0}\le 2p-4$ and $2\le k_i\le p-3$ for $1\le i\le f-1$. 
\end{abstract}
	
\maketitle

\setcounter{tocdepth}{1}
\tableofcontents

\section{Introduction}\label{intro-sec}

In \cite{Guz1} we computed semi-simple reductions of crystalline representations of the absolute Galois group $\Gal(\overline{\Q}_p/\Qpf)$ when the $p$-adic valuations of certain $a_2^{(i)}$-parameters were sufficiently large with respect to the labeled Hodge-Tate weights. Our strategy, surveyed in Section \ref{prelim-chap}, was to find simple representatives of the isomorphism classes of such representations and explicitly construct a Kisin module over an enlarged coefficient ring from the data of the associated weakly admissible filtered $\varphi$-module. Under the aforementioned restrictions, we were then able to descend the coefficients of these Kisin modules so as to allow us to compute reductions. We refer the reader to the introduction of \textit{loc. cit.} for the global considerations motivating such computations.

Where the aforementioned article restricted itself to the case of large valuations on certain parameters, the purpose of this article is to specialize these arguments so that we may loosen the restrictions on the $p$-adic valuations with the trade off being the enforcement of harsher restrictions on the labeled Hodge-Tate weights:
\begin{align*}
	p+2&\le k_{0}\le 2p-4& 2&\le k_{i}\le p-3 \hspace{.2cm}\text{\normalfont for all}\hspace{.2cm} 1\le i\le f-1.
\end{align*}

\noindent This quest involves much more delicate computations than those performed in \cite{Guz1} due to their `by hand' nature. Our methods are not foolproof however and indeed, some cases are out of reach using our strategy for technical reasons. When we are able to use our methods, we explicitly compute reductions in Chapter \ref{comp-red-chap} in the following cases:
\begin{itemize}
	\item $\nu_p(a_2^{(0)})\ge 1$;
	\item $\nu_p(a_2^{(0)})<1$ and \# Type $\I$ is even if we assume $\nu_p(a_2^{(i)})\ge 1-\nu_p(a_2^{(0)})$ for select $i\in\Z/f\Z$;
	\item $\nu_p(a_2^{(0)})<1$ and \# Type $\I$ is odd if we assume $\nu_p(a_2^{(i)})\ge 1$ for select $i\in\Z/f\Z$;
	\item In all cases for $f=2$.
\end{itemize}

%%%%%%%%%%%%%%%%%%%%%%%%%%%%%%%%%%%%%%%%%%%%%%%%%%%%%%%%%%%%%%%%%%%%%%%%%%%%%%%%%%%%%%%
\subsection{Notation}\label{notation-sec}
Let $p>2$ be an odd prime and equip $\Qp$ with a $p$-adic valuation $\nu_p$ normalized so that $\nu_p(p)=1$. Upon fixing an algebraic closure $\overline{\Q}_p$, we denote the $p$-adic completion of $\overline{\Q}_p$ by $\C_p$. 

Let $K$ be the unique unramified extension of $\Qp$ with (inertial) degree $f\ge1$; that is, we have an isomorphism $K\cong \Qpf$. Denote its ring of integers by $\oh_K$ with uniformizer $\pi_K$ and its residue field $k=\oh_K/\pi_K\oh_K\cong\Fpf$. Let $\sigma_K$ denote the absolute Frobenius on $K$ induced from the natural Frobenius on the residue field $k$. We define the absolute Galois group of $K$ to be $G_K=\Gal(\overline{K}/K)$ and we will let $I_K$ denote the inertia subgroup of $G_K$.

Fix a uniformizer $\pi_K=-p$ of $K$ so that the Eisenstein polynomial of $\pi_K$ in the formal variable $u$ may be denoted $E\coloneqq E(u)=u+p\in\oh_K[u]$. We fix once and for all a $p$-power compatible sequence $\pi=(\pi_0,\pi_1,\pi_2,\dots)$ in $\overline{K}$ such that $\pi_0=\pi_K$ and $\pi_n^p=\pi_{n-1}$. We define the field $K_\infty$ to be the compositum $K_\infty=\cup_nK(\pi_n)$ contained in $\overline{K}$ and the absolute Galois group of $K_\infty$ is denoted $G_\infty\coloneqq\Gal(\overline{K}/K_\infty)$. 

Let $\Lambda\subset K_0\pu$ denote the ring of rigid analytic functions on the open $p$-adic unit disc in $K$ and set $\Sig\coloneqq W(k)\pu\subset \Lambda$. The ring $K_0\pu$ admits a Frobenius action $\varphi$ which acts on coefficients by the absolute Frobenius $\sigma_K$ on $K$ and acts on the formal variable $u$ by $\varphi(u)=u^p$. Observe that the rings $\Sig\subset\Lambda\subset K_0\pu$ are $\varphi$-stable under this definition.

Acting as linear coefficients, let $F$ be a finite extension of $\Qp$ taken large enough to admit an embedding $\tau_0:K\hookrightarrow F$. Fix a uniformizer $\varpi$ of $F$ and define its ring of integers by $\oh_F$ with residue field $k_F\coloneqq\oh_F/\varpi\oh_F$. Extending the scalars of our rigid analytic functions, let $\Sig_F\coloneqq \Sig\otimes_{\Zp}\oh_F$ and $\Lambda_F\coloneqq \Lambda\otimes_{\Q_p}F$. Extending the $\varphi$-action on $K_0\pu$ by $F$-linearity, we obtain $\varphi$-stable rings $\Sig_F\subset\Lambda_F\subset (F\otimes_{\Q_p}K_0)\pu$. 

\subsection{Acknowledgments} This article alongside \cite{Guz1} were borne out of the authors Ph.D. thesis at the University of Arizona under the direction of Brandon Levin, whose generosity and knowledge made this work possible. It will quickly become evident the influence that Tong Liu and John Bergdall have in this paper, we thank them for the insightful conversations related to this project. 

%%%%%%%%%%%%%%%%%%%%%%%%%%%%%%%%%%%%%%%%%%%%%%%%%%%%%%%%%%%%%%%%%%%%%%%%%%%%%%%%%%%%%%%
%%%%%%%%%%%%%%%%%%%%%%%%%%%%%%%%%%%%%%%%%%%%%%%%%%%%%%%%%%%%%%%%%%%%%%%%%%%%%%%%%%%%%%%
\section{Preliminaries}\label{prelim-chap}

\subsection{Models for Crystalline Representations}\label{models-subsec}
For the readers convenience, let us recall the results we will require from \cite{Guz1} which are heavily inspired by \cite{Liu21, BL20, BLL22, Zhu09}. Due to the bounding of labeled Hodge-Tate weights $k_i<2p$ in this article, we will state specialized versions when appropriate.

Suppose $V$ is a two-dimensional crystalline $F$-representation of $G_K$ and denote the category of such objects by $\Rep_{\cris/F}(G_K)$. The associated weakly admissible filtered $\varphi$-module $D\coloneqq\D_\cris^*(V)$ is valued over the tensor $F\otimes_{\Q_p}K$ and in order to get a hold of its structure, we utilize the $f$-embeddings $\tau_i\coloneqq \tau_0\circ\sigma_K^{-i}:K\hookrightarrow F$ to decompose $D=\prod_i D^{(i)}$ into pieces, so we can write down explicit partial Frobenius matrices and filtration structures on each individual piece, see \cite[Prop. 2.1.2]{Guz1}. There exists a strongly divisible lattice $L$ over $\oh_F\otimes_{\Z_p}\oh_K$ inside $D$ by \cite[Prop. 3.1.6]{Guz1} and this strong divisibility structure will descend along the decomposition $D=\prod D^{(i)}$ so that we may write $L=\prod L^{(i)}$ where each $L^{(i)}$ may be viewed as a strongly divisible lattice inside $D^{(i)}$, see \cite[Prop. 3.1.2]{Guz1}.

We may choose a basis so that the filtration structure is fixed for our strongly divisible lattice $L$, and hence it will be fixed on $D$. In this basis, the structure of $L$ is completely determined by a $f$-tuple of matrices $(A)=(A^{(0)},A^{(1)},\dots, A^{(f-1)})\in\GL_2(\oh_F)^f$ in that we write $L=L(A)=\prod L(A)^{(i)}$ where the partial Frobenius matrix on $L(A)^{(i)}$ in this chosen basis is given by $\varphi^{(i)}=A^{(i)}\cdot\Diag(p^{k_{i-1}},1)$. The result is that the isomorphism classes of strongly divisible lattices are determined by the $f$-tuple $(A^{(i)})\in\GL_2(\oh_F)^f$. Moreover, we give a concrete method to describe such isomorphism classes in terms of parabolic equivalence classes on $\GL_2(\oh_F)^f$ via a bijection $\Theta:[(A)]\leftrightarrow L(A)$, see \cite[Thm. 3.2.2]{Guz1}. 

We then utilize a matrix simplifying algorithm inspired by to find a `nice' representative of $[(A)]$ in $\GL_2(\oh_F)^f$ so that each $A^{(i)}$ takes the form of one of two matrix \textit{Types}:
\begin{itemize}
	\item Type $\I$:	$A^{(i)}=\begin{pmatrix}0 & a_{1}^{(i)} \\ 1 & a_{2}^{(i)} \end{pmatrix}$ where $a_{1}^{(i)}\in\oh_F^*$ and $a_{2}^{(i)}\in\oh_F$.
	\item Type $\II$:	$A^{(i)}=\begin{pmatrix}a_{1}^{(i)} & 0 \\ a_{2}^{(i)} & 1\end{pmatrix}$ where $a_{1}^{(i)}\in\oh_F^*$ and $a_{2}^{(i)}\in\varpi\oh_F$.
\end{itemize}

\noindent Since the equivalence class $[(A)]$ determines an isomorphism class of strongly divisible lattices which in turn dictate the structure of weakly admissible filtered $\varphi$-modules, then $D(A)\coloneqq L(A)\otimes_{\oh_F}F$ will be a `simple' representative of its isomorphism class. We write $V(A)$ to be the crystalline representation whose image under $\D_\cris^*$ is $D(A)$. 

\begin{thm}\label{model-irred-cris-reps-thm}
	Let $V\in\Rep_{\cris/F}(G_K)$ be an irreducible, two dimensional crystalline representation of $G_K$ with labeled Hodge-Tate weights $\{k_i,0\}_{i\in\Z/f\Z}$ where $k_i>0$. Then $V\cong V(A)$ where $D(A)\coloneqq\D^*_\cris(V(A))$ with each $A^{(i)}$ being Type $\I$ or $\II$. In particular, $D(A)=\prod_iD(A)^{(i)}$ has a basis $\{\eta_1^{(i)},\eta_2^{(i)}\}$ such that
	\begin{align*}
		[\varphi^{(i)}]_\eta &=\begin{dcases}
			\begin{pmatrix}
				0 & a_1^{(i)} \\ p^{k_{i-1}} & a_2^{(i)}
			\end{pmatrix} & \text{\normalfont Type $\I$} \\
			\begin{pmatrix}
				a_1^{(i)}p^{k_{i-1}} & 0 \\ a_2^{(i)}p^{k_{i-1}} & 1
			\end{pmatrix} & \text{\normalfont Type $\II$}
		\end{dcases} & \Fil^j D(A)^{(i)}&=\begin{dcases}
			D(A)^{(i)} & j\le 0 \\
			F(\eta_1^{(i)}) &  0<j\le k_i \\
			0 & k_i<j.
		\end{dcases}
	\end{align*}
	\begin{proof}
		See \cite[Thm. 4.2.3]{Guz1}
	\end{proof}
\end{thm}

Hence, to compute reductions, we need only compute reductions of those representations of the form $V(A)$. Moreover, we may use irreducibility to enforce bounds on the $p$-adic valuations of the $a_2^{(i)}$-parameters.  We identify two disjoint subsets
\begin{align*}
	\Ss&\coloneqq\{i\in\Z/f\Z: \As^{(i)} \hspace{.1cm}\text{\normalfont is Type I}\hspace{.1cm}\}&\Ts&\coloneqq\{i\in\Z/f\Z: \As^{(i)} \hspace{.1cm}\text{\normalfont is Type II}\}.
\end{align*}
so that $\Z/f\Z=\Ss\sqcup \Ts$.

\begin{prop}\label{irred-vals-prop}
	Let $V(A)$ be as defined in Theorem \ref{model-irred-cris-reps-thm}. Then irreducibility implies that:
	\begin{enumerate}
		\item $a_2^{(i)}\in\varpi\oh_F$ for all $\Ts$;
		\item $\prod_{i\in\Ss}a_2^{(i)}\in\varpi\oh_F$.
	\end{enumerate}
	\begin{proof}
		Point $(a)$ follows from the definition of Type $\II$. For $(b)$, see \cite[Prop. 4.2.1]{Guz1}.
	\end{proof}
\end{prop}

\subsection{Constructing and Descending Kisin Modules}\label{const-desc-subsec}
With an understanding of the structure of $V(A)$ in place, we then move on to the task of computing the semi-simple modulo $\varpi$ reductions of irreducible, two-dimensional crystalline representations $V(A)$. By this, we mean the computation of semi-simple $k_F$-representations $\overline{V(A)}=(T/\varpi T)^{\semi}$ where $T\subset V(A)$ is a $G_K$-stable $\oh_F$-lattice inside of $V(A)$. To do so, we employ the notation of a Kisin module.

\begin{defin}\label{kisin-mod-def}
	A \textit{Kisin module over $\Sig_F$ of ($E$-)height $\le h$} is a $\varphi$-module $\Mf$ over $\Sig_F$ such that the linearization map
	\[1\otimes\varphi_\Mf:\varphi^*\Mf\coloneqq \Sig_F\otimes_{\varphi,\Sig_F}\Mf\rightarrow\Mf\]
	has cokernel killed by $E^h$. Denote the category of such objects by $\Mod_{\Sig_F}^{\varphi,\le h}$.
\end{defin}

The motivation for doing so is that a Kisin module $\Mf$ over $\Sig_F$ may be used to compute $\overline{V(A)}$ via the data provided by the modulo $\varpi$ reduction $\overline{\Mf}$ \textit{if} the Kisin module $\Mf$ is canonically associated to a Galois-stable $\oh_F$-lattice $T\subset V(A)$. Let $\V_{k_F}:\Mod_{k_F(\!(u)\!)}^{\varphi,\et}\rightarrow\Rep_{/k_F}(G_\infty)$ denote Fountain's \'{e}tale $\varphi$-module functor, see \cite[\S~A1]{Fon91}.

\begin{prop}\label{kisin-red-prop}
	Let $\Mf\in\Mod_{\Sig_F}^{\varphi,\le h}$ and $V(A)\in\Rep_{\cris/F}(G_K)$ be such that $\Mf=\Mf(T)$ is the canonically associated Kisin module to a $G_\infty$-stable $\oh_F$-lattice $T\subset V(A)$. Then \[(\V_{k_F}^*(\Mf\otimes_{\oh_F}k_F[1/u]))^\semi\cong \overline{V(A)}|_{G_\infty}.\]
	\begin{proof}
		See \cite[Prop 2.2.9]{Guz1}.
	\end{proof}
\end{prop}

The issue that persists is that it is very difficult to determine the structure of $\Mf$ from $V(A)$. To get around this, we detail an algorithm to explicitly construct a Kisin module $\Ms(\As)$ over an enlarged coefficient ring $S_F= \Sig_F[\![E^p/p]\!]$ from the data of the weakly admissible filtered $\varphi$-module $D(A)$. Define $\gamma=\varphi(E)/p\in S_F^*$ and	\[\lambda_b=\prod_{n\ge0}\varphi^{bn}(\gamma)\in S_F^*.\]

\begin{prop}\label{const-kisin-prop}
	There exists a finite height Kisin module $\Ms(\As)=\prod_i\Ms(\As)^{(i)}$ admitting a $S_F$-basis $\{\eta_1^{(i)},\eta_2^{(i)}\}$ such that $\varphi_{\Ms(A)}^{(i)}(\eta_1^{(i-1)},\eta_2^{(i-1)})=(\eta_1^{(i)},\eta_2^{(i)})\As^{(i)}$ where
	\begin{align*}
		\As^{(i)}=\begin{dcases}
			\begin{pmatrix}
				0 & E^{k_i}a_1^{(i)} \\ 1 & a_2^{(i)}\lambda_b^{h^{(i)}(\varphi)-g^{(i)}(\varphi)}
			\end{pmatrix} & \textit{if	} i\in\Ss \\
			\begin{pmatrix}
				E^{k_i}a_1^{(i)} & 0 \\ a_2^{(i)}\lambda_b^{h^{(i)}(\varphi)-g^{(i)}(\varphi)} & 1
			\end{pmatrix} & \textit{if	} i\in\Ts
		\end{dcases}
	\end{align*}
	\begin{proof}
		The existence of $\Ms(\As)$ follows from \cite[Prop 5.1.2]{Guz1}. The particular basis for which $\As^{(i)}$ is of this form is by \cite[Prop 5.2.2]{Guz1}.
	\end{proof}
\end{prop}

To compute reductions in the sense of Proposition \ref{kisin-red-prop}, we must find a descent $\Mf(\Af)$ of $\Ms(\As)$ to $\Sig_F$ in the sense that
\[\Mf(\Af)\otimes_{\Sig_F}S_F=\Ms(\As).\]
Moreover, we must canonically associated the descent $\Mf(\Af)$ to some $G_\infty$-stable $\oh_F$-lattice $T\subset V(A)$. All of this is accomplished by our \textit{descent algorithm}.

\begin{rem}
	In the following, we specialize results of \cite[\S 6]{Guz1} to the case of $c_{max}=2$.
\end{rem}

\noindent Define an ideal $I=\varpi p S_F+\varpi\Fil^{2p}S_F+E\Fil^{2p}S_F$ of $S_F$. We write $X^{(i)}*_\varphi A^{(i)}=X^{(i)}A^{(i)}\varphi(A^{(i-1)})$. Let us establish a set of assumptions under which descent may take place.

\begin{assump}[Descent Assumptions]\label{descent-assump}
	Let $\Ms(\As)$ be as in Proposition \ref{const-desc-subsec}. Suppose there exists a sequence of base changes $(X_n)$ so that by setting $(X_n^{(i)})*_\varphi(\As^{(i)})=(\As_n^{(i)})$, the following holds:
	\begin{enumerate}
		\item We have $k_i\le 2p-4$ for all $i\in\Z/f\Z$;
		\item Either $X_n^{(i)}\in\GL_2(S_F[1/p])$ with $\det(X_n^{(i)})=1$ or $X_n^{(i)}=\Diag(x,y)\in\GL_2(F)$;
		\item There exists a finite $m>0$ such that 
		\[\As_m^{(i)}=\Af_0^{(i)}+C^{(i)}\]
		where $\Af_0^{(i)}\in\Mat_2(\Sig_F)$ and $C^{(i)}\in\Mat_2(I)$.
	\end{enumerate}
\end{assump}

\begin{thm}\label{desc-alg-thm}
	Suppose $p>2$. If $\Ms(\As)$ satisfies the Descent Assumptions \ref{descent-assump}, then there exists a descent $\Mf(\Af)$ over $\Sig_F$ of $\Ms(\As)$ which is canonically associated to a $G_\infty$-stable $\oh_F$-lattice $T\subset V(A)$ such that $\Af^{(i)}\equiv\Af_0^{(i)}\pmod\varpi$.
	\begin{proof}
		For existence of $\Mf(\Af)$, see \cite[Thm. 6.2.4]{Guz1}. The canonical association to a Galois stable lattice $T\subset V(A)$ is via \cite[Prop 7.1.1]{Guz1}.
	\end{proof}
\end{thm}

\subsection{Computing Reductions}\label{comp-red-subsec}

Let $\Mf(\Af)$ be the descent of $\Ms(\As)$ canonically associated to a $G_\infty$-stable $\oh_F$-lattice $T\subset V(A)$ as given in Theorem \ref{desc-alg-thm}. By setting $\overline{\Mf(\Af)}\coloneqq\Mf(\Af)/\varpi\Mf(\Af)$, Proposition \ref{kisin-red-prop} says that the reduction $\overline{V(A)}\vert_{G_\infty}=(T/\varpi T)^\semi\vert_{G_\infty}$ is completely determined by the $k_F$-representation of $G_\infty$ given by
\[\overline{V(A)}|_{G_\infty}\cong \V_{k_F}^*\left(\overline{\Mf(\Af)}\left[\frac{1}{u}\right]\right).\]

To compute the value of this functor, we need to reduce our $f$-many partial Frobenius matrices to a single one. We may write $\overline{\Mf(\Af)}=\prod_{i}\overline{\Mf(\Af)}^{(i)}$ so that each $\overline{\Mf(\Af)}^{(i)}$ is a $\varphi$-module with Frobenius action given by $\overline{\Af}^{(i)}\coloneqq \Af^{(i)}\pmod\varpi\equiv\Af_0^{(i)}\pmod\varpi$. We endow the $k_F[\![u]\!]$-module $\overline{\Mf(\Af)}^{(0)}$ with the structure of an \'{e}tale $\varphi^f$-module by defining its Frobenius action to be given by the product
\[\prod_{i\in\Z/f\Z}\Af^{(i)}=\overline{\Af}^{(0)}\varphi\left(\overline{\Af}^{(1)}\right)\varphi^2\left(\overline{\Af}^{(2)}\right)\cdots\varphi^{f-1}\left(\overline{\Af}^{(f-1)}\right).\] 
Such objects are defined just as in \cite[\S 2]{CDM14}.
\begin{lem}\label{red-in-first-piece-lem}
	Using the notation from above, let $(M^{(0)},\varphi^f)$ be the \'{e}tale $\varphi^f$-module over $k_F(\!(u)\!)$ associated to an \'{e}tale $\varphi$-module $(M,\varphi)$ over $k_F(\!(u)\!)$. 
	Then
	\[\V_{k_F}^*\left(M\right)\cong\Hom_{k_F(\!(u)\!),\varphi^f}\left(M^{(0)},\oh_{\C_p}^\flat\right).\]
	\begin{proof}
		See \cite[Lem. 7.1.2]{Guz1}
	\end{proof}
\end{lem}

Recall that $I_K$ denotes the inertial subgroup of $G_K$ and choose for any $r\ge1$ a $p^r-1$-th root of $\pi_K=-p$ denoted $\xi_r=(-p)^{1/p^r-1}\in \overline{K}$. This choice of $\xi_r$ defines a character $\omega_{\xi_r}:I_K\rightarrow \oh_K^*$. By composing with our fixed embedding $\tau_0:\oh_K\hookrightarrow\oh_F$, we give rise to a \textit{fundamental character of niveau $r$}
\[\widetilde{\omega}_{\xi_r}:I_K\rightarrow\oh_F^*.\]
The residual character of $\widetilde{\omega}_{\xi_r}$, denoted $\omega_r$, is then a \textit{fundamental character of level $r$}.

\begin{prop}\label{compute-Hom-prop}
	Let $\overline{\Mf}=\prod_{i}\overline{\Mf}^{(i)}$ be a finite height Kisin module over $k_F[\![u]\!]$ with basis $\{e_1^{(i)},\dots,e_d^{(i)}\}_{i\in\Z/f\Z}$. Suppose that $\varphi^f(e_j^{(0)})=u^{a_j}e_{j+1}^{(0)}$ with $e_{d+1}^{(0)}=e_1^{(0)}$ and $0\le a_j\le p^f-1$. Then 
	\[\V_{k_F}^*(\overline{\Mf})|_{G_\infty}=\ind_{G_{\Q_{p^{df}}}}^{G_{\Q_{p^f}}}\omega_{df}^{\sum_{j=1}^{d}p^{d-j}a_j}.\]
	\begin{proof}	
		See \cite[Prop. 7.1.3]{Guz1}. Here we may enforce $a_j\le p^f-1$ due to the fact that $k_i<2p$.
	\end{proof}
\end{prop}

With this, we have our primary method of computing the $k_F$-representation associated to a finite height Kisin module. However, we will soon see that the condition that $\varphi^f(e_j^{(0)})=u^{a_j}e_{j+1}^{(0)}$ is often not satisfied out of the box. In order to find the basis where this behavior holds, we rely on a new lemma based on \cite[Lemma 5.2]{Liu21}. The idea is that if a row operation is equivalent to identity modulo $u^2$, then we may `straightened' the $\varphi$-twisted conjugation to simply be the single row operation.

\begin{lem}\label{straightening-lem}
	For $i\in\Z/f\Z$, suppose $\Af^{(i)}\in\Mat_d^{h^{(i)}}(k_F[\![u]\!])$ for $h^{(i)}\le 2p-3$. If $(X^{(i)})\in \GL_d(k_F[\![u]\!])^f$ is such that $X^{(i)}\equiv I_d\pmod{u^2}$ for all $i$, then there exists $(Y^{(i)})\in\GL_d(k_F[\![u]\!])$ such that
	\[Y^{(i)}*_\varphi\Af^{(i)}=X^{(i)}\Af^{(i)}.\]
	
	\begin{proof}
		Since $X^{(i)}\equiv I_d\pmod{u^2}$ for all $i\in\Z/f\Z$, then 
		\[X^{(i)}*_\varphi\Af^{(i)}=X^{(i)}\Af^{(i)}\varphi(X^{(i-1)})^{-1}=X^{(i)}\Af^{(i)}+u^{2p}Z^{(i)}\]
		for some $Z^{(i)}\in\Mat_d(k_F[\![u]\!])$. Since each $\Af^{(i)}\in\Mat_d^{h^{(i)}}(k_F[\![u]\!])$ then there exists a matrix $B^{(i)}\in \Mat_d^{h^{(i)}}(k_F[\![u]\!])$ such that $\Af^{(i)}B^{(i)}=B^{(i)}\Af^{(i)}=u^{h^{(i)}}I_d$. Hence, we may write 
		\begin{align*}
			X^{(i)}*_\varphi\Af^{(i)}&=X^{(i)}\Af^{(i)}+u^{2p-h^{(i)}}Z^{(i)}B^{(i)}(X^{(i)})^{-1}X^{(i)}\Af^{(i)} \\
			&=(I_d+u^{h_1^{(i)}}Z_1^{(i)})X^{(i)}\Af^{(i)}
		\end{align*}
		
		\noindent where $h_1^{(i)}=2p-h^{(i)}\ge 3$ and $Z_1^{(i)}=Z^{(i)}B^{(i)}(X^{(i)})^{-1}$. We set $Y_1^{(i)}=(I_d+u^{h_1^{(i)}}Z_1^{(i)})^{-1}\in\GL_d(k_F[\![u]\!])$ so that by the same reasoning as the above arguments, we have
		\begin{align*}
			(Y_1X)^{(i)}*_\varphi\Af^{(i)}&=Y_1^{(i)}X^{(i)}\Af^{(i)}\varphi\left((X^{(i-1)})^{-1}(Y_1^{(i-1)})^{-1}\right) \\
			&=Y_1^{(i)}(Y_1^{(i)})^{-1}X^{(i)}\Af^{(i)}\varphi(Y_1^{(i-1)})^{-1} \\
			&=X^{(i)}\Af^{(i)}+u^{ph_1^{(i-1)}}Z'^{(i)} \\
			&=X^{(i)}\Af^{(i)}+u^{ph_1^{(i-1)}-h^{(i)}}Z'^{(i)}B^{(i)}(X^{(i)})^{-1}X^{(i)}\Af^{(i)} \\
			&=(I_d+u^{h_2^{(i)}}Z_1^{(i)})X^{(i)}\Af^{(i)}
		\end{align*}
		
		\noindent where $h_2^{(i)}=ph_1^{(i-1)}-h^{(i)}\ge p+3$ and $Z_2^{(i)}=Z'^{(i)}B^{(i)}(X^{(i)})^{-1}$. Repeating this process will give us a base change
		\[(Y_nY_{n-1}\cdots Y_1X)^{(i)}*_\varphi\Af^{(i)}=(I_d+u^{h_n^{(i)}}Z_n^{(i)})X^{(i)}\Af^{(i)}\]
		with $h_n^{(i)}=ph_{n-1}^{(i-1)}-h^{(i)}\ge 3+\sum_{j=1}^np^j$. Hence we see that $h_n^{(i)}\rightarrow \infty$ and by setting $Y^{(i)}=\left(\prod_{n=1}^{\infty}Y_n\right)^{(i)}X^{(i)}\in\GL_d(k_F[\![u]\!])$ then our result follows.
	\end{proof}
\end{lem}

We now need a way to encode information from $\overline{\Mf(\Af)}$ that is necessary to compute $\overline{V(A)}$. Let $\ident$ denote the identity matrix and let $S=\begin{psmallmatrix}
	0 & 1 \\ 1 & 0
\end{psmallmatrix}$. For a pair $\lambda_i=(n_i,m_i)\in\Z^2$ with each $0\le n_i,m_i\le p-1$, we will set $u^{\lambda_i}=\Diag(u^{n_i},u^{m_i})$. In the coming sections, we will describe a $k_F[\![u]\!]$-basis such that each $\overline{\Af}^{(i)}$ is either $\ident u^{\lambda_i}$ or $Su^{\lambda_i}$ up to scalars in $k_F$. We will then summarize this data in the form $\mu_i=(M_i,\lambda_i)$ with $M_i\in\{\ident,S\}$ and $\lambda_i$ is as defined above. Let us also decompose $i\in\Z/f\Z$ into two disjoint subsets $\Z/f\Z=\Is\sqcup\Ps$ where $\Is$ denotes the set of $i$ such that $A_i=\ident$ and $\Ps$ denotes the set of $i$ with $A_i=S$.

\begin{defin}\label{red-data-def}
	We will call the $f$-tuple of pairs $\mu=(\mu_i)=(M_i,\lambda_i)_i$ derived from $\overline{\Mf(\Af)}$ to be its \textit{reduction data}.
\end{defin}

\noindent As the name suggests, this data suffices to compute the reduction $\overline{V(A)}|_{G_\infty}$ explicitly given specific data, up to an unramified twist in the irreducible case or restriction to inertia in the reducible case.

\begin{prop}\label{compute-hom-with-data-prop}
	Using the notation from above, let $\overline{\Mf(\Af)}$ be a finite height, rank-two Kisin module with reduction data $\mu=(\mu_i)$ canonically associated to $T\subset V(A)$. Then there exists $v_i,w_i$ with $0\le v_i,w_i\le p-1$ such that one of two cases hold:
	\begin{enumerate}
		\item If the order of the set $|\Ps|$ for $\mu$ is even or zero, then we have reducible reduction
		\[\overline{V(A)}|_{I_K}=\omega_f^{\sum_{j=0}^{f-1} p^{j}v_j}\oplus\omega_f^{\sum_{j=0}^{f-1} p^{j}w_j};\]
		\item If the order of the set $|\Ps|$ is odd, then we have an irreducible reduction up to unramified twist
		\[\overline{V(A)}|_{G_\infty}=\ind^{G_{\Qpf}}_{G_{\Q_{p^{2f}}}}\omega_{2f}^{\sum_{j=0}^{f-1} p^{j+1}w_j+\sum p^{j}v_j}.\]
	\end{enumerate}
	\begin{proof}
		See \cite[Prop. 7.1.5]{Guz1} with the reducible case in $(b)$ of \textit{loc. cit.} being absent due to restricting $v_i,w_i\le p-1$.
	\end{proof}
\end{prop}
	
%%%%%%%%%%%%%%%%%%%%%%%%%%%%%%%%%%%%%%%%%%%%%%%%%%%%%%%%%%%%%%%%%%%%%%%%%%%%%%%%%%%%%%%
%%%%%%%%%%%%%%%%%%%%%%%%%%%%%%%%%%%%%%%%%%%%%%%%%%%%%%%%%%%%%%%%%%%%%%%%%%%%%%%%%%%%%%%
\section{Preparatory Computations for Descent}\label{prep-comp-chap}

As an implementation of the Descent Algorithm \ref{desc-alg-thm}, we would like to display the computations necessary to show that the Descent Assumptions \ref{descent-assump} hold for $\Ms(\As)$ with labeled heights in the range
\begin{align*}
	p+2&\le k_{0}\le 2p-4& 2&\le k_{i}\le p-3 \hspace{.2cm}\text{\normalfont for all}\hspace{.2cm} 1\le i\le f-1
\end{align*}

\noindent and the $p$-adic valuations $\nu_p(a_2^{(i)})$ are sufficiently large for select $i\in\Z/f\Z$. Moreover, we will provide a complete description in the case of $f=2$. Our strategy will be to perform a finite number of successive row and scaling operations $(X_n)*_\varphi(\As)=(\As_n)$ until each entry of $\As_n^{(i)}$ is in $\Sig_F+I$ at which point the Descent Assumptions \ref{descent-assump} will be satisfied.

%%%%%%%%%%%%%%%%%%%%%%%%%%%%%%%%%%%%%%%%%%%%%%%%%%%%%%%%%%%%%%%%%%%%%%%%%%%%%%%%%%%%%%%
\subsection{The Large Valuation Case}\label{prep-large-val-sec}

We begin with the case of $\nu_p(a_2^{(0)})\ge 1$ which we collectively call the \textit{large valuation} cases, following the ideas found in \cite[\S 6.3]{Guz1}. Note that we will not make this assumption quite yet as this first step is valid for any valuation and will form the first step in all proceeding arguments cases. 

Since $\lambda_b\in S_F^*$, let us write $\lambda_b^{g^{(i)}(\varphi)}=\sum_{j=0}^{\infty}\alpha^{(i)}_j (E^p/p)^j$ so we may define elements in $S_F[1/p]$ given by
\[x^{(i)}=\begin{dcases}
	-\left(\frac{a_2}{a_1}\right)^{(i)}\sum_{j=2}^{\infty}\alpha^{(i)}_j \frac{E^{jp-k_i}}{p^j} & \text{if	}i= 0\\
	-\left(\frac{a_2}{a_1}\right)^{(i)}\sum_{j=1}^{\infty}\alpha^{(i)}_j \frac{E^{jp-k_i}}{p^j} & \text{for	}1\le i\le f-1.
\end{dcases}\]

Set $X_1^{(i)}=\begin{psmallmatrix}
	1 & 0 \\ x^{(i)} & 1
\end{psmallmatrix}\in\GL_2(S_F[1/p])$ for all $i\in\Z/f\Z$ and observe that by setting $(X_1)*_\varphi(\As)=(\As_1)$ we will obtain
\begin{align*}
		\As_1^{(0)}&=\begin{dcases}
		\begin{pmatrix}
			E^{k_{0}}a_1^{(0)}\varphi(x^{(f-1)}) & E^{k_{0}}a_1^{(0)} \\ 1+((a_2\alpha_0)^{(0)}+(a_2\alpha_1)^{(0)}\frac{E^p}{p})\varphi(x^{(f-1)}) & (a_2\alpha_0)^{(0)}+(a_2\alpha_1)^{(0)}\frac{E^p}{p}
		\end{pmatrix} & \text{if	} 0\in\Ss\\
		\begin{pmatrix}
			E^{k_{0}}a_1^{(0)} & 0 \\ (a_2\alpha_0)^{(0)}+(a_2\alpha_1)^{(0)}\frac{E^p}{p}+\varphi(x^{(f-1)}) & 1
		\end{pmatrix}& \text{if	} 0\in\Ts 
	\end{dcases} \\
	\As_1^{(i)}&=\begin{dcases}
		\begin{pmatrix}
			E^{k_i}a_1^{(i)}\varphi(x^{(i-1)}) & E^{k_i}a_1^{(i)} \\ 1+(a_2\alpha_0)^{(i)}\varphi(x^{(i-1)}) & (a_2\alpha_0)^{(i)}
		\end{pmatrix} & \text{for	} i\neq 0 \text{ and } i \in\Ss\\
		\begin{pmatrix}
			E^{k_i}a_1^{(i)} & 0 \\ (a_2\alpha_0)^{(i)}+\varphi(x^{(i-1)}) & 1
		\end{pmatrix}& \text{for	} i\neq 0 \text{ and } i \in\Ts. \\
	\end{dcases} 
\end{align*}

\noindent In order to ensure that none of the $\varphi(x^{(i)})$-terms introduced as a result of $\varphi$-conjugation are non-integral or non-$I$ factors, we need a slightly modified version of \cite[Lem. 6.3.2]{Guz1}.

\begin{lem}\label{phi-admis-lem}
	If $i\neq 0$ and $k_i\le p-h$ for $h>0$, then $\varphi(x^{(i)})\in a_2^{(i)}p^{h-1}S_F$. Likewise, if $k_{0}\le 2p-4$ then $\varphi(x^{(0)})\in a_2^{(0)}p^2S_F$.
	\begin{proof}
		Follows by the same calculations present in \cite[Lem. 6.3.2]{Guz1}.
	\end{proof}
\end{lem}

\begin{rem}\label{bound-rem}
	The bound $k_i\le p-3$ for $1\le i\le f-1$ differs from those found in \cite{Guz1} which depend on multiples of $p-2$. The reason for this is that we can no longer guarantee that $\nu_p(a_2^{(i)})>0$ for all $i\in\Z/f\Z$, see Proposition \ref{irred-vals-prop} and the below discussion.
\end{rem}

Using the fact that $a_2^{(i)}p\in \varpi\oh_F$ for all $i\in\Z/f\Z$, we can see that by Lemma \ref{phi-admis-lem}, when $\nu_p(a_2^{(i)})>0$ we can write $\varphi(x^{(i)})\in p^2\varpi S_F\subset pI$ and whenever $\nu_p(a_2^{(i)})=0$, we have $\varphi(x^{(i)})\in p\varpi S_F\subset I$. To reduce notation, let us define a nonzero valuation indicator function in the following way
\[v(i)=\begin{dcases}
	1 & \text{if	} \nu_p(a_2^{(i)})>0 \\
	0 & \text{if	} \nu_p(a_2^{(i)})=0.
\end{dcases}\]
We then have that $\varphi(x^{(i)})\in \varpi^{v(i)} p^2 S_F\subset p^{v(i)}I$. Let us also begin to gather elements of the ideal $I$ into terms denoted $\varepsilon^{(i)}_{jk}$ and $\delta^{(i)}_{jk}$. Since Theorem \ref{desc-alg-thm} allows us to disregard such elements when computing reductions, then there is not harm in making this simplifying notation. By the above discussion, we write
\begin{align*}
	\As_1^{(0)}&=\begin{dcases}
		\begin{pmatrix}
			\varepsilon_{11}^{(0)} & E^{k_{0}}a_1^{(0)} \\ 1+\varepsilon_{21}^{(0)} & (a_2\alpha_0)^{(0)}+(a_2\alpha_1)^{(0)}\frac{E^p}{p}
		\end{pmatrix} & \text{if	} 0\in\Ss\\
		\begin{pmatrix}
			E^{k_{0}}a_1^{(0)} & 0 \\ (a_2\alpha_0)^{(0)}+(a_2\alpha_1)^{(0)}\frac{E^p}{p}+\delta_{21}^{(0)} & 1
		\end{pmatrix}& \text{if	} 0\in\Ts 
	\end{dcases} \\
	\As_1^{(i)}&=\begin{dcases}
	\begin{pmatrix}
		\varepsilon^{(i)}_{11} & E^{k_i}a_1^{(i)} \\ 1+\varepsilon^{(i)}_{21} & (a_2\alpha_0)^{(i)}
	\end{pmatrix} & \text{if	} i\neq 0 \text{ and } i \in\Ss\\
	\begin{pmatrix}
		E^{k_i}a_1^{(i)} & 0 \\ (a_2\alpha_0)^{(i)}+\delta^{(i)}_{21} & 1
	\end{pmatrix}& \text{if	} i\neq 0 \text{ and } i \in\Ts
\end{dcases} \end{align*}

\begin{itemize}
	\item $\varepsilon^{(i)}_{11}=E^{k_i}a_1^{(i)}\varphi(x^{(i-1)})\in E^{k_i}p^{v(i-1)}I$;
	
	\item $\varepsilon^{(i)}_{21}=(a_2\alpha_0)^{(i)}\varphi(x^{(i-1)})\in \varpi^{v(i)+v(i-1)}p^2S_F\subset I$;
	
	\item $\delta^{(i)}_{21}=\varphi(x^{(i-1)})\in\varpi^{v(i-1)}p^2S_F\subset p^{v(i-1)}I$;
	
	\item $\varepsilon_{21}^{(0)}=((a_2\alpha_0)^{(0)}+(a_2\alpha_1)^{(0)}\frac{E^p}{p})\varphi(x^{(f-1)})\in \varpi^{v(0)+v(f-1)}p^2S_F\subset I$.
\end{itemize}

\noindent With this, we observe  that the only non-integral and non-$I$ term is $(a_2\alpha_1)^{(0)}\frac{E^p}{p}$ in  $\As_1^{(0)}$. If we assume that $\nu_p(a_2^{(0)})\ge 1$ then this term becomes integral and we get our first case of satisfied descent.

\begin{prop}\label{dec-assump-large-val-prop}
	Suppose $\Ms(\As)$ is a rank-two Kisin module over $S_F$ with labeled heights $p+2\le k_{0}\le 2p-4$ and $2\le k_i\le p-3$ for $i\neq 0$. If $\nu_p(a_2^{(0)})\ge 1$, then the Frobenius matrix $(\As)$ satisfies the Descent Assumptions \ref{descent-assump} with
	\begin{align*}
				\Af_0^{(0)}&=\begin{dcases}
			\begin{pmatrix}
				0 & E^{k_{0}}a_1^{(0)} \\ 1 & (a_2\alpha_0)^{(0)}+\left(\frac{a_2^{(0)}}{p}\right)\alpha_1^{(0)}E^p
			\end{pmatrix} & \text{if	} 0\in\Ss\\
			\begin{pmatrix}
				E^{k_{0}}a_1^{(0)} & 0 \\ (a_2\alpha_0)^{(0)}+\left(\frac{a_2^{(0)}}{p}\right)\alpha_1^{(0)}E^p & 1
			\end{pmatrix}& \text{if	} 0\in\Ts
		\end{dcases}\\
		\Af_0^{(i)}&=\begin{dcases}
			\begin{pmatrix}
				0 & E^{k_i}a_1^{(i)} \\ 1 & (a_2\alpha_0)^{(i)}
			\end{pmatrix} & \text{if	} i\neq 0 \text{ and } i \in\Ss\\
			\begin{pmatrix}
				E^{k_i}a_1^{(i)} & 0 \\ (a_2\alpha_0)^{(i)} & 1
			\end{pmatrix}& \text{if	} i\neq 0 \text{ and } i \in\Ts.
		\end{dcases} 
	\end{align*}
\end{prop}

%%%%%%%%%%%%%%%%%%%%%%%%%%%%%%%%%%%%%%%%%%%%%%%%%%%%%%%%%%%%%%%%%%%%%%%%%%%%%%%%%%%%%%%
\subsection{The Small Valuation Case}\label{prep-small-val-sec}

We now assume that $0\le \nu_p(a_2^{(0)})<1$ and will continue from $\As_1^{(i)}$. Recall that at this point we have 
\begin{align*}
	\As_1^{(0)}&=\begin{dcases}
		\begin{pmatrix}
			\varepsilon_{11}^{(0)} & E^{k_{0}}a_1^{(0)} \\ 1+\varepsilon_{21}^{(0)} & (a_2\alpha_0)^{(0)}+(a_2\alpha_1)^{(0)}\frac{E^p}{p}
		\end{pmatrix} & \text{if	} 0\in\Ss\\
		\begin{pmatrix}
			E^{k_{0}}a_1^{(0)} & 0 \\ (a_2\alpha_0)^{(0)}+(a_2\alpha_1)^{(0)}\frac{E^p}{p}+\delta_{21}^{(0)} & 1
		\end{pmatrix}& \text{if	} 0\in\Ts 
	\end{dcases} \\
	\As_1^{(i)}&=\begin{dcases}
		\begin{pmatrix}
			\varepsilon^{(i)}_{11} & E^{k_i}a_1^{(i)} \\ 1+\varepsilon^{(i)}_{21} & (a_2\alpha_0)^{(i)}
		\end{pmatrix} & \text{if	} i\neq 0 \text{ and } i \in\Ss\\
		\begin{pmatrix}
			E^{k_i}a_1^{(i)} & 0 \\ (a_2\alpha_0)^{(i)}+\delta^{(i)}_{21} & 1
		\end{pmatrix}& \text{if	} i\neq 0 \text{ and } i \in\Ts
\end{dcases} 
\end{align*}

\noindent where $\varepsilon_{11}^{(i)}\in E^{k_i}p^{v(i-1)}I$, $\varepsilon_{21}^{(i)}\in \varpi^{(v(i)+v(i-1))}p^2S_F\subset I$ and $\delta_{21}^{(i)}\in p^{v(i-1)}I$. Our strategy will be to kill the $E^{k_{0}}$-term followed by a scaling operation and concluding (if necessary) with a row operation designed to kill the $E^p/p$ term which, as we may recall, is the only non-integral and non-$I$ term in the matrices above. The real danger in this process is the aforementioned scaling operation as we will have to be very careful to not introduce non-integral factors which will require careful estimates of certain terms.

Notice that since $0\le \nu_p(a_2^{(0)})<1$, we can write $\nu_p(p/a_2^{(0)})>0$. Consider the row operation over $\GL_2(\Sig_F)$ designed to kill the $E^{k_{0}}$-term in $\As_1^{(0)}$ given by
\begin{align*}
	X_2^{(0)}&=\begin{pmatrix}
		1 & -\frac{a_1^{(0)}p}{(a_2\alpha_1)^{(0)}}E^{k_{0}-p} \\  0 & 1
	\end{pmatrix}& X_2^{(i)}=\begin{pmatrix}
		1 & 0 \\ 0 & 1
	\end{pmatrix}\text{ for }i\neq 0.
\end{align*}

\noindent We see that this base change is isolated as to only affect $\As_1^{(0)}$ and $\As_1^{(1)}$ so the result of $(X_2)*_\varphi(\As_1)$ will be given by $\As_2^{(i)}=\As_1^{(i)}$ for $2\le i\le f-1$ and
\begin{align*}
	\As_2^{(0)}&=\begin{dcases}
		\begin{pmatrix}
			-\left(\frac{a_1p}{a_2\alpha_1}\right)^{(0)}E^{k_{0}-p}+\tilde{\varepsilon}_{11}^{(0)} & -\left(\frac{a_1\alpha_0p}{\alpha_1}\right)^{(0)}E^{k_{0}-p} \\ 1+\varepsilon_{21}^{(0)} & (a_2\alpha_0)^{(0)}+(a_2\alpha_1)^{(0)}\frac{E^p}{p}
		\end{pmatrix} & \text{if	} 0\in\Ss\\
		\begin{pmatrix}
			-\left(\frac{a_1\alpha_0p}{\alpha_1}\right)^{(0)}E^{k_{0}-p}+\delta_{11}^{(0)} & -\left(\frac{a_1p}{a_2\alpha_1}\right)^{(0)}E^{k_{0}-p} \\ (a_2\alpha_0)^{(0)}+(a_2\alpha_1)^{(0)}\frac{E^p}{p}+\delta_{21}^{(0)} & 1
		\end{pmatrix}& \text{if	} 0\in\Ts
	\end{dcases}\\
	\As_2^{(1)}&=\begin{dcases}
		\begin{pmatrix}
			\varepsilon^{(1)}_{11} & E^{k_{1}}a_1^{(1)}+\varepsilon_{12}^{(1)} \\ 1+\varepsilon^{(1)}_{21} & (a_2\alpha_0)^{(1)}+\varepsilon_{22}^{(1)}
		\end{pmatrix} & \text{if	} 1\in\Ss\\
		\begin{pmatrix}
			E^{k_{1}}a_1^{(1)} & \delta_{12}^{(1)} \\ (a_2\alpha_0)^{(1)}+\delta^{(1)}_{21} & 1+\delta_{22}^{(1)}
		\end{pmatrix}& \text{if	} 1\in\Ts
	\end{dcases}
\end{align*}

\begin{itemize}
	\item Both $\varepsilon_{j2}^{(1)}$ and $\delta_{j2}^{(1)}$ are the result of adding a $(p/a_2^{(0)})\varphi(E)^{k_{0}-p}$-multiple of the first column of $\As_2^{(1)}$ to the second column. Since $\varphi(E)\in pS_F$, $p/a_2^{(0)}\in \varpi \oh_F$ and $k_{0}-p\ge 2$ then $(p/a_2^{(0)})\varphi(E)^{k_{0}-p}\in \varpi p^2 S_F\subset pI$. Hence, $\varepsilon_{j2}^{(1)},\delta_{j2}^{(1)}\in pI$;
	
	\item $\tilde{\varepsilon}_{11}^{(0)}=\varepsilon_{11}^{(0)}-\varepsilon_{21}^{(0)}\left(\frac{a_1p}{a_2\alpha_1}\right)^{(0)}E^{k_{0}-p}\in I$;
	
	\item $\delta_{11}^{(0)}=-\delta_{21}^{(0)}\left(\frac{a_1p}{a_2\alpha_1}\right)^{(0)}E^{k_{0}-p}\in p I$ since $\delta_{21}^{(0)}\in p^2S_F$ and $p/a_2^{(0)}\in\varpi\oh_F$.
\end{itemize}

As mentioned above, our next step will be to scale the anti-diagonals of all $f$-matrices by an appropriate element, either $p$ or $p/a_2^{(0)}$ depending on the situation. The difficulty here is that any Type $\I$ matrix will be made non-integral if we do this as they stand now. Our answer to this problem will be to perform row and column permutation operations depending on the Type combinations present in order to `line up' the anti-diagonals so they have sufficient powers of $p$ as to not become non-integral after scaling. We note that we are not able to cover all cases with our methods however.\footnote{As we will see later in Section \ref{complete-f=2-sec}, there are ways around this problem when $f=2$ but those methods require matrix by matrix considerations that are impractical for the general $f>2$ case.}

\begin{defin}\label{sandwich-subset-def}
	Recall that $\Ss\subset \Z/f\Z$ denotes the set of those $i$ for which $\As^{(i)}$ is a Type $\I$ Frobenius matrix. We define a subset $\Ws\subset \Z/f\Z$ with the following steps:
	\begin{enumerate}[label*=\arabic*.]
		\item Label the elements of $\Ss\setminus \{0\}$ in $\Z/f\Z$ by $r_j$ so that $r_1<r_2<\cdots<r_d$ for some $d$ in the range $1\le d\le f-2$. Note that this set may be empty, in which case $\Ws$ is empty.
		
		\item Include $i\in\Ws$ for all $r_{1}\le i< r_2$.
		
		\item Now include $i\in\Ws$ for all $r_{3}\le i< r_4$.
		
		\item Continue this pattern in that we include $i\in\Ws$ for all $r_{1+n}\le i< r_{2+n}$ for $1< n< d-1$ and $n$ even.
		
		\item If $d$ is odd, include $i\in\Ws$ for all $r_d\le i\le f-1$.
	\end{enumerate}
\end{defin}

The effect of these steps is to include, starting at $i=1$, all $i\in\Ws$ which correspond to type combinations of the form $(\I,\II,,\II,\dots,\II,\I)$ excluding the index corresponding to the right most Type $\I$. We will lovingly refer to these combinations as \textit{Type $\II$ sandwiches}. When $d$ is odd then there is a left over Type $\I$ leaving us with an open-faced Type $\II$ sandwich meaning that the right most Type $\I$ does not exist. As a result we are forced to include all $i$ up to $i=f-1$ and perform a column permutation on the $i=0$ matrix.

\begin{example}
	Suppose $f=9$ and $\Ms(\As)$ has Type combination given by
	\[(\I,\I,\II,\II,\I,\II,\I,\I,\I)\]
	Then $d=5$ with $r_1=1$, $r_2=4$, $r_3=6$, $r_4=7$ and $r_5=8$. Hence, we get the set $\Ws=\{1,2,3,6,8\}\subset\Z/9\Z$ with $i=8$ being included since $d$ is odd.
\end{example}

Consider the base change over $\GL_2(F)$, intended to line up the anti-diagonals of each matrix given by $X_3^{(i)}=\begin{psmallmatrix}
	0 & 1 \\ 1 & 0
\end{psmallmatrix}$ for $i\in\Ws$ and $X_3^{(i)}=I_2$ for $i\notin\Ws$. With heartfelt apologies to the reader, the result of $(X_3)*_\varphi (\As_2)$ will be given by
\begin{align*}
		\As_3^{(0)}&=\begin{dcases}
		\begin{pmatrix}
			-\left(\frac{a_1\alpha_0p}{\alpha_1}\right)^{(0)}E^{k_{0}-p} & -\left(\frac{a_1p}{a_2\alpha_1}\right)^{(0)}E^{k_{0}-p}+\tilde{\varepsilon}_{11}^{(0)} \\ (a_2\alpha_0)^{(0)}+(a_2\alpha_1)^{(0)}\frac{E^p}{p} & 1+\varepsilon_{21}^{(0)}
		\end{pmatrix} & \text{if	} 0\in\Ss\text{ and } f-1\in\Ws\\
		\begin{pmatrix}
			-\left(\frac{a_1p}{a_2\alpha_1}\right)^{(0)}E^{k_{0}-p}+\tilde{\varepsilon}_{11}^{(0)} & -\left(\frac{a_1\alpha_0p}{\alpha_1}\right)^{(0)}E^{k_{0}-p} \\ 1+\varepsilon_{21}^{(0)} & (a_2\alpha_0)^{(0)}+(a_2\alpha_1)^{(0)}\frac{E^p}{p}
		\end{pmatrix} & \text{if	} 0\in\Ss\text{ and } f-1\notin\Ws\\
		\begin{pmatrix}
			-\left(\frac{a_1p}{a_2\alpha_1}\right)^{(0)}E^{k_{0}-p} & 	-\left(\frac{a_1\alpha_0p}{\alpha_1}\right)^{(0)}E^{k_{0}-p}+\delta_{11}^{(0)}\\ 1 & (a_2\alpha_0)^{(0)}+(a_2\alpha_1)^{(0)}\frac{E^p}{p}+\delta_{21}^{(0)}
		\end{pmatrix}& \text{if	} 0\in\Ts \text{ and } f-1\in\Ws \\
		\begin{pmatrix}
			-\left(\frac{a_1\alpha_0p}{\alpha_1}\right)^{(0)}E^{k_{0}-p}+\delta_{11}^{(0)} & -\left(\frac{a_1p}{a_2\alpha_1}\right)^{(0)}E^{k_{0}-p} \\ (a_2\alpha_0)^{(0)}+(a_2\alpha_1)^{(0)}\frac{E^p}{p}+\delta_{21}^{(0)} & 1
		\end{pmatrix}& \text{if	} 0\in\Ts \text{ and } f-1\notin\Ws \\
	\end{dcases}\\
	\As_3^{(1)}&=\begin{dcases}
		\begin{pmatrix}
			1+\varepsilon^{(1)}_{21} & (a_2\alpha_0)^{(1)}+\varepsilon_{22}^{(1)} \\ \varepsilon^{(1)}_{11} & E^{k_{1}}a_1^{(1)}+\varepsilon_{12}^{(1)}
		\end{pmatrix} & \text{if	} 1\in\Ss \text{ so that } 1\in\Ws\\
		\begin{pmatrix}
			E^{k_{1}}a_1^{(1)} & \delta_{12}^{(1)} \\ (a_2\alpha_0)^{(1)}+\delta^{(1)}_{21} & 1+\delta_{22}^{(1)}
		\end{pmatrix}& \text{if	} 1\in\Ts \text{ so that } 1\notin\Ws
	\end{dcases}\\
	\As_3^{(i)}&=\begin{dcases}
		\begin{pmatrix}
			1+\varepsilon^{(i)}_{21} & (a_2\alpha_0)^{(i)} \\ \varepsilon^{(i)}_{11} & E^{k_i}a_1^{(i)}
		\end{pmatrix} & \text{if	} 2\le i\le f-1 \text{ with } i \in\Ss \text{ and } i\in\Ws\\
		\begin{pmatrix}
			E^{k_i}a_1^{(i)} & 	\varepsilon^{(i)}_{11}\\ (a_2\alpha_0)^{(i)} & 1+\varepsilon^{(i)}_{21}
		\end{pmatrix}& \text{if	} 2\le i\le f-1 \text{ with } i \in\Ss \text{ and } i\notin\Ws \\
		\begin{pmatrix}
			1 & (a_2\alpha_0)^{(i)}+\delta^{(i)}_{21} \\  0 & E^{k_i}a_1^{(i)}
		\end{pmatrix}& \text{if	} 2\le i\le f-1 \text{ with } i \in\Ts \text{ and } i\in\Ws \\
		\begin{pmatrix}
			E^{k_i}a_1^{(i)} & 0 \\ (a_2\alpha_0)^{(i)}+\delta^{(i)}_{21} & 1
		\end{pmatrix}& \text{if	} 2\le i\le f-1 \text{ with } i \in\Ts \text{ and } i\notin\Ws.
	\end{dcases}
\end{align*}

\noindent As hinted in step 4 of Definition \ref{sandwich-subset-def}, there is a dichotomy with respect to the even or oddness of the order of the set $\Ss$. Indeed, these two distinct cases must be handled separately.

\subsubsection{Subcase: $|\Ss|$ Even} \label{even-subsec}
When $|\Ss|$ is even, then either $d$ in Definition \ref{sandwich-subset-def} is even or $d$ is odd with $\As^{(0)}$ being Type $\I$. As a result, we may write
\begin{align*}
	\As_3^{(0)}&=\begin{dcases}
		\begin{pmatrix}
			-\left(\frac{a_1\alpha_0p}{\alpha_1}\right)^{(0)}E^{k_{0}-p} & -\left(\frac{a_1p}{a_2\alpha_1}\right)^{(0)}E^{k_{0}-p}+\tilde{\varepsilon}_{11}^{(0)} \\ (a_2\alpha_0)^{(0)}+(a_2\alpha_1)^{(0)}\frac{E^p}{p} & 1+\varepsilon_{21}^{(0)}
		\end{pmatrix} & \text{if	} 0\in\Ss\\
		\begin{pmatrix}
			-\left(\frac{a_1\alpha_0p}{\alpha_1}\right)^{(0)}E^{k_{0}-p}+\delta_{11}^{(0)} & -\left(\frac{a_1p}{a_2\alpha_1}\right)^{(0)}E^{k_{0}-p} \\ (a_2\alpha_0)^{(0)}+(a_2\alpha_1)^{(0)}\frac{E^p}{p}+\delta_{21}^{(0)} & 1
		\end{pmatrix}& \text{if	} 0\in\Ts. \\
	\end{dcases}
\end{align*}

We may notice that regardless of Type, scaling $\As_3^{(0)}$ by $p/a_2^{(0)}\in\varpi\oh_F$ along the diagonal will do the job of making all the terms either integral or lie the ideal $I$. Consider the base change over $\GL_2(F)$ intended to do just this given by $X_4^{(i)}=\Diag(1,p/a_2^{(0)})$ for all $i\in\Z/f\Z$. The result of $(X_4)*_\varphi (\As_3)$ will be given by
\begin{align*}
		\As_4^{(0)}&=\begin{dcases}
		\begin{pmatrix}
			-\left(\frac{a_1\alpha_0p}{\alpha_1}\right)^{(0)}E^{k_{0}-p} & -\left(\frac{a_1}{\alpha_1}\right)^{(0)}E^{k_{0}-p}+\frac{a_2^{(0)}}{p}\tilde{\varepsilon}_{11}^{(0)} \\ p\alpha_0^{(0)}+\alpha_1^{(0)}E^p & 1+\varepsilon_{21}^{(0)}
		\end{pmatrix} & \text{if	} 0\in\Ss\\
		\begin{pmatrix}
			-\left(\frac{a_1\alpha_0p}{\alpha_1}\right)^{(0)}E^{k_{0}-p}+\delta_{11}^{(0)} & -\left(\frac{a_1}{\alpha_1}\right)^{(0)}E^{k_{0}-p} \\ p\alpha_0^{(0)}+\alpha_1^{(0)}E^p+\frac{p}{a_2^{(0)}}\delta_{21}^{(0)} & 1
		\end{pmatrix}& \text{if	} 0\in\Ts \\
	\end{dcases}\\
	\As_4^{(1)}&=\begin{dcases}
		\begin{pmatrix}
			1+\varepsilon^{(1)}_{21} & \frac{a_2^{(0)}}{p}(a_2\alpha_0)^{(1)}+\frac{a_2^{(0)}}{p}\varepsilon_{22}^{(1)} \\ \frac{p}{a_2^{(0)}}\varepsilon^{(1)}_{11} & E^{k_{1}}a_1^{(1)}+\varepsilon_{12}^{(1)}
		\end{pmatrix} & \text{if	} 1\in\Ss \\
		\begin{pmatrix}
			E^{k_{1}}a_1^{(1)} & \frac{a_2^{(0)}}{p}\delta_{12}^{(1)} \\ \frac{p}{a_2^{(0)}}(a_2\alpha_0)^{(1)}+\frac{p}{a_2^{(0)}}\delta^{(1)}_{21} & 1+\delta_{22}^{(1)}
		\end{pmatrix}& \text{if	} 1\in\Ts
	\end{dcases}\\
	\As_4^{(i)}&=\begin{dcases}
		\begin{pmatrix}
			1+\varepsilon^{(i)}_{21} & \frac{a_2^{(0)}}{p}(a_2\alpha_0)^{(i)} \\ \frac{p}{a_2^{(0)}}\varepsilon^{(i)}_{11} & E^{k_i}a_1^{(i)}
		\end{pmatrix} & \text{if	} 2\le i\le f-1 \text{ and } i \in\Ss \text{ and } i\in\Ws\\
		\begin{pmatrix}
			E^{k_i}a_1^{(i)} & 	\frac{a_2^{(0)}}{p}\varepsilon^{(i)}_{11}\\ \frac{p}{a_2^{(0)}}(a_2\alpha_0)^{(i)} & 1+\varepsilon^{(i)}_{21}
		\end{pmatrix}& \text{if	} 2\le i\le f-1 \text{ and } i \in\Ss \text{ and } i\notin\Ws \\
		\begin{pmatrix}
			1 & \frac{a_2^{(0)}}{p}(a_2\alpha_0)^{(i)}+\frac{a_2^{(0)}}{p}\delta^{(i)}_{21} \\  0 & E^{k_i}a_1^{(i)}
		\end{pmatrix}& \text{if	} 2\le i\le f-1 \text{ and } i \in\Ts \text{ and } i\in\Ws \\
		\begin{pmatrix}
			E^{k_i}a_1^{(i)} & 0 \\ \frac{p}{a_2^{(0)}}(a_2\alpha_0)^{(i)}+\frac{p}{a_2^{(0)}}\delta^{(i)}_{21} & 1
		\end{pmatrix}& \text{if	} 2\le i\le f-1 \text{ and } i \in\Ts \text{ and } i\notin\Ws.
	\end{dcases}
\end{align*}

Since $p/a_2^{(0)}$ is integral then we need only be concerned with the top right entry of each matrix. To this end, let us now make the following assumption:
\[\nu_p(a_2^{(i)})\ge 1-\nu_p(a_2^{(0)})\hspace{.5cm}\text{\normalfont for all}\hspace{.5cm} i\in\Ws.\]
Notice that under this assumption, we now have $v(i)=1$ for all $i\in\Ws$ where $v(i)$ denotes the non-zero valuation indicator function for $a_2^{(i)}$. Make the observation that if $i\notin\Ws$ but $i\in\Ss$ then we are guaranteed to have $i-1\in\Ws$. We analyze the potentially problematic top right entries individually:

\begin{itemize}
	\item $\frac{a_2^{(0)}}{p}\varepsilon^{(i)}_{11}\in \frac{a_2^{(0)}}{p}E^{k_i}p^{v(i-1)}I\subset I$ since $i-1\in\Ws$;
	
	\item $\frac{a_2^{(0)}}{p}\delta^{(i)}_{21}\in \frac{a_2^{(0)}}{p}p^{v(i-1)}I\subset I$ since $i-1\in\Ws$;
	
	\item $\frac{a_2^{(0)}}{p}\varepsilon_{22}^{(f-2)},\frac{a_2^{(0)}}{p}\delta_{12}^{(f-2)}\in \frac{a_2^{(0)}}{p}p I\subset I$;
	
	\item $\frac{a_2^{(0)}}{p}\tilde{\varepsilon}_{11}^{(0)}=\frac{a_2^{(0)}}{p}(\varepsilon_{11}^{(0)}-\varepsilon_{21}^{(0)}\left(\frac{a_1p}{a_2\alpha_1}\right)^{(0)}E^{k_{0}-p})$ so that the only possible problem would occur in $\frac{a_2^{(0)}}{p}\varepsilon_{11}^{(0)}\in \frac{a_2^{(0)}}{p}E^{k_{0}}p^{v(f-1)}I$. However, we know that $f-1\in\Ws$ since $d$ is odd. Hence, $\frac{a_2^{(0)}}{p}\tilde{\varepsilon}_{11}^{(0)}\in I$.
\end{itemize}

Reapplying the $\Ws$-derived row swap $X_5^{(i)}=X_3^{(i)}$ has the effect of putting ourselves back in the original orientation prior to the lineup of anti-diagonals for scaling. As a result, we have our next case of satisfied descent assumptions.

\begin{prop}\label{dec-assump-even-prop}
	Suppose $\Ms(\As)$ is a rank-two Kisin module over $S_F$ with labeled heights $p+2\le k_{0}\le 2p-4$ and $2\le k_i\le p-3$ for $i\neq 0$. If $|\Ss|$ is even and $\nu_p(a_2^{(i)})\ge 1-\nu_p(a_2^{(0)})$ for all $i\in\Ws$, then the Frobenius matrix $(\As)$ satisfies the Descent Assumptions \ref{descent-assump} with
\begin{align*}
		\Af_0^{(0)}&=
	\begin{pmatrix}
		-\left(\frac{a_1\alpha_0p}{\alpha_1}\right)^{(0)}E^{k_{0}-p} & -\left(\frac{a_1}{\alpha_1}\right)^{(0)}E^{k_{0}-p} \\ p\alpha_0^{(0)}+\alpha_1^{(0)}E^p & 1
	\end{pmatrix}\\
	\Af_0^{(i)}&=\begin{dcases}
		\begin{pmatrix}
			1 & \frac{a_2^{(0)}}{p}(a_2\alpha_0)^{(i)} \\ 0 & E^{k_i}a_1^{(i)}
		\end{pmatrix} & \text{if	} 1\le i\le f-1 \text{ and } i\in\Ws\\
		\begin{pmatrix}
			E^{k_i}a_1^{(i)} & 0\\ \frac{p}{a_2^{(0)}}(a_2\alpha_0)^{(i)} & 1
		\end{pmatrix}& \text{if	} 1\le i\le f-1 \text{ and } i\notin\Ws.
	\end{dcases}
\end{align*}
\end{prop}

\subsubsection{Subcase: $|\Ss|$ Odd} \label{odd-subsec}
Let us return ourselves to $\As_3^{(i)}$ as stated prior to the even subcase \ref{even-subsec}. When $|\Ss|$ is odd, then either $d$ in Definition \ref{sandwich-subset-def} is odd or $d$ is even with $\As^{(0)}$ being Type $\I$. Hence, we have that
\begin{align*}
		\As_3^{(0)}&=\begin{dcases}
		\begin{pmatrix}
			-\left(\frac{a_1p}{a_2\alpha_1}\right)^{(0)}E^{k_{0}-p}+\tilde{\varepsilon}_{11}^{(0)} & -\left(\frac{a_1\alpha_0p}{\alpha_1}\right)^{(0)}E^{k_{0}-p} \\ 1+\varepsilon_{21}^{(0)} & (a_2\alpha_0)^{(0)}+(a_2\alpha_1)^{(0)}\frac{E^p}{p}
		\end{pmatrix} & \text{if	} 0\in\Ss\\
		\begin{pmatrix}
			-\left(\frac{a_1p}{a_2\alpha_1}\right)^{(0)}E^{k_{0}-p} & 	-\left(\frac{a_1\alpha_0p}{\alpha_1}\right)^{(0)}E^{k_{0}-p}+\delta_{11}^{(0)}\\ 1 & (a_2\alpha_0)^{(0)}+(a_2\alpha_1)^{(0)}\frac{E^p}{p}+\delta_{21}^{(0)}
		\end{pmatrix}& \text{if	} 0\in\Ts\\
	\end{dcases}\\
		\As_3^{(1)}&=\begin{dcases}
		\begin{pmatrix}
			1+\varepsilon^{(1)}_{21} & (a_2\alpha_0)^{(1)}+\varepsilon_{22}^{(1)} \\ \varepsilon^{(1)}_{11} & E^{k_{1}}a_1^{(1)}+\varepsilon_{12}^{(1)}
		\end{pmatrix} & \text{if	} 1\in\Ss\\
		\begin{pmatrix}
			E^{k_{1}}a_1^{(1)} & \delta_{12}^{(1)} \\ (a_2\alpha_0)^{(1)}+\delta^{(1)}_{21} & 1+\delta_{22}^{(1)}
		\end{pmatrix}& \text{if	} 1\in\Ts
	\end{dcases}\\
	\As_3^{(i)}&=\begin{dcases}
		\begin{pmatrix}
			1+\varepsilon^{(i)}_{21} & (a_2\alpha_0)^{(i)} \\ \varepsilon^{(i)}_{11} & E^{k_i}a_1^{(i)}
		\end{pmatrix} & \text{if	} 2\le i\le f-1 \text{ and } i \in\Ss \text{ and } i\in\Ws\\
		\begin{pmatrix}
			E^{k_i}a_1^{(i)} & 	\varepsilon^{(i)}_{11}\\ (a_2\alpha_0)^{(i)} & 1+\varepsilon^{(i)}_{21}
		\end{pmatrix}& \text{if	} 2\le i\le f-1 \text{ and } i \in\Ss \text{ and } i\notin\Ws \\
		\begin{pmatrix}
			1 & (a_2\alpha_0)^{(i)}+\delta^{(i)}_{21} \\  0 & E^{k_i}a_1^{(i)}
		\end{pmatrix}& \text{if	} 2\le i\le f-1 \text{ and } i \in\Ts \text{ and } i\in\Ws \\
		\begin{pmatrix}
			E^{k_i}a_1^{(i)} & 0 \\ (a_2\alpha_0)^{(i)}+\delta^{(i)}_{21} & 1
		\end{pmatrix}& \text{if	} 2\le i\le f-1 \text{ and } i \in\Ts \text{ and } i\notin\Ws.
	\end{dcases}
\end{align*}

Our strategy will be to scale by the anti-diagonals of each matrix by $p$ followed by a row operation that uses the $E^{k_{0}-p}$-term in $\As_3^{(0)}$ to kill the non-integral $E^p/p$-term. Indeed, consider the base change over $\GL_2(F)$ given by $X_4^{(i)}=\Diag(1,p)$ for all $i\in\Z/f\Z$. The result of $(X_4)*_\varphi(\As_3)$ will be given by
\begin{align*}
		\As_4^{(0)}&=\begin{dcases}
		\begin{pmatrix}
			-\left(\frac{a_1p}{a_2\alpha_1}\right)^{(0)}E^{k_{0}-p}+\tilde{\varepsilon}_{11}^{(0)} & -\left(\frac{a_1\alpha_0}{\alpha_1}\right)^{(0)}E^{k_{0}-p} \\ p+p\varepsilon_{21}^{(0)} & (a_2\alpha_0)^{(0)}+(a_2\alpha_1)^{(0)}\frac{E^p}{p}
		\end{pmatrix} & \text{if	} 0\in\Ss\\
		\begin{pmatrix}
			-\left(\frac{a_1p}{a_2\alpha_1}\right)^{(0)}E^{k_{0}-p} & 	-\left(\frac{a_1\alpha_0}{\alpha_1}\right)^{(0)}E^{k_{0}-p}+\frac{1}{p}\delta_{11}^{(0)}\\ p & (a_2\alpha_0)^{(0)}+(a_2\alpha_1)^{(0)}\frac{E^p}{p}+\delta_{21}^{(0)}
		\end{pmatrix}& \text{if	} 0\in\Ts\\
	\end{dcases}\\
		\As_4^{(1)}&=\begin{dcases}
		\begin{pmatrix}
			1+\varepsilon^{(1)}_{21} & \frac{1}{p}(a_2\alpha_0)^{(1)}+\frac{1}{p}\varepsilon_{22}^{(1)} \\ p\varepsilon^{(1)}_{11} & E^{k_{1}}a_1^{(1)}+\varepsilon_{12}^{(1)}
		\end{pmatrix} & \text{if	} 1\in\Ss\\
		\begin{pmatrix}
			E^{k_{1}}a_1^{(1)} & \frac{1}{p}\delta_{12}^{(1)} \\ p(a_2\alpha_0)^{(1)}+p\delta^{(1)}_{21} & 1+\delta_{22}^{(1)}
		\end{pmatrix}& \text{if	} 1\in\Ts
	\end{dcases} \\
	\As_4^{(i)}&=\begin{dcases}
		\begin{pmatrix}
			1+\varepsilon^{(i)}_{21} & \frac{1}{p}(a_2\alpha_0)^{(i)} \\ p\varepsilon^{(i)}_{11} & E^{k_i}a_1^{(i)}
		\end{pmatrix} & \text{if	} 2\le i\le f-1 \text{ and } i \in\Ss \text{ and } i\in\Ws\\
		\begin{pmatrix}
			E^{k_i}a_1^{(i)} & 	\frac{1}{p}\varepsilon^{(i)}_{11}\\ p(a_2\alpha_0)^{(i)} & 1+\varepsilon^{(i)}_{21}
		\end{pmatrix}& \text{if	} 2\le i\le f-1 \text{ and } i \in\Ss \text{ and } i\notin\Ws \\
		\begin{pmatrix}
			1 & \frac{1}{p}(a_2\alpha_0)^{(i)}+\frac{1}{p}\delta^{(i)}_{21} \\  0 & E^{k_i}a_1^{(i)}
		\end{pmatrix}& \text{if	} 2\le i\le f-1 \text{ and } i \in\Ts \text{ and } i\in\Ws \\
		\begin{pmatrix}
			E^{k_i}a_1^{(i)} & 0 \\ p(a_2\alpha_0)^{(i)}+p\delta^{(i)}_{21} & 1
		\end{pmatrix}& \text{if	} 2\le i\le f-1 \text{ and } i \in\Ts \text{ and } i\notin\Ws.
	\end{dcases}
\end{align*}

As with the analogous scaling operation we did in the even case \ref{even-subsec}, the only threats to non-integrality lie in the upper right hand entries of each matrix. As a result, we will now make the following assumption:
\[\nu_p(a_2^{(i)})\ge 1\hspace{.5cm}\text{\normalfont for all}\hspace{.5cm} i\in\Ws.\]
Under these assumptions, we need only ensure that the $I$-terms remain in $I$ after division by $p$. The analysis of these terms is largely the same as in the even case \ref{even-subsec} save for the $i=0$ matrix which is actually easier to see. Indeed, it remains true that $\frac{1}{p}\varepsilon^{(i)}_{11}, \frac{1}{p}\delta^{(i)}_{21}, \frac{1}{p}\varepsilon_{22}^{(1)},  \frac{1}{p}\delta_{12}^{(1)}\in I$ for $2\le i$. Since $\delta_{11}^{(0)}\in pI$ then it is easy to see that $\frac{1}{p}\delta_{11}^{(0)}\in I$.

As mentioned prior, we want to kill the $E^p/p$-term in $\As_4^{(0)}$ which by the above analysis, is the only non-integral or non-$I$ term remaining. Consider the base change over $\GL_2(S_F)$ given by
\begin{align*}
	X_5^{(0)}&=\begin{pmatrix}
		1 & 0 \\ (\frac{a_2\alpha_1^2}{a_1\alpha_0p})^{(0)}E^{2p-k_{0}} & 1
	\end{pmatrix} & X_5^{(i)}=\begin{pmatrix}
	1 & 0 \\ 0 & 1
	\end{pmatrix}\text{ for }i\neq 0.
\end{align*}

\noindent The result of $(X_5)*_\varphi(\As_4)$ will leave $\As_5^{(i)}=\As_4^{(i)}$ for $2\le i\le f-1$ and 
\begin{align*}
	\As_5^{(0)}&=\begin{dcases}
		\begin{pmatrix}
			-\left(\frac{a_1p}{a_2\alpha_1}\right)^{(0)}E^{k_{0}-p}+\tilde{\varepsilon}_{11}^{(0)} & -\left(\frac{a_1\alpha_0}{\alpha_1}\right)^{(0)}E^{k_{0}-p} \\ p-(\frac{\alpha_1}{\alpha_0})^{(0)}E^p+\tilde{\varepsilon}_{21}^{(0)} & (a_2\alpha_0)^{(0)}
		\end{pmatrix} & \text{if	} 0\in\Ss\\
		\begin{pmatrix}
			-\left(\frac{a_1p}{a_2\alpha_1}\right)^{(0)}E^{k_{0}-p} & 	-\left(\frac{a_1\alpha_0}{\alpha_1}\right)^{(0)}E^{k_{0}-p}+\frac{1}{p}\delta_{11}^{(0)}\\ p-(\frac{\alpha_1}{\alpha_0})^{(0)}E^p & (a_2\alpha_0)^{(0)}+\tilde{\delta}_{22}^{(0)}
		\end{pmatrix}& \text{if	} 0\in\Ts
	\end{dcases}\\
	\As_5^{(1)}&=\begin{dcases}
		\begin{pmatrix}
			1+\tilde{\varepsilon}_{11}^{(1)} & \frac{1}{p}(a_2\alpha_0)^{(1)}+\frac{1}{p}\varepsilon_{22}^{(1)} \\ \tilde{\varepsilon}_{21}^{(1)} & E^{k_{1}}a_1^{(1)}+\varepsilon_{12}^{(1)}
		\end{pmatrix} & \text{if	} 1\in\Ss\\
		\begin{pmatrix}
			E^{k_{1}}a_1^{(1)}+\tilde{\delta}_{11}^{(1)} & \frac{1}{p}\delta_{12}^{(1)} \\ p(a_2\alpha_0)^{(1)}+\tilde{\delta}^{(1)}_{21} & 1+\delta_{22}^{(1)}
		\end{pmatrix}& \text{if	} 1\in\Ts
	\end{dcases}
\end{align*}

\begin{itemize}
	\item Each of $\tilde{\varepsilon}_{11}^{(1)}, \tilde{\varepsilon}_{21}^{(1)}, \tilde{\delta}_{11}^{(1)}, \tilde{\delta}^{(1)}_{21}$ result from adding a $-\varphi(\frac{a_2\alpha_1^2}{a_1\alpha_0p})^{(0)}\varphi(E)^{2p-k_{0}}$ multiple of the second columns of $\As_4^{(1)}$ to the first column. Since $\varphi(E)^{2p-k_{0}}\in p^4 S_F$ due to the fact that $2p-k_{0}\ge 4$, then all of the described terms lie in $\varpi p^3S_F\subset I$.
	
	\item $\tilde{\varepsilon}_{21}^{(0)}=p\varepsilon_{21}^{(0)}-\tilde{\varepsilon}_{11}^{(0)}(\frac{a_2\alpha_1^2}{a_1\alpha_0p})^{(0)}E^{2p-k_{0}}$. Since $\tilde{\varepsilon}_{11}^{(0)}=\varepsilon_{11}^{(0)}-\varepsilon_{21}^{(0)}\left(\frac{a_1p}{a_2\alpha_1}\right)^{(0)}E^{k_{0}-p}$ then we need only be concerned with the resulting $\varepsilon_{11}^{(0)}(\frac{a_2\alpha_1^2}{a_1\alpha_0p})^{(0)}E^{2p-k_{0}}$-term. However, recalling that $\varepsilon_{11}^{(0)}\in E^{k_{0}}I$ will imply that this term lies in $ E^{2p}p^{-1}I$ but $p\mid E^{2p}$ in $S_F$ so $\tilde{\varepsilon}_{21}^{(0)}\in I$.
	
	\item $\tilde{\delta}_{22}^{(0)}=\delta_{21}^{(0)}-\frac{1}{p}\delta_{11}^{(0)}(\frac{a_2\alpha_1^2}{a_1\alpha_0p})^{(0)}E^{2p-k_{0}}$. Since $\delta_{11}^{(0)}=\delta_{21}^{(0)}\left(\frac{a_1p}{a_2\alpha_1}\right)^{(0)}E^{k_{0}-p}$ and $\delta_{21}^{(0)}\in pI$ due to $0\in\Ws$, then it is easy to see that $\tilde{\delta}_{22}^{(0)}\in I$.
\end{itemize}

\noindent As a result, we have yet another case of satisfied descent assumptions.

\begin{prop}\label{dec-assump-odd-prop}
	Suppose $\Ms(\As)$ is a rank-two Kisin module over $S_F$ with labeled heights $p+2\le k_{0}\le 2p-4$ and $2\le k_i\le p-3$ for $i\neq 0$. If $|\Ss|$ is odd and $\nu_p(a_2^{(i)})\ge 1$ for all $i\in\Ws$, then the Frobenius matrix $(\As)$ satisfies the Descent Assumptions \ref{descent-assump} with
	\begin{align*}
				\Af_0^{(0)}&=
		\begin{pmatrix}
			-\left(\frac{a_1p}{a_2\alpha_1}\right)^{(0)}E^{k_{0}-p} & -\left(\frac{a_1\alpha_0}{\alpha_1}\right)^{(0)}E^{k_{0}-p} \\ p-(\frac{\alpha_1}{\alpha_0})^{(0)}E^p & (a_2\alpha_0)^{(0)}
		\end{pmatrix}\\
		\Af_0^{(i)}&=\begin{dcases}
			\begin{pmatrix}
				1 & \frac{1}{p}(a_2\alpha_0)^{(i)} \\ 0 & E^{k_i}a_1^{(i)}
			\end{pmatrix} & \text{if	} 1\le i\le f-1 \text{ and }i\in\Ws\\
			\begin{pmatrix}
					 E^{k_i}a_1^{(i)} & 0\\  p(a_2\alpha_0)^{(i)} & 1
			\end{pmatrix}& \text{if	} 1\le i\le f-1 \text{ and }i\notin\Ws.
		\end{dcases}
	\end{align*}
\end{prop}

%%%%%%%%%%%%%%%%%%%%%%%%%%%%%%%%%%%%%%%%%%%%%%%%%%%%%%%%%%%%%%%%%%%%%%%%%%%%%%%%%%%%%%%
\subsection{\texorpdfstring{Complete Description in $f=2$}{Complete Description in f=2}}\label{complete-f=2-sec}
Restricting to the case of $f=2$ gives us much more control meaning a full description can be obtained. Using Prepositions \ref{irred-vals-prop} and \ref{const-kisin-prop}, we may obtain an explicit description of finite height, rank-two Kisin modules $\Ms(\As,\Bs)$ over $S_F$ in terms of the three irreducible Type combinations given below:
\begin{itemize}
	\item Type $(\I,\I)$: 	\begin{align*}
				\As& =\begin{pmatrix}
							0 & E^{k_0}a_1 \\ 1 & a_2\lambda_2^{k_1-k_0\varphi}
						\end{pmatrix} & \Bs&=\begin{pmatrix}
							0 & E^{k_1}b_1 \\ 1 & b_2\lambda_2^{k_0-k_1\varphi}
						\end{pmatrix}
			\end{align*}
	with $a_1,b_1\in\oh_F^*$ and $a_2\cdot b_2\in\varpi\oh_F$.
	
	\item Type $(\I,\II)$: 	\begin{align*}
				\As&=\begin{pmatrix}
							0 & E^{k_0}a_1 \\ 1 & a_2\lambda_4^{(1-\varphi^2)(k_1+k_0\varphi)}\end{pmatrix}&\Bs&=\begin{pmatrix}
							E^{k_1}b_1 & 0 \\ b_2\lambda_4^{(1-\varphi^2)(k_1\varphi-k_0)} & 1
						\end{pmatrix}
			\end{align*}
	with $a_1,b_1\in\oh_F^*$ and $a_2,b_2\in\varpi\oh_F$.
	
	\item Type $(\II,\I)$: 	\begin{align*}
				\As&=\begin{pmatrix}
							E^{k_0}a_1 & 0 \\ a_2\lambda_4^{(1-\varphi^2)(k_0\varphi-k_1)} & 1
						\end{pmatrix}&\Bs&=\begin{pmatrix}
							0 & E^{k_1}b_1 \\ 1 & b_2\lambda_4^{(1-\varphi^2)(k_1\varphi+k_0)}\end{pmatrix}
			\end{align*}
	with $a_1,b_1\in\oh_F^*$ and $a_2,b_2\in\varpi\oh_F$.
\end{itemize}
Recall that Type $(\II,\II)$ is reducible by \cite[Prop. 4.2.1]{Guz1}. Here we will continue to assume that we have labeled heights in the range $p+2\le k_0\le 2p-4$ and $2\le k_1\le p-3$. It is then clear that the results of Sections \ref{prep-large-val-sec} and \ref{prep-small-val-sec} will provide satisfied descent assumptions outside of two cases. The first being Type $(\I,\I)$ with $0\le \nu_p(b_2)<1-\nu_p(a_2)$ and the second being Type $(\II,\I)$ with $0<\nu_p(b_2)<1$. Note that Type $(\I,\II)$ is covered by Proposition \ref{dec-assump-odd-prop} as $\Ws=\emptyset$ in this Type combination. Hence, to get a full description in the case of $f=2$, we need to treat these two cases separately.

\subsubsection{Type $(\I,\I)$ with $0\le \nu_p(a_2)<1-\nu_p(b_2)$}\label{Type-I,I-subsec}
Let us recall the following Type $\I$ matrices $\As_3^{(i)}$ found in Section \ref{even-subsec},
\begin{align*}	
	\As_3^{(0)}&=\begin{pmatrix}
		-\left(\frac{a_1\alpha_0p}{\alpha_1}\right)^{(0)}E^{k_{0}-p} & -\left(\frac{a_1p}{a_2\alpha_1}\right)^{(0)}E^{k_{0}-p}+\tilde{\varepsilon}_{11}^{(0)} \\ (a_2\alpha_0)^{(0)}+(a_2\alpha_1)^{(0)}\frac{E^p}{p} & 1+\varepsilon_{21}^{(0)}
	\end{pmatrix}\\
	\As_3^{(1)}&=\begin{pmatrix}
		1+\varepsilon^{(1)}_{21} & (a_2\alpha_0)^{(1)}+\varepsilon_{22}^{(1)} \\ \varepsilon^{(1)}_{11} & E^{k_{1}}a_1^{(1)}+\varepsilon_{12}^{(1)}
	\end{pmatrix}.
\end{align*}

We will change notation to ease the decomposition superscripts in the following manner: Set $\As_3=\As_3^{(0)}$ and drop all of the superscripts $(0)$. Also set $\Bs_3=\As_3^{(1)}$ and make the replacements $b_j=a_j^{(1)}$, $\beta_j=\alpha_j^{(1)}$ and $\delta_{jk}=\varepsilon_{jk}^{(1)}$. Hence after performing this change, we are left with
\begin{align*}
		\As_3&=\begin{pmatrix}
			-\left(\frac{a_1\alpha_0p}{\alpha_1}\right)E^{k_{0}-p} & -\left(\frac{a_1p}{a_2\alpha_1}\right)E^{k_{0}-p}+\tilde{\varepsilon}_{11} \\ a_2\alpha_0+(a_2\alpha_1)\frac{E^p}{p} & 1+\varepsilon_{21}
		\end{pmatrix} \\
		\Bs_3&=\begin{pmatrix}
		1+\delta_{21} & b_2\beta_0+\delta_{22} \\ \delta_{11} & E^{k_{1}}b_1+\delta_{12}
	\end{pmatrix}.
\end{align*}

\noindent Since we know that the product $a_2b_2\in\varpi\oh_F$ by Proposition \ref{irred-vals-prop} then we can write $\delta_{2j},\delta_{12},\varepsilon_{21}\in pI$ and $\delta_{11},\tilde{\varepsilon}_{11}\in I$. 

We will now stray away from our normal strategy to instead apply a row operation on $\Bs_3$ intending to kill the $E^{k_1}b_1$-term. The reason for this deviation is that the $p$-adic valuation of $b_2$ is no longer large enough to absorb a division by $p$. As a result, we need to create another term that can be divided by $p$. To prepare for this however, we first scale the anti-diagonals of both matrices by $b_2\in\oh_F$ using the following base change over $\GL_2(F)$ given by $X_4=Y_4=\Diag(1,b_2)$. The result of $(X_4,Y_4)*_\varphi(\As_3,\Bs_3)$ will then be given by 
\begin{align*}
	\As_4&=\begin{pmatrix}
		-\left(\frac{a_1\alpha_0p}{\alpha_1}\right)E^{k_{0}-p} & -\left(\frac{a_1p}{a'_2\alpha_1}\right)E^{k_{0}-p}+\frac{1}{b_2}\tilde{\varepsilon}_{11} \\ a'_2\beta_0+(a'_2\alpha_1)\frac{E^p}{p} & 1+\varepsilon_{21}
	\end{pmatrix} \\
	\Bs_4&=\begin{pmatrix}
		1+\delta_{21} & \beta_0+\frac{1}{b_2}\delta_{22} \\ b_2\delta_{11} & E^{k_{1}}b_1+\delta_{12}
	\end{pmatrix}
\end{align*}

\begin{itemize}
	\item $a'_2\coloneqq a_2b_2\in\varpi\oh_F$ so that $p/a'_2\in\varpi\oh_F$ under our assumptions;
	\item $\frac{1}{b_2}\delta_{22}\in I$ since $\delta_{22}\in pI$;
	\item $\frac{1}{b_2}\tilde{\varepsilon}_{11}=\frac{1}{b_2}(\varepsilon_{11}^{(0)}-\varepsilon_{21}^{(0)}\left(\frac{a_1p}{a_2\alpha_1}\right)^{(0)}E^{k_{0}-p})$ using the notation as in Section \ref{prep-small-val-sec}. Here we have that in our new notation, $\varepsilon_{11}^{(0)}\in E^{k_0}b_2p^2S_F$ by Lemma \ref{phi-admis-lem}. Let us write $E^{k_0}b_2p^2S_F=E^{k_0-p}(E^p)b_2p^2S_F$ and observe that $(pE^p)S_F\in p\Sig_F+p^2\Fil^{2p}S_F$. Hence, $\varepsilon_{11}^{(0)}\in pb_2E^{k_{0}-p}(p\Sig_F+p^2\Fil^{2p}S_F)\subset b_2(p^2\Sig_F+p^3E\Fil^{2p}S_F)\subset b_2(p^2\Sig_F+p^3I)$. Moreover, since $p/(a_2b_2)\in \varpi\oh_F$ then we see that $\frac{1}{b_2}\tilde{\varepsilon}_{11}\in p^2\Sig_F+pI$ after recalling that $\varepsilon_{21}^{(0)}\in pI$.
\end{itemize}

Let us now kill the $E^{k_1}$-term in $\Bs_4$ via the following base change over $\GL_2(\Sig_F)$,
\begin{align*}
	X_5&=\begin{pmatrix}
	1 & 0 \\ 0 & 1
	\end{pmatrix} & Y_5&=\begin{pmatrix}
	1 & 0 \\ -\frac{b_1}{\beta_0}E^{k_1} & 1
	\end{pmatrix}.
\end{align*}

\noindent The result of $(X_5,Y_5)*_\varphi(\As_4,\Bs_4)$ will then be given by
\begin{align*}
		\As_5&=\begin{pmatrix}
		-\left(\frac{a_1\alpha_0p}{\alpha_1}\right)E^{k_{0}-p}+\varepsilon'_{11} & -\left(\frac{a_1p}{a'_2\alpha_1}\right)E^{k_{0}-p}+\frac{1}{b_2}\tilde{\varepsilon}_{11} \\ a'_2\alpha_0+(a'_2\alpha_1)\frac{E^p}{p}+\varepsilon'_{21} & 1+\varepsilon_{21}
	\end{pmatrix} \\
	\Bs_5&=\begin{pmatrix}
		1+\delta_{21} & \beta_0+\frac{1}{b_2}\delta_{22} \\  -\frac{b_1}{\beta_0}E^{k_1}+\delta'_{12} & \delta'_{22}
	\end{pmatrix}
\end{align*}

\begin{itemize}
	\item $\delta'_{12}=b_2\delta_{11}-\delta_{21}\frac{b_1}{\beta_0}E^{k_1}\in I$.
	\item $\delta'_{22}=\delta_{12}-\frac{1}{b_2}\delta_{22}\frac{b_1}{\beta_0}E^{k_1}\in p^2\Sig_F+pI$.
	\item $\varepsilon'_{j1}\in I$ since they both result from adding a $\varphi(E)^{k_1}\in p^2S_F\subset I$ multiple of the second column of $\As_3$ to first column.
\end{itemize}

\begin{rem}\label{why-f=2-nice-rem}
	When $f>2$, this step of killing the $E^{k_i}$-term can be done as we have performed it here. However, we would at some point be required to scale by $a_2^{(i)}$ for all $i\in\Ws$ to clear denominators which leads to a proliferation of subcases depending on the valuations of each $a_2^{(i)}$ in relation to one another. It is for this reason that we restrict ourselves to $f=2$ when computing these very small valuation cases.
\end{rem}

We will now return to the strategy employed in Section \ref{odd-subsec}. In particular, let us swap the columns of $\As_5$ which swaps the rows of $\Bs_5$ and scale the anti-diagonals of both matrices by $p$ via $X_6=Y_6=\Diag(1,p)$. The result is
\begin{align*}
		\As_6&=\begin{pmatrix}
		-\left(\frac{a_1p}{a'_2\alpha_1}\right)E^{k_{0}-p}+\frac{1}{b_2}\tilde{\varepsilon}_{11} & -\left(\frac{a_1\alpha_0}{\alpha_1}\right)E^{k_{0}-p}+\frac{1}{p}\varepsilon'_{11} \\ p+p\varepsilon_{21} & a'_2\alpha_0+(a'_2\alpha_1)\frac{E^p}{p}+\varepsilon'_{21}
	\end{pmatrix} \\
	\Bs_6&=\begin{pmatrix}
		-\frac{b_1}{\beta_0}E^{k_1}+\delta'_{12} & \frac{1}{p}\delta'_{22} \\ p+p\delta_{21} & \beta_0+\frac{1}{b_2}\delta_{22}
	\end{pmatrix}
\end{align*}

\begin{itemize}
	\item $\frac{1}{p}\delta'_{22}\in p\Sig_F+I$ since $\delta'_{22}\in p^2\Sig_F+pI$;
	
	\item $\frac{1}{p}\varepsilon'_{11}=\frac{1}{p}(-(\frac{a_1p}{a'_2\alpha_1})E^{k_{0}-p}+\frac{1}{b_2}\tilde{\varepsilon}_{11})\frac{b_1}{\varphi(\beta_0)}\varphi(E)^{k_1}\in I$ since $\varphi(E)^{k_1}\in p^2S_F$ and $p/a'_2\in\varpi\oh_F$.
\end{itemize}

Our final operation will be to kill the $E^p/p$-term in $\As_6$ via an operation analogous to $X_5^{(i)}$ in Section \ref{odd-subsec}. Indeed, consider the base change over $\GL_2(S_F)$ given by
\begin{align*}
	X_7&=\begin{pmatrix}
		1 & \frac{a'_2\alpha_1^2}{a_1\alpha_0p}E^{2p-k_0} \\ 0 & 1
	\end{pmatrix} &
	Y_7&=\begin{pmatrix}
		1 & 0 \\ 0 & 1
	\end{pmatrix}.
\end{align*}

\noindent Then the result of $(X_7,Y_7)*_\varphi(\As_6,\Bs_6)$ will be given by 
\begin{align*}
		\As_6&=\begin{pmatrix}
		-\left(\frac{a_1p}{a'_2\alpha_1}\right)E^{k_{0}-p}+\frac{1}{b_2}\tilde{\varepsilon}_{11} & -\left(\frac{a_1\alpha_0}{\alpha_1}\right)E^{k_{0}-p}+\frac{1}{p}\varepsilon'_{11} \\ p-\frac{\alpha_1}{\alpha_0}E^p+\varepsilon''_{21} & a'_2\alpha_0+\varepsilon_{22}
	\end{pmatrix} \\
	\Bs_6&=\begin{pmatrix}
		-\frac{b_1}{\beta_0}E^{k_1}+\delta'_{12} & \delta''_{12} \\ p+p\delta_{21} & \beta_0+\delta''_{22}
	\end{pmatrix}
\end{align*}

\begin{itemize}
	\item Both $\delta''_{j2}$ result from adding a $-a'_2/p\varphi(E)^{2p-k_0}$ multiple of the first columns of $\Bs_5$ to the second column. Since $\varphi(E)^{2p-k_0}\in p^4S_F$ then $\delta''_{22}\in I$ and $\delta''_{12}\in p\Sig_F+I$;
	\item $\varepsilon''_{21}=p\varepsilon_{21}+\frac{1}{b_2}\tilde{\varepsilon}_{11}\frac{a'_2\alpha_1^2}{a_1\alpha_0p}E^{2p-k_0}\in I$ after recalling that $\frac{1}{b_2}\tilde{\varepsilon}_{11}\in p^2\Sig_F+pI$ and observing that $\varpi p\Sig_F\subset I$;
	
	\item $\varepsilon_{22}=\varepsilon'_{21}+\frac{1}{p}\varepsilon'_{11}\frac{a'_2\alpha_1^2}{a_1\alpha_0p}E^{2p-k_0}$. By the discussion of $\frac{1}{p}\varepsilon'_{11}$ above, the only possible complication could arise from the product $\left(\frac{a_1p}{a'_2\alpha_1}\right)E^{k_{0}-p}\frac{a'_2\alpha_1^2}{a_1\alpha_0p}E^{2p-k_0}=\frac{\alpha_1}{\alpha_0}E^p$ which is divisible by $p$ in $S_F$ so that $\varepsilon_{22}\in I$.
\end{itemize}

\noindent With this, we have achieved the penultimate case of satisfied descent assumptions.

\begin{prop}\label{dec-assump-I,I-prop}
	Suppose $\Ms(\As,\Bs)$ is a rank-two Kisin module over $S_F$ of Type $(\I,\I)$ with labeled heights $p+2\le k_{0}\le 2p-4$ and $2\le k_1\le p-3$. If $\nu_p(b_2)<1-\nu_p(a_2)$, then the Frobenius matrix $(\As,\Bs)$ satisfies the Descent Assumptions \ref{descent-assump} with
	\begin{align*}
		\Af_0&=\begin{pmatrix}
			-\left(\frac{a_1p}{a'_2\alpha_1}\right)E^{k_{0}-p}+p\Sig_F & -\left(\frac{a_1\alpha_0}{\alpha_1}\right)E^{k_{0}-p} \\ p-\frac{\alpha_1}{\alpha_0}E^p & a'_2\beta_0
		\end{pmatrix} \\
		\Bf_0&=\begin{pmatrix}
			-\frac{b_1}{\beta_0}E^{k_1} & p\Sig_F \\ p & \beta_0
		\end{pmatrix}
	\end{align*}
	where $p\Sig_F$ denotes an arbitrary element in said ring.
\end{prop}

\subsubsection{Type $(\II,\I)$ with $0< \nu_p(b_2)<1$}\label{Type-II,I-subsec}
Let us recall the following matrices $\As_3^{(i)}$ from section \ref{odd-subsec},
\begin{align*}
	\As_3^{(0)}&=\begin{pmatrix}
		-\left(\frac{a_1p}{a_2\alpha_1}\right)^{(0)}E^{k_{0}-p} & 	-\left(\frac{a_1\alpha_0p}{\alpha_1}\right)^{(0)}E^{k_{0}-p}+\delta_{11}^{(0)}\\ 1 & (a_2\alpha_0)^{(0)}+(a_2\alpha_1)^{(0)}\frac{E^p}{p}+\delta_{21}^{(0)}
	\end{pmatrix}\\
	\As_3^{(1)}&=
		\begin{pmatrix}
			1+\varepsilon^{(1)}_{21} & (a_2\alpha_0)^{(1)}+\varepsilon_{22}^{(1)} \\ \varepsilon^{(1)}_{11} & E^{k_{1}}a_1^{(1)}+\varepsilon_{12}^{(1)}
		\end{pmatrix}
\end{align*}

\noindent which are of Type $\II$ and Type $\I$ respectively. As we did before, we will change notation in the following way: Set $\As_3=\As_3^{(0)}$ and drop all of the superscripts $(0)$. Also set $\Bs_3=\As_3^{(1)}$ and make the replacements $b_j=a_j^{(1)}$, $\beta_j=\alpha_j^{(1)}$ and $\varepsilon_{jk}=\varepsilon_{jk}^{(1)}$. Hence, after performing this replacements we are left with
\begin{align*}
	\As_3&=\begin{pmatrix}
		-\left(\frac{a_1p}{a_2\alpha_1}\right)E^{k_{0}-p} & 	-\left(\frac{a_1\alpha_0p}{\alpha_1}\right)E^{k_{0}-p}+\delta_{11}\\ 1 & a_2\alpha_0+(a_2\alpha_1)\frac{E^p}{p}+\delta_{21}
	\end{pmatrix} \\
	\Bs_3&=\begin{pmatrix}
		1+\varepsilon_{21} & b_2\beta_0+\varepsilon_{22} \\ \varepsilon_{11} & E^{k_{1}}b_1+\varepsilon_{12}
	\end{pmatrix}.
\end{align*}

\noindent Since we know that both $a_2,b_2\in\varpi\oh_F$ by Proposition \ref{irred-vals-prop} then we can write that all of $\varepsilon_{jk},\delta_{j1}\in pI$.

We follow much of the same processes that we undertook in the Type $(\I,\I)$ case \ref{Type-I,I-subsec}. Indeed, since $p/b_2\in\varpi\oh_F$, then our first step will be to scale by $b_2$ along the anti-diagonals of both matrices via $X_4=Y_4=\Diag(1,b_2)\in\GL_2(F)$ so that
\begin{align*}
	\As_4&=\begin{pmatrix}
		-\left(\frac{a_1p}{a_2\alpha_1}\right)E^{k_{0}-p} & 	-\left(\frac{a_1\alpha_0p}{b_2\alpha_1}\right)E^{k_{0}-p}+\frac{1}{b_2}\delta_{11}\\ b_2 & a_2\alpha_0+(a_2\alpha_1)\frac{E^p}{p}+\delta_{21}
	\end{pmatrix} \\
	\Bs_4&=\begin{pmatrix}
		1+\varepsilon_{21} & \beta_0+\frac{1}{b_2}\varepsilon_{22} \\ b_2\varepsilon_{11} & E^{k_{1}}b_1+\varepsilon_{12}
	\end{pmatrix}
\end{align*}

\noindent where both $\frac{1}{a_2}\varepsilon_{22},\frac{1}{a_2}\delta_{11}\in I$ since $\varepsilon_{22},\delta_{11}\in pI$. 

As before, let us now kill the $E^{k_1}b_1$ term in $\Bs_4$ via the following base change over $\GL_2(\Sig_F)$
\begin{align*}
	X_5&=\begin{pmatrix}
		1 & 0 \\ 0 & 1
	\end{pmatrix} &
	Y_5&=\begin{pmatrix}
		1 & 0 \\ -\frac{b_1}{\beta_0}E^{k_1} & 1
	\end{pmatrix}.
\end{align*}

\noindent The result of $(X_5,Y_5)*_\varphi(\As_4,\Bs_4)$ will then be given by
\begin{align*}
	\As_5&=\begin{pmatrix}
		-\left(\frac{a_1p}{a_2\alpha_1}\right)E^{k_{0}-p}+\delta'_{11} & 	-\left(\frac{a_1\alpha_0p}{b_2\alpha_1}\right)E^{k_{0}-p}+\frac{1}{b_2}\delta_{11}\\ b_2+\delta'_{21} & a_2\alpha_0+(a_2\alpha_1)\frac{E^p}{p}+\delta_{21}
	\end{pmatrix} \\
	\Bs_5&=\begin{pmatrix}
		1+\varepsilon_{21} & \beta_0+\frac{1}{b_2}\varepsilon_{22} \\ -\frac{b_1}{\beta_0}E^{k_1}+\varepsilon'_{21} & \varepsilon'_{22}
	\end{pmatrix}
\end{align*}

\begin{itemize}
	\item $\varepsilon_{21}'=b_2\varepsilon_{11}-\varepsilon_{21}\frac{b_1}{\beta_0}E^{k_1}\in I$;
	\item $\varepsilon'_{22}=\varepsilon_{12}-\frac{1}{b_2}\varepsilon_{22}\frac{b_1}{\beta_0}E^{k_1}\in I$;
	\item Both $\delta'_{11}$ and $\delta'_{21}$ arise from adding a $\frac{b_1}{\varphi(\beta_0)}\varphi(E)^{k_1}$-multiple of the second columns of $\As_3$ to the first column. Since $b_1\varphi(E)^{k_1}\in b_1p^2S_F\subset \varpi p S_F\subset I$ then both $\delta'_{11},\delta'_{21}\in I$.
\end{itemize}

We will now swap the columns of $\As_5$ which swaps the rows of $\Bs_5$ and scale the anti-diagonals of the resulting matrices by $p/a_2\in\varpi\oh_F$ via $X_6=Y_6=\Diag(1,p/a_2)$ so that 
\begin{align*}
	\As_6&=\begin{pmatrix}
		-\left(\frac{a_1\alpha_0p}{b_2\alpha_1}\right)E^{k_{0}-p}+\frac{1}{b_2}\delta_{11} & -\left(\frac{a_1}{\alpha_1}\right)E^{k_{0}-p}+\frac{a_2}{p}\delta'_{11}\\ p\alpha_0+\alpha_1E^p+\frac{p}{a_2}\delta_{21} & b_2+\delta'_{21}
	\end{pmatrix}\\
	\Bs_6&=\begin{pmatrix}
		-\frac{b_1}{\beta_0}E^{k_1}+\varepsilon'_{21} & \frac{a_2}{p}\varepsilon'_{22} \\ \frac{p}{a_2}+\frac{p}{a_2}\varepsilon_{21} & \beta_0+\frac{1}{b_2}\varepsilon_{22}
	\end{pmatrix}
\end{align*}

\begin{itemize}
	\item $\frac{a_2}{p}\varepsilon'_{22}=\frac{a_2}{p}(\varepsilon_{12}-\frac{1}{b_2}\varepsilon_{22}\frac{b_1}{\beta_0}E^{k_1})$. Since $\varepsilon_{12}\in pI$ then we need only concern ourselves with $\frac{a_2}{p}(\frac{1}{b_2}\varepsilon_{22}\frac{b_1}{\beta_0}E^{k_1})$. However, we may recall that $\varepsilon_{22}=\varepsilon_{22}^{(1)}$ results as a multiple of $\frac{p}{a_2}\varphi(E)^{k_0-p}$ so that because $\varphi(E)^{k_0-p}\in p^2S_F$ then $\frac{a_2}{p}(\frac{1}{b_2}\varepsilon_{22}\frac{b_1}{\beta_0}E^{k_1})\in \frac{p^2}{b_2}S_F\subset \varpi p S_F\subset I$. Hence, $\frac{a_2}{p}\varepsilon'_{22}\in I$.
	\item $\frac{a_2}{p}\delta'_{11}\in I$ since $\delta'_{11}$ is divisible by $\varphi(E)^{k_1}\in p^2S_F$ so that $\frac{a_2}{p}\varphi(E)^{k_1}\in \varpi p S_F\subset I$.
\end{itemize}

\noindent At long last, we have reached our final case of satisfied descent assumptions.

\begin{prop}\label{dec-assump-II,I-prop}
	Suppose $\Ms(\As,\Bs)$ is a rank-two Kisin module over $S_F$ of Type $(\II,\I)$ with labeled heights $p+2\le k_{0}\le 2p-4$ and $2\le k_1\le p-3$. If $\nu_p(b_2)<1$, then the Frobenius matrix $(\As,\Bs)$ satisfies the Descent Assumptions \ref{descent-assump} with
\begin{align*}
	\Af_0&=\begin{pmatrix}
		-\left(\frac{a_1\alpha_0p}{b_2\alpha_1}\right)E^{k_{0}-p} & -\left(\frac{a_1}{\alpha_1}\right)E^{k_{0}-p} \\ p\alpha_0+\alpha_1E^p & b_2
	\end{pmatrix} \\
	\Bf_0&=\begin{pmatrix}
		-\frac{b_1}{\beta_0}E^{k_1} & 0 \\ \frac{p}{a_2} & \beta_0
	\end{pmatrix}.
\end{align*}
\end{prop}

%%%%%%%%%%%%%%%%%%%%%%%%%%%%%%%%%%%%%%%%%%%%%%%%%%%%%%%%%%%%%%%%%%%%%%%%%%%%%%%%%%%%%%%
%%%%%%%%%%%%%%%%%%%%%%%%%%%%%%%%%%%%%%%%%%%%%%%%%%%%%%%%%%%%%%%%%%%%%%%%%%%%%%%%%%%%%%%
\section{Computing Explicit Reductions}\label{comp-red-chap}

We are now ready to explicitly compute the reductions of irreducible, two-dimensional crystalline representations $V(A)$ of $G_K$ which are associated to the Kisin modules $\Ms(\As)$ from which the Descent Assumptions \ref{descent-assump} have been shown to be satisfied in Section \ref{prep-comp-chap}. This is done in each of the cases detailed in the previous section.

%We begin in Section \ref{mod-rep-sec} by detailing how the descent algorithm of Section \ref{red-alg-chap} produces a Kisin module $\Mf(\Af)$ which is canonically associated to $V(A)$. We also display how to use $\Mf(\Af)$ to compute the reduction of $V(A)$ along with various results that assist us in doing so. Sections \ref{red-large-val-sec}, \ref{red-small-val-sec} and \ref{red-f=2-sec} then us these results to compute reductions in each of the cases from Section \ref{prep-comp-chap}. We then use Section \ref{summary-sec} to summarize the results of the previous sections in a concise manner for the readers reference.

%%%%%%%%%%%%%%%%%%%%%%%%%%%%%%%%%%%%%%%%%%%%%%%%%%%%%%%%%%%%%%%%%%%%%%%%%%%%%%%%%%%%%%%
\subsection{Reductions for the Large Valuation Case}\label{red-large-val-sec}

Let $\Ms(\As)$ be a Kisin module over $S_F$ with labeled heights in the range
\begin{align*}
	p+2&\le k_{0}\le 2p-4& 2&\le k_{i}\le p-3 \hspace{.2cm}\text{\normalfont for all}\hspace{.2cm} 1\le i\le f-1
\end{align*}

\noindent and $\nu_p(a_2^{(0)})\ge 1$ as in Section \ref{prep-large-val-sec}. The results of Proposition \ref{dec-assump-large-val-prop} combined with Theorem \ref{desc-alg-thm} allow use to identify a descent of $\Ms(\As)$ to $\Sig_F$ denoted $\Mf(\Af)$. The reduction of which $\overline{\Mf(\Af)}$ necessarily has partial Frobenius matrices given by
	\begin{align*}
	\overline{\Af}^{(0)}&=\begin{dcases}
			\begin{pmatrix}
				0 & u^{k_{0}}\overline{a}_1^{(0)} \\ 1 & \left(\overline{a_2^{(0)}/p}\right)\overline{\alpha}_1^{(0)}u^p
			\end{pmatrix} & \text{if	} 0\in\Ss\\
			\begin{pmatrix}
				u^{k_{0}}\overline{a}_1^{(0)} & 0 \\ \left(\overline{a_2^{(0)}/p}\right)\overline{\alpha}_1^{(0)}u^p & 1
			\end{pmatrix}& \text{if	} 0\in\Ts 
		\end{dcases}\\
	\overline{\Af}^{(i)}&=\begin{dcases}
		\begin{pmatrix}
			0 & u^{k_i}\overline{a}_1^{(i)} \\ 1 & \overline{a}_2^{(i)}
		\end{pmatrix} & \text{if	} i\neq 0 \text{ and } i \in\Ss\\
		\begin{pmatrix}
			u^{k_i}\overline{a}_1^{(i)} & 0 \\ 0 & 1
		\end{pmatrix}& \text{if	} i\neq 0 \text{ and } i \in\Ts
	\end{dcases}
\end{align*}

\noindent due to the fact that $E=u+p$, $a_2^{(i)}\in\varpi\oh_F$ for all $i\in\Ts$ and the following lemma regarding the reduction of the $\lambda_b$-coefficients $\alpha_0^{(i)}$ and $\alpha_1^{(i)}$.

\begin{lem}\label{lambda-estim-lem}
	Let $\gamma=\varphi(E)/p$ and $\lambda_b^{f(\varphi)}=\prod_{n\ge0}\varphi^{bn}(\gamma)^{f(\varphi)}=\sum_{j\ge0}\alpha_j\left(\frac{E^p}{p}\right)^j$. We have the following estimates:
	\begin{enumerate}
		\item $\gamma\equiv \frac{E^p}{p}+1\pmod p$;
		
		\item $\alpha_0\equiv 1\pmod p$;
		
		\item $\overline{\alpha}_1\equiv f(0)^{\pm1}\pmod p$.
	\end{enumerate}
	\begin{proof}
		$(a)$ follows from the fact that $\frac{E^p}{p}\equiv \frac{u^p}{p}\pmod p$ and that $\gamma=1+\frac{u^p}{p}$. For $(b)$ and $(c)$, write $\lambda_b^{f(\varphi)}=\gamma^{f(0)}\lambda_b^{f(\varphi)-f(0)}$. Here we have $\lambda_b^{f(\varphi)-f(0)}=1+\sum_{r\ge2,\ell\ge1} a_\ell(\frac{u^{rp}}{p})^\ell$ where the sum only contributes to coefficients of $\left(\frac{E^p}{p}\right)^j$ for $j\ge2$. Hence, $\alpha_0$ and $\alpha_1$ only depend on $\gamma^{f(0)}$ and the result follows from $(a)$ after performing a binomial expansion on $\gamma^{f(0)}\equiv(\frac{E^p}{p}+1)^{f(0)}$. 
	\end{proof}
\end{lem}

In preparation, let us define subsets of $\Ss$ by
\[\Vs^=\coloneqq\{i\in\Ss:\nu_p(a_2^{(i)})=0\}\hspace{1cm}\Vs^>\coloneqq\{i\in\Ss:\nu_p(a_2^{(i)})>0\}.\]
That is, we have decomposed the set $\Ss=\Vs^>\sqcup\Vs^=$ so that $\Vs^>$ gives the $i\in\Ss$ with partial Frobenius $\overline{\Af}^{(i)}$ being anti-diagonal. We will repeatedly see that the elements of these sets dictate the shape of the reduction data.

%%%%%%%%%%%%%%%%%%%%%%%%%%%%%%%%%%%%%%%%%%%%%%%%%%%%%%%%%%%%%%%%%%%%%%%%%%%%%%%%%%%%%%%
\subsubsection{Reduction Data for $\nu_p(a_2^{(0)})>1$}\label{verylarge-red-data-subsec}

Let us begin with the case that $\nu_p(a_2^{(0)})>1$ so that we can take $0\in\Vs^>$ if $0\in\Ss$. Hence, we can simplify our reduced partial Frobenius matrices to
\begin{align*}
	\overline{\Af}^{(i)}&=\begin{dcases}
		\begin{pmatrix}
			0 & u^{k_i}\overline{a}_1^{(i)} \\ 1 & \overline{a}_2^{(i)}
		\end{pmatrix} & \text{if	}  i \in\Vs^=\\
		\begin{pmatrix}
			0 & u^{k_i}\overline{a}_1^{(i)} \\ 1 & 0
		\end{pmatrix} & \text{if	}  i \in\Vs^>\\
		\begin{pmatrix}
			u^{k_i}\overline{a}_1^{(i)} & 0 \\ 0 & 1
		\end{pmatrix}& \text{if	}  i \in\Ts.
	\end{dcases}
\end{align*}

\noindent We have two subcases depending on whether or not the set $\Vs^=$ is empty.

\subsubsection*{Subcase: $\Vs^=$ Empty}
If $\Vs^=$ is empty then all the matrices are diagonal or anti-diagonal so our only obstacle to displaying the reduction data is the fact that the powers of $u$ in $\overline{\Af}^{(0)}$ are not smaller than $p$. To fix this, let us apply the base change over $k_F(\!(u)\!)$ given by $X^{(f-1)}=\Diag(1,u)$ if $0\in\Ss$ and $X^{(f-1)}=\Diag(u,1)$ if $0\in\Ts$ and $X^{(i)}=I_2$ everywhere else. The result of $X^{(i)}*_\varphi\overline{\Af}^{(i)}$ will leave $\overline{\Af}^{(i)}$ unaffected for $0<i<f-1$ and provide
	\begin{align*}
		\overline{\Af}^{(0)}&=\begin{dcases}
			\begin{pmatrix}
				0 & u^{k_{0}-p}\overline{a}_1^{(0)} \\ 1 & 0
			\end{pmatrix} & \text{if	} 0\in\Ss\\
			\begin{pmatrix}
				u^{k_{0}-p}\overline{a}_1^{(0)} & 0 \\ 0 & 1
			\end{pmatrix}& \text{if	} 0\in\Ts
		\end{dcases}\\
		\overline{\Af}^{(f-1)}&=\begin{dcases}
			\begin{pmatrix}
				0 & u^{k_{f-1}}\overline{a}_1^{(f-1)} \\ u & 0
			\end{pmatrix} & \text{if	} f-1 \in\Ss \text{ and }0\in\Ss\\
						\begin{pmatrix}
				0 & u^{k_{f-1}+1}\overline{a}_1^{(f-1)} \\ 1 & 0
			\end{pmatrix} & \text{if	} f-1 \in\Ss \text{ and }0\in\Ts\\
			\begin{pmatrix}
				u^{k_{f-1}}\overline{a}_1^{(f-1)} & 0 \\ 0 & u
			\end{pmatrix}& \text{if	} f-1 \in\Ts \text{ and } 0\in\Ss \\
				\begin{pmatrix}
				u^{k_{f-1}+1}\overline{a}_1^{(f-1)} & 0 \\ 0 & 1
			\end{pmatrix}& \text{if	} f-1 \in\Ts \text{ and } 0\in\Ts.
		\end{dcases}
\end{align*}

\noindent Up to an unramified twist or restriction to inertia, the reduction does not depend on the $k_F$-elements $\overline{a}_1^{(i)}$. Hence, in the case that $\nu_p(a_2^{(f-1)})>1$ and $\nu_p(a_2^{(i)})>0$ for all $i<f-1$, we have the reduction data $\mu=(\mu_i)$ where
\begin{align*}
	\mu_{0}&=\begin{dcases}
		(S, (0,k_{0}-p))& \text{if	} 0 \in\Ss\\
		(I, (k_{0}-p, 0)) & \text{if	}  0\in\Ts
	\end{dcases}\\
	\mu_i&=\begin{dcases}
		(S, (0,k_i))& \text{if	} 0<i< f-1 \text{ and } i \in\Ss\\
		(I, (k_i, 0)) & \text{if	} 0<i< f-1 \text{ and } i \in\Ts\\
	\end{dcases}\\
	\mu_{f-1}&=\begin{dcases}
		(S, (1,k_{f-1})) & \text{if	}f-1\in\Ss \text{ and } 0\in \Ss \\
		(S, (0,k_{f-1}+1)) & \text{if	}f-1\in\Ss \text{ and } 0\in \Ts \\
		(I, (k_{f-1},1)) & \text{if	}f-1\in\Ts \text{ and } 0\in \Ss \\
		(I, (k_{f-1}+1,0)) & \text{if	}f-1\in\Ts \text{ and } 0\in \Ts.
	\end{dcases}
\end{align*}

\noindent Notice that the reduction given by Proposition \ref{compute-hom-with-data-prop} will be irreducible if and only if $|\Ss|$ is odd.

\subsubsection*{Subcase: $\Vs^=$ Nonempty}
Now let us return ourselves to the original setting of $\nu_p(a_2^{(0)})>1$ and suppose $\Vs^=$ is nonempty so we have reduced partial Frobenius given by
	\begin{align*}
	\overline{\Af}^{(i)}&=\begin{dcases}
		\begin{pmatrix}
			0 & u^{k_i}\overline{a}_1^{(i)} \\ 1 & \overline{a}_2^{(i)}
		\end{pmatrix} & \text{if	}  i \in\Vs^=\\
			\begin{pmatrix}
			0 & u^{k_i}\overline{a}_1^{(i)} \\ 1 & 0
		\end{pmatrix} & \text{if	}  i \in\Vs^>\\
		\begin{pmatrix}
			u^{k_i}\overline{a}_1^{(i)} & 0 \\ 0 & 1
		\end{pmatrix}& \text{if	}  i \in\Ts.
	\end{dcases}
\end{align*}

We will need to permute rows and columns so as to make each matrix lower triangular in order to perform a semi-simplification. To do this, we define a subset $\Xs\subset\Z/f\Z$ in a similar manner to how we defined $\Ws$ in Definition \ref{sandwich-subset-def}.

\begin{defin}\label{X-sandwich-subset-def}
	We define a subset $\Xs\subset \Z/f\Z$ with the following steps:
	\begin{enumerate}[label*=\arabic*.]
		\item Label the elements of $\Ss=\{s_1,\dots,s_d\}$ and assign $s_j=r_j$ whenever $s_j\in\Vs^>$ so that $\Vs^>=\{r_j\}_j$. If $\Vs^{>}=\emptyset$ then set $\Xs=\emptyset$.
		
		\item Choose any element $s_m\in\Vs^=$. Let $s_\ell=r_\ell$ denote the smallest element of $\Vs^>$ with $s_m<r_\ell$. Include $i\in\Xs$ for all $r_{\ell} \le i< s_{\ell+1}$. 
		
		\item Now pick $r_k$ to be the smallest element of $\Vs^>$ with $s_{\ell+1}<r_k$. Include $i\in\Xs$ for all $r_k\le i<s_{k+1}$.
		
		\item Repeat this process, possibly taking the labels modulo $f$. After a finite amount of steps, we will return to $s_m$ and the process will terminate as $s_m\notin\Vs^>$.
	\end{enumerate}
\end{defin}

We will now permute the rows of $\overline{\Af}^{(i)}$ according to the subset $\Xs$ by making $X^{(i)}=\begin{psmallmatrix}
	0 & 1 \\ 1 & 0
\end{psmallmatrix}$ for all $i\in\Xs$ and $X^{(i)}=I_2$ for all $i\notin\Xs$. The result of $X^{(i)}*_\varphi\overline{\Af}^{(i)}$ will be 
	\begin{align*}
	\overline{\Af}^{(i)}&=\begin{dcases}
		\begin{pmatrix}
		u^{k_i}\overline{a}_1^{(i)} & 0 \\ \overline{a}_2^{(i)} & 1
		\end{pmatrix} & \text{if	}  i \in\Vs^=\text{ and }i-1\in\Xs\\
		\begin{pmatrix}
			0 & u^{k_i}\overline{a}_1^{(i)} \\ 1 & \overline{a}_2^{(i)}
		\end{pmatrix} & \text{if	}  i \in\Vs^=\text{ and }i-1\notin\Xs\\
		\begin{pmatrix}
			1 & 0 \\ 0 & u^{k_i}\overline{a}_1^{(i)}
		\end{pmatrix} & \text{if	}  i \notin\Vs^=\text{ and }i\in\Xs\\
		\begin{pmatrix}
			u^{k_i}\overline{a}_1^{(i)} & 0 \\ 0 & 1
		\end{pmatrix}& \text{if	}  i \notin\Vs^=\text{ and }i\notin\Xs.
	\end{dcases}
\end{align*}

We can now invoke Lemma \ref{straightening-lem} by base changing over $k_F\pu$ via
\[X^{(i)}=\begin{pmatrix}
	1 & -\frac{\overline{a}_1^{(i)}}{\overline{a}_2^{(i)}}u^{k_i} \\ 0 & 1
\end{pmatrix}\in\GL_2(k_F\pu)\]
whenever $i\in\Vs^{=}$ \textit{and} $i-1\notin\Xs$ and setting all other $X^{(i)}$ to be the identity. The result of $X^{(i)}*_\varphi\overline{\Af}^{(i)}$ will be 
\begin{align*}
	\overline{\Af}^{(i)}&=\begin{dcases}
		\begin{pmatrix}
			u^{k_i}\overline{a}_1^{(i)} & 0 \\ \overline{a}_2^{(i)} & 1
		\end{pmatrix} & \text{if	}  i \in\Vs^=\text{ and }i-1\in\Xs\\
		\begin{pmatrix}
			-\frac{\overline{a}_1^{(i)}}{\overline{a}_2^{(i)}}u^{k_i} & 0 \\ 1 & \overline{a}_2^{(i)}
		\end{pmatrix} & \text{if	}  i \in\Vs^=\text{ and }i-1\notin\Xs\\
		\begin{pmatrix}
			1 & 0 \\ 0 & u^{k_i}\overline{a}_1^{(i)}
		\end{pmatrix} & \text{if	}  i \notin\Vs^=\text{ and }i\in\Xs\\
		\begin{pmatrix}
			u^{k_i}\overline{a}_1^{(i)} & 0 \\ 0 & 1
		\end{pmatrix}& \text{if	}  i \notin\Vs^=\text{ and }i\notin\Xs.
	\end{dcases}
\end{align*}

Since $\overline{a}_2^{(i)}\in k_F^*$ for all $i\in\Vs^=$, then we may semi-simplify along the $\varphi$-stable subspace given by $\varphi^{(i)}(e_1^{(i-1)})$ to disregard the lower left entry in the $\overline{\Af}^{(i)}$ corresponding to $i\in\Vs^=$. As before, we may also disregard the $k_F$-elements $\overline{a}_1^{(i)}$ and $\overline{a}_2^{(i)}$ when computing reductions up to unramified twists so that 
\begin{align*}
	\overline{\Af}^{(i)}&=\begin{dcases}
		\Diag\left(1,u^{k_i}\right) & \text{if	}  i\in\Xs\\
		\Diag\left(u^{k_i},1\right)& \text{if	}  i\notin\Xs
	\end{dcases}
\end{align*}

\noindent due to the fact that any element of $\Vs^=$ cannot lie in $\Xs$.

Our final step is the same as the $\Vs^=$ empty case. In particular, we take $X^{(f-1)}=\Diag(1,u)$ if $0\in\Xs$ and $X^{(f-1)}=\Diag(u,1)$ if $0\notin\Xs$ with $X^{(i)}=I_2$ for all $i\neq 0$. The result of $X^{(i)}*_\varphi\overline{\Af}^{(i)}$ will leave $\overline{\Af}^{(i)}$ unaffected for $0<i<f-1$ and give 
\begin{align*}
	\overline{\Af}^{(0)}&=\begin{dcases}
		\Diag\left(1,u^{k_0-p}\right) & \text{if	}  0\in\Xs\\
		\Diag\left(u^{k_0-p},1\right) & \text{if	}  0\notin\Xs
	\end{dcases}\\
		\overline{\Af}^{(f-1)}&=\begin{dcases}
		\Diag\left(1,u^{k_{f-1}+1}\right) & \text{if	}  0,f-1\in\Xs\\
		\Diag\left(u,u^{k_{f-1}}\right) & \text{if	}  f-1\in\Xs\text{ and }0\notin\Xs\\
		\Diag\left(u^{k_{f-1}},u\right)& \text{if	}  f-1\notin\Xs \text{ and }0\in\Xs \\
		\Diag\left(u^{k_{f-1}+1},1\right)& \text{if	}  0,f-1\notin\Xs.
	\end{dcases}
\end{align*}

\noindent Hence, we get reduction data up to restriction to inertia given by $\mu=(\mu_i)$ with
\begin{align*}
	\mu_{0}&=\begin{dcases}
		(I, (0, k_{0}-p)) & \text{if	}  0\in\Xs\\
		(I, (k_{0}-p, 0))& \text{if	}  0\notin\Xs
	\end{dcases}\\
	\mu_i&=\begin{dcases}
		(I, (0, k_i)) & \text{if	}0<i<f-1\text{ and }  i\in\Xs\\
		(I, (k_i,0))& \text{if	}  0<i<f-1\text{ and }i\notin\Xs
	\end{dcases}\\
	\mu_{f-1}&=\begin{dcases}
		(I, (0,k_{f-1}+1)) & \text{if }0,f-1\in\Xs \\
		(I, (1, k_{f-1})) & \text{if }f-1\in\Xs\text{ and }0\notin\Xs\\
		(I, (k_{f-1},1)) & \text{if }f-1\in\Xs\text{ and }0\in\Xs \\
		(I, (k_{f-1}0+1,0)) & \text{if }0,f-1\notin\Xs.
	\end{dcases}
\end{align*}

\noindent Note that the reduction given by Proposition \ref{compute-hom-with-data-prop} will necessarily be reducible.

%%%%%%%%%%%%%%%%%%%%%%%%%%%%%%%%%%%%%%%%%%%%%%%%%%%%%%%%%%%%%%%%%%%%%%%%%%%%%%%%%%%%%%%

\subsubsection{Reduction Data for $\nu_p(a_2^{(0)})=1$}\label{integral-red-data-subsec}
Recall that we have reduced partial Frobenius given by
	\begin{align*}
	\overline{\Af}^{(0)}&=\begin{dcases}
			\begin{pmatrix}
				0 & u^{k_{0}}\overline{a}_1^{(0)} \\ 1 & \left(\overline{a_2^{(0)}/p}\right)\overline{\alpha}_1^{(0)}u^p
			\end{pmatrix} & \text{if	} 0\in\Ss\\
			\begin{pmatrix}
				u^{k_{0}}\overline{a}_1^{(0)} & 0 \\ \left(\overline{a_2^{(0)}/p}\right)\overline{\alpha}_1^{(0)}u^p & 1
			\end{pmatrix}& \text{if	} 0\in\Ts
		\end{dcases}\\
	\overline{\Af}^{(i)}&=\begin{dcases}
		\begin{pmatrix}
			0 & u^{k_i}\overline{a}_1^{(i)} \\ 1 & \overline{a}_2^{(i)}
		\end{pmatrix} & \text{if	} i\neq 0 \text{ and } i \in\Ss\\
		\begin{pmatrix}
			u^{k_i}\overline{a}_1^{(i)} & 0 \\ 0 & 1
		\end{pmatrix}& \text{if	} i\neq 0 \text{ and } i \in\Ts.
	\end{dcases}
\end{align*}

\noindent We will again consider the sets $\Vs^>$ and $\Vs^=$ but will thankfully not have to split into subcases depending on the emptiness of these sets. Unfortunately, we will have to define a new subset $\Ys\subset\Z/f\Z$ that is slightly different from the $\Xs$ we defined in Definition \ref{X-sandwich-subset-def}.

\begin{defin}\label{Y-sandwich-subset-def}
		We define a subset $\Ys\subset \Z/f\Z$ with the following steps:
	\begin{enumerate}[label*=\arabic*.]
		\item Label the elements of $\Ss=\{s_1,\dots,s_d\}$ and assign $s_j=r_j$ whenever $s_j\in\Vs^>$ so that $\Vs^>=\{r_j\}_j$. If $\Vs^{>}=\emptyset$ then set $\Ys=\emptyset$.
		
		\item Let $r_\ell$ denote the smallest element of $\Vs^>$. Include $i\in\Ys$ for all $r_\ell\le i<s_{\ell+1}$.
		
		\item Now pick $r_k$ to be the smallest element of $\Vs^>$ with $s_{\ell+1}<r_k$. Include $i\in\Ys$ for all $r_k\le i<s_{k+1}$.
		
		\item Repeat this process until we run out of elements of $\Vs^{>}$. If $r_d=s_d$ exists, include $r_d\le i\le f-1$.
	\end{enumerate}
\end{defin}

As before, we swap the rows of $\overline{\Af}^{(i)}$ whenever $i\in\Ys$ via $X^{(i)}=\begin{psmallmatrix}
	0 & 1 \\ 1 & 0
\end{psmallmatrix}$ for all $i\in\Ys$ and $X^{(i)}$ being the identity otherwise. For any $i\in\Vs^=$ with $i-1\notin\Ys$, we can invoke Lemma \ref{straightening-lem} by setting
\[X^{(i)}=\begin{pmatrix}
1 & -\frac{\overline{a}_1^{(i)}}{\overline{a}_2^{(i)}}u^{k_i} \\ 0 & 1
\end{pmatrix}\in\GL_2(k_F\pu)\]
whenever $i\in\Vs^{=}$ and $i-1\in\notin\Ys$ and setting all other $X^{(i)}$ to be the identity. The result of all of this will be
	\begin{align*}
	\overline{\Af}^{(0)}&=\begin{dcases}
			\begin{pmatrix}
				0 & u^{k_{0}}\overline{a}_1^{(0)} \\ 1 & \tilde{a}_2^{(0)}\overline{\alpha}_1^{(0)}u^p
			\end{pmatrix} & \text{if	} 0\in\Ss \text{ and }f-1\notin\Ys \text{ or } 0\in\Ts\text{ and }f-1\in\Ys\\
			\begin{pmatrix}
				u^{k_{0}}\overline{a}_1^{(0)} & 0 \\ \tilde{a}_2^{(0)}\overline{\alpha}_1^{(0)}u^p & 1
			\end{pmatrix}& \text{if	} 0\in\Ts \text{ and }f-1\notin\Ys \text{ or }0\in\Ss \text{ and }f-1\in\Ys
		\end{dcases} \\
	\overline{\Af}^{(i)}&=\begin{dcases}
		\begin{pmatrix}
			u^{k_i}\overline{a}_1^{(i)} & 0 \\ \overline{a}_2^{(i)} & 1
		\end{pmatrix} & \text{if	}  i \in\Vs^=\text{ and }i-1\in\Ys\\
		\begin{pmatrix}
			-\frac{\overline{a}_1^{(i)}}{\overline{a}_2^{(i)}}u^{k_i} & 0 \\ 1 & \overline{a}_2^{(i)}
		\end{pmatrix} & \text{if	}  i \in\Vs^=\text{ and }i-1\notin\Ys\\
		\begin{pmatrix}
			1 & 0 \\ 0 & u^{k_i}\overline{a}_1^{(i)}
		\end{pmatrix} & \text{if	}  i \notin\Vs^=\text{ and }i\in\Ys\\
		\begin{pmatrix}
			u^{k_i}\overline{a}_1^{(i)} & 0 \\ 0 & 1
		\end{pmatrix}& \text{if	}  i \notin\Vs^=\text{ and }i\notin\Ys
	\end{dcases}
\end{align*}

\noindent where we begin to write $\tilde{a}_2^{(0)}=(\overline{a_2^{(0)}/p})$. At this point, we are required to diverge into two separate subcases depending on the shape of the $\overline{\Af}^{(0)}$ matrix.

\subsubsection*{Subcase: $0\in\Ss$ and $f-1\notin\Ys$ or $0\in\Ts$ and $f-1\in\Ys$}
We begin by using Lemma \ref{straightening-lem} with $k_{0}-p\ge 4$ to kill the $u^{k_{0}}$-term in $\overline{\Af}^{(0)}$. Indeed, we set
\[X^{(0)}=\begin{pmatrix}
	1 & -\left(\frac{\overline{a}_1}{\tilde{a}_2\overline{\alpha}_1}\right)^{(0)}u^{k_{0}-p} \\ 0 & 1
\end{pmatrix}\in\GL_2(k_F\pu)\]
and $X^{(i)}=I_2$ for all $i\neq 0$. The result of $X^{(i)}*_\varphi\overline{\Af}^{(i)}$ will leave $\overline{\Af}^{(i)}$ unaffected for $1\le i\le f-1$ and give
	\begin{align*}
	\overline{\Af}^{(0)}&=
		\begin{pmatrix}
			-\left(\frac{\overline{a}_1}{\tilde{a}_2\overline{\alpha}_1}\right)^{(0)}u^{k_{0}-p} & 0 \\ 1 & \tilde{a}_2^{(0)}\overline{\alpha}_1^{(0)}u^p
		\end{pmatrix}.
\end{align*}

\noindent Hence, we may semi-simplify to disregard the lower left entries of each matrix. Our last operation will be to scale by $u$ via $X^{(f-1)}=\Diag(1,u)$ so that up restriction to inertia, our reduction data will be given by $\mu=(\mu_i)$ with
\begin{align*}
		\mu_{0}&=(I, (k_{0}-p, 0)).\\
		\mu_i&=\begin{dcases}
			(I, (k_i, 0)) & \text{if }1\le i<f-1 \text{ and }i\notin\Ys \\
			(I, (0, k_i)) & \text{if }1\le i<f-1 \text{ and }i\in \Ys
		\end{dcases}\\
		\mu_{f-1}&=\begin{dcases}
		(I, (k_{f-1}, 1)) & \text{if }f-1\notin\Ys \\
		(I, (0, k_{f-1}+1)) & \text{if }f-1\in \Ys.
	\end{dcases}
\end{align*}

\noindent Note that this case will necessarily result in a reducible reduction by Proposition \ref{compute-hom-with-data-prop}.

\subsubsection*{Subcase: $0\in\Ts$ and $f-1\notin\Ys$ or $0\in\Ss$ and $f-1\in\Ys$}
Recall that we have reduced partial Frobenius given by 
	\begin{align*}
	\overline{\Af}^{(0)}&=
		\begin{pmatrix}
			u^{k_{0}}\overline{a}_1^{(0)} & 0 \\ \tilde{a}_2^{(0)}\overline{\alpha}_1^{(0)}u^p & 1
		\end{pmatrix}\\
	\overline{\Af}^{(i)}&=\begin{dcases}
		\begin{pmatrix}
			u^{k_i}\overline{a}_1^{(i)} & 0 \\ \overline{a}_2^{(i)} & 1
		\end{pmatrix} & \text{if	}  i \in\Vs^=\text{ and }i-1\in\Ys\\
		\begin{pmatrix}
			-\frac{\overline{a}_1^{(i)}}{\overline{a}_2^{(i)}}u^{k_i} & 0 \\ 1 & \overline{a}_2^{(i)}
		\end{pmatrix} & \text{if	}  i \in\Vs^=\text{ and }i-1\notin\Ys\\
		\begin{pmatrix}
			1 & 0 \\ 0 & u^{k_i}\overline{a}_1^{(i)}
		\end{pmatrix} & \text{if	}  i \notin\Vs^=\text{ and }i\in\Ys\\
		\begin{pmatrix}
			u^{k_i}\overline{a}_1^{(i)} & 0 \\ 0 & 1
		\end{pmatrix}& \text{if	}  i \notin\Vs^=\text{ and }i\notin\Ys.
	\end{dcases}
\end{align*}

\noindent Let us scale by $u$ via $X^{(f-1)}=\Diag(u,1)\in\GL_2(k_F(\!(u)\!))$ and $X^{(i)}=I_2$ for $i\neq f-1$ so that $X^{(i)}*_\varphi\overline{\Af}^{(i)}$ leaves $\overline{\Af}^{(i)}$ unaffected for $1\le i<f-1$ and 
	\begin{align*}
	\overline{\Af}^{(0)}&=
		\begin{pmatrix}
			u^{k_{0}-p}\overline{a}_1^{(0)} & 0 \\ \tilde{a}_2^{(0)}\overline{\alpha}_1^{(0)} & 1
		\end{pmatrix}\\
	\overline{\Af}^{(f-1)}&=\begin{dcases}
		\begin{pmatrix}
			u^{k_{f-1}+1}\overline{a}_1^{(f-1)} & 0 \\ \overline{a}_2^{(f-1)} & 1
		\end{pmatrix} & \text{if	}  f-1 \in\Vs^=\text{ and }f-2\in\Ys\\
		\begin{pmatrix}
			-\frac{\overline{a}_1^{(f-1)}}{\overline{a}_2^{(f-10)}}u^{k_{f-1}+1} & 0 \\ 1 & \overline{a}_2^{(f-1)}
		\end{pmatrix} & \text{if	}  f-1 \in\Vs^=\text{ and }f-2\notin\Ys\\
		\begin{pmatrix}
			u & 0 \\ 0 & u^{k_{f-1}}\overline{a}_1^{(f-1)}
		\end{pmatrix} & \text{if	}  f-1 \notin\Vs^=\text{ and }f-1\in\Ys\\
		\begin{pmatrix}
			u^{k_{f-1}+1}\overline{a}_1^{(f-1)} & 0 \\ 0 & 1
		\end{pmatrix}& \text{if	}  f-1 \notin\Vs^=\text{ and }f-1\notin\Ys.
	\end{dcases}
\end{align*}

\noindent Hence, up to semi-simplification and restriction to inertia, we will have reduction data given by
\begin{align*}
		\mu_{0}&=(I, (k_{0}-p, 0))\\
		\mu_i&=\begin{dcases}
			(I, (k_i, 0)) & \text{if }1\le i<f-1 \text{ and }i\notin\Ys \\
			(I, (0, k_i)) & \text{if }1\le i<f-1 \text{ and }i\in \Ys
		\end{dcases}\\
		\mu_{f-1}&=\begin{dcases}
		(I, (k_{f-1}+1, 0)) & \text{if }f-1\notin\Ys \\
		(I, (1, k_{f-1})) & \text{if }f-1\in \Ys.
	\end{dcases}
\end{align*}

\noindent Note that this case will necessarily result in a reducible reduction by Proposition \ref{compute-hom-with-data-prop}.

%%%%%%%%%%%%%%%%%%%%%%%%%%%%%%%%%%%%%%%%%%%%%%%%%%%%%%%%%%%%%%%%%%%%%%%%%%%%%%%%%%%%%%%
\subsection{Reductions for the Small Valuation Case}\label{red-small-val-sec}

Let $\Ms(\As)$ be a Kisin module over $S_F$ with labeled heights in the range
\begin{align*}
	p+2&\le k_{0}\le 2p-4& 2&\le k_{i}\le p-3 \hspace{.2cm}\text{\normalfont for all}\hspace{.2cm} 1\le i\le f-1
\end{align*}
and $0\le\nu_p(a_2^{(0)})< 1$ as in Section \ref{prep-small-val-sec}. We recall that we encountered a dichotomy in behavior with respect to the even or oddness of the number of Type $\I$ matrices present. Indeed, we will need to take these two cases separately.

\subsubsection{Reduction Data for $|\Ss|$ Even}\label{red-data-even-subseec}
Under the assumption that $\nu_p(a_2^{(i)})\ge 1-\nu_p(a_2^{(0)})$ for all $i\in\Ws$, the results of Proposition \ref{dec-assump-even-prop} combined with Theorem \ref{desc-alg-thm} allow use to identify
a descent of $\Ms(\As)$ to $\Sig_F$ denoted $\Mf(\Af)$. The reduction of which $\overline{\Mf(\Af)}$ necessarily has partial Frobenius matrices given by
\begin{align*}
	\overline{\Af}^{(0)}&=
	\begin{pmatrix}
		0 & -(\overline{a_1/\alpha_1})^{(0)}u^{k_{0}-p} \\ \overline{\alpha}_1^{(0)}u^p & 1
	\end{pmatrix}\\
	\overline{\Af}^{(i)}&=\begin{dcases}
		\begin{pmatrix}
			1 & (\overline{a_2^{(0)}a_2^{(i)}/p}) \\ 0 & u^{k_i}\overline{a}_1^{(i)}
		\end{pmatrix} & \text{if	} 1\le i\le f-1 \text{ and } i\in\Ws\\
		\begin{pmatrix}
			u^{k_i}\overline{a}_1^{(i)} & 0\\ 0 & 1
		\end{pmatrix}& \text{if	} 1\le i \le f-1 \text{ and } i\notin\Ws
	\end{dcases} 
\end{align*}

\noindent since $p/a_2^{(0)}\in\varpi\oh_F$ and applying Lemma \ref{lambda-estim-lem}. As before, let us redefine a subset $\Vs\subset\Ws$ by  
\[\Vs=\{i\in\Ws: \nu_p(a_2^{(i)})= 1-\nu_p(a_2^{(0)})\}.\]
That is, the set of $i\in\Ws$ such that the upper right hand entry of $\overline{\Af}^{(i)}$ is nonzero. Just as before, we have two subcases depending on whether or not the set $\Vs$ is empty.

\subsubsection*{Subcase: $\Vs$ Empty}
In the case that the set $\Vs$ is empty then our only obstruction to displaying the reduction data is $\overline{\Af}^{(0)}$. Since $k_{0}-p\ge 2$ then we may evoke Lemma \ref{straightening-lem} with 
\[X^{(0)}=\begin{pmatrix}
	1 & (\overline{a_1/\alpha_1})^{(0)}u^{k_{0}-p} \\ 0 & 1
\end{pmatrix}\in\GL_2(k_F\pu)\]
and $X^{(i)}=I_2$ for all $i\neq 0$. The result of $X^{(i)}*_\varphi\overline{\Af}^{(i)}$ will leave $\overline{\Af}^{(i)}$ unaffected for $1\le i\le f-1$ and give
\begin{align*}
	\overline{\Af}^{(0)}&=
	\begin{pmatrix}
		\overline{a}_1^{(0)}u^{k_{0}} & 0  \\ \overline{\alpha}_1^{(0)}u^p & 1
	\end{pmatrix}.
\end{align*}

\noindent We will now scale by $u$ via $X^{(f-1)}=\Diag(u,1)$ and $X^{(i)}=I_2$ for $i\neq f-1$ so that the result of $X^{(i)}*_\varphi\overline{\Af}^{(i)}$ will leave $\overline{\Af}^{(i)}$ unaffected for $0<i<f-1$ and give
\begin{align*}
	\overline{\Af}^{(0)}&=
	\begin{pmatrix}
		\overline{a}_1^{(0)}u^{k_{0}-p} & 0  \\ \overline{\alpha}_1^{(0)} & 1
	\end{pmatrix}\\
		\overline{\Af}^{(f-1)}&=\begin{dcases}
		\begin{pmatrix}
			u & 0 \\ 0 & u^{k_{f-1}}\overline{a}_1^{(f-1)}
		\end{pmatrix} & \text{if	}  f-1\in\Ws\\
		\begin{pmatrix}
			u^{k_{f-1}+1}\overline{a}_1^{(f-1)} & 0\\ 0 & 1
		\end{pmatrix}& \text{if	} f-1\notin\Ws.
	\end{dcases}
\end{align*}

\noindent Since $\overline{\alpha}_1^{(0)}\in k_F^*$ by Lemma \ref{lambda-estim-lem}, then up to semi-simplification, the lower right hand entry of $\overline{\Af}^{(0)}$ may be disregarded to get reduction data given by $\mu=(\mu_i)$ with
\begin{align*}
	\mu_{0}&=(I, (k_{0}-p,0))\\
	\mu_i&=\begin{dcases}
			(I, (0, k_i)) & \text{if	} 0<i< f-1 \text{ and }i\in\Ws \\
			(I, (k_i,0))&\text{if	} 0<i< f-1 \text{ and }i\notin\Ws 
		\end{dcases}\\
	\mu_{f-1}&=\begin{dcases}
		(I, (1, k_{f-1})) & \text{if	} f-1\in\Ws \\
		(I, (k_{f-1}+1,0)) &\text{if	} f-1\notin\Ws.
	\end{dcases}
\end{align*}

\noindent Note that this case will necessarily result in a reducible reduction by Proposition \ref{compute-hom-with-data-prop}.

\subsubsection*{Subcase: $\Vs$ Nonempty}
Let us return ourselves to the original setting with the assumption that the set $\Vs$ is nonempty,
\begin{align*}
	\overline{\Af}^{(0)}&=
	\begin{pmatrix}
		0 & -(\overline{a_1/\alpha_1})^{(0)}u^{k_{0}-p} \\ \overline{\alpha}_1^{(0)}u^p & 1
	\end{pmatrix} \\
	\overline{\Af}^{(i)}&=\begin{dcases}
		\begin{pmatrix}
			1 & \tilde{a}_2^{(i)} \\ 0 & u^{k_i}\overline{a}_1^{(i)}
		\end{pmatrix} & \text{if	} 1\le i\le f-1 \text{ and } i\in\Vs\\
		\begin{pmatrix}
			1 & 0 \\ 0 & u^{k_i}\overline{a}_1^{(i)}
		\end{pmatrix} & \text{if	} 1\le i\le f-1 \text{ and } i\notin\Vs \text{ but }i\in\Ws\\
		\begin{pmatrix}
			u^{k_i}\overline{a}_1^{(i)} & 0\\ 0 & 1
		\end{pmatrix}& \text{if	} 1\le i\le f-1 \text{ and } i\notin\Ws
	\end{dcases}
\end{align*}

\noindent where we write $\tilde{a}_2^{(i)}=(\overline{a_2^{(0)}a_2^{(i)}/p})$. Let $r$ denote the minimal element in $\Vs$. Since $k_r\ge 2$ and $\tilde{a}_2^{(r)}\in k_F^*$ then we may evoke Lemma \ref{straightening-lem} with 
\[X^{(r)}=\begin{pmatrix}
	1 & 0 \\ -(\overline{a_1/\tilde{a}_2})^{(r)}u^{k_r} & 1
\end{pmatrix}\in\GL_2(k_F\pu)\]
and $X^{(i)}=I_2$ for $i\neq r$. The result will leave all $\overline{\Af}^{(i)}$ the same except for $i=r$ where we have
\[\overline{\Af}^{(r)}=\begin{pmatrix}
	1 & \tilde{a}_2^{(r)} \\ -(\overline{a_1/\tilde{a}_2})^{(r)}u^{k_r} & 0
\end{pmatrix}.\]
Now let us swap rows and columns by $X^{(i)}=\begin{psmallmatrix}
	0 & 1 \\ 1 & 0
\end{psmallmatrix}$ for all $0\le i<r$ and $X^{(j)}=I_2$ for all $r\le j\le f-1$. The result will be
\begin{align*}
	\overline{\Af}^{(0)}&=
	\begin{pmatrix}
		\overline{\alpha}_1^{(0)}u^p & 1 \\ 0 & -(\overline{a_1/\alpha_1})^{(0)}u^{k_{0}-p}
	\end{pmatrix}\\
	\overline{\Af}^{(j)}&=\begin{dcases}
		\begin{pmatrix}
			u^{k_j}\overline{a}_1^{(j)} & 0 \\ 0 & 1
		\end{pmatrix} & \text{if	} 0<j<r \text{ and } j\in\Ws\\
		\begin{pmatrix}
			1 & 0\\ 0 & u^{k_j}\overline{a}_1^{(j)}
		\end{pmatrix}& \text{if	} 0<j<r \text{ and } j\notin\Ws
	\end{dcases}\\
		\overline{\Af}^{(r)}&=\begin{pmatrix}
		\tilde{a}_2^{(r)} & 1 \\ 0 & -(\overline{a_1/\tilde{a}_2})^{(r)}u^{k_r}
	\end{pmatrix}\\
	\overline{\Af}^{(i)}&=\begin{dcases}
		\begin{pmatrix}
			1 & \tilde{a}_2^{(i)} \\ 0 & u^{k_i}\overline{a}_1^{(i)}
		\end{pmatrix} & \text{if	} r< i\le f-1 \text{ and } i\in\Vs\\
		\begin{pmatrix}
			1 & 0 \\ 0 & u^{k_i}\overline{a}_1^{(i)}
		\end{pmatrix} & \text{if	} r< i\le f-1 \text{ and } i\notin\Vs\text{ but }i\in\Ws\\
		\begin{pmatrix}
			u^{k_i}\overline{a}_1^{(i)} & 0\\ 0 & 1
		\end{pmatrix}& \text{if	} r< i\le f-1 \text{ and } i\notin\Ws.
	\end{dcases}
\end{align*}

\noindent Since $\tilde{a}_2^{(i)}\in k_F^*$ for all $i\in\Vs$ then up to semi-simplification we may disregard all upper right entries. Our last operation is to scale by $u$ via $X^{(f-1)}=\Diag(u,1)$ so we get reduction data given by $\mu=(\mu_i)$ with
\begin{align*}
	\mu_{0}&=(I, (0,k_{0}-p))\\
	\mu_j&=\begin{dcases}
		(I, (k_j, 0))& \text{if }0<j<r \text{ and }j\in\Ws \\
		(I, (0,k_j))& \text{if }0<j<r \text{ and }j\notin\Ws
	\end{dcases}\\
	\mu_i&=\begin{dcases}
		(I, (0, k_i)) & \text{if }r\le i<f-1 \text{ and }i\in\Ws \\
		(I, (k_i, 0))& \text{if	}r\le i<f-1 \text{ and }i\notin\Ws
	\end{dcases}\\
	\mu_{f-1}&=\begin{dcases}
		(I, (1, k_{f-1})) & \text{if }f-1\in\Ws \\
		(I, (k_{f-1}+1, 0))& \text{if	}f-1\notin\Ws.
	\end{dcases}
\end{align*}

\noindent Note that this case will necessarily result in a reducible reduction by Proposition \ref{compute-hom-with-data-prop}.

\subsubsection{Reduction Data for $|\Ss|$ Odd}\label{red-data-odd-subsec}
Under the assumption that $\nu_p(a_2^{(i)})\ge 1$ for all $i\in\Ws$, the results of Proposition \ref{dec-assump-odd-prop} combined with Theorem \ref{desc-alg-thm} allow use to identify
a descent of $\Ms(\As)$ to $\Sig_F$ denoted $\Mf(\Af)$. The reduction of which $\overline{\Mf(\Af)}$ necessarily has partial Frobenius matrices given by
	\begin{align*}
	\overline{\Af}^{(0)}&=
		\begin{pmatrix}
			0 & -(\overline{a_1/\alpha_1})^{(0)}u^{k_{0}-p} \\ -\overline{\alpha}_1^{(0)}u^p & \overline{a}_2^{(0)}
		\end{pmatrix}\\
	\overline{\Af}^{(i)}&=\begin{dcases}
		\begin{pmatrix}
			1 & (\overline{a_2^{(i)}/p}) \\ 0 & u^{k_i}\overline{a}_1^{(i)}
		\end{pmatrix} & \text{if	} 1\le i \le f-1 \text{ and }i\in\Ws\\
		\begin{pmatrix}
			u^{k_i}\overline{a}_1^{(i)} & 0\\  0 & 1
		\end{pmatrix}& \text{if	} 1\le i \le f-1 \text{ and }i\notin\Ws
	\end{dcases}
\end{align*}

\noindent since $p/a_2^{(0)}\in\varpi\oh_F$ and applying Lemma \ref{lambda-estim-lem}.

We follow very similar ideas as in Section \ref{red-data-even-subseec} with the only exception being the possibility that $\nu_p(a_2^{(0)})>0$ so that $\overline{\Af^{(0)}}$ is anti-diagonal. For this reason, we will leave out the details in the completely similar cases. As before, let us redefine a subset $\Vs\subset\Ws$ by  
\[\Vs=\{i\in\Ws: \nu_p(a_2^{(i)})= 1\}.\]

\subsubsection*{Subcase: $\Vs$ Empty and $\nu_p(a_2^{(0)})>0$}
Here we have each reduced partial Frobenius matrix is either diagonal or anti-diagonal so we can write down the reduction data right away after scaling by $u$ via $X^{(f-1)}=\Diag(u,1)$ so that 
\begin{align*}
	\mu_{0}&=(S, (0 ,k_{0}-p))\\
	\mu_i&=\begin{dcases}
		(I, (0, k_i)) & \text{if } 0<i< f-1 \text{ and }i\in\Ws\\
		(I, (k_i, 0)) & \text{if } 0<i< f-1 \text{ and }i\notin\Ws
	\end{dcases}\\
		\mu_{f-1}&=\begin{dcases}
		(I, (1, k_{f-1})) & \text{if	} f-1\in\Ws\\
		(I, (k_{f-1}+1, 0)) & \text{if	} f-1\notin\Ws.
	\end{dcases}
\end{align*}

\noindent Note that Proposition \ref{compute-hom-with-data-prop} will necessarily provide an irreducible reduction from this data.

\subsubsection*{Subcase: $\Vs$ Empty and $\nu_p(a_2^{(0)})=0$}
This case proceeds precisely as in the $\Vs$ empty case of Section \ref{red-data-even-subseec} after making the necessary coefficient adjustments to account for changing $1$ into $\overline{a}_2^{(0)}$ and $(\overline{a_2^{(0)}a_2^{(i)}/p})$ into $(\overline{a_2^{(i)}/p})$. Up to restriction to inertia, these coefficients do not change the reduction so we may copy the reduction data exactly as it is displayed in that case.
\begin{align*}
	\mu_{0}&=(I, (k_{0}-p,0))\\
	\mu_i&=\begin{dcases}
		(I, (0, k_i)) & \text{if	} 0<i< f-1 \text{ and }i\in\Ws \\
		(I, (k_i,0))&\text{if	} 0<i< f-1 \text{ and }i\notin\Ws 
	\end{dcases}\\
	\mu_{f-1}&=\begin{dcases}
		(I, (1, k_{f-1})) & \text{if	} f-1\in\Ws \\
		(I, (k_{f-1}+1,0)) &\text{if	} f-1\notin\Ws.
	\end{dcases}
\end{align*}

\noindent Note that Proposition \ref{compute-hom-with-data-prop} will necessarily provide a reducible reduction from this data.

\subsubsection*{Subcase: $\Vs$ Nonempty}
Just as before, this case is completely similar to the $\Vs$ nonempty case of Section \ref{red-data-even-subseec} after making the necessary adjustments to various coefficients. Hence, up to restriction to inertia we have our reduction data given by 
\begin{align*}
	\mu_{0}&=(I, (0,k_{0}-p))\\
	\mu_j&=\begin{dcases}
		(I, (k_j, 0))& \text{if }0<j<r \text{ and }j\in\Ws \\
		(I, (0,k_j))& \text{if }0<j<r \text{ and }j\notin\Ws
	\end{dcases}\\
	\mu_i&=\begin{dcases}
		(I, (0, k_i)) & \text{if }r\le i<f-1 \text{ and }i\in\Ws \\
		(I, (k_i, 0))& \text{if	}r\le i<f-1 \text{ and }i\notin\Ws
	\end{dcases}\\
	\mu_{f-1}&=\begin{dcases}
		(I, (1, k_{f-1})) & \text{if }f-1\in\Ws \\
		(I, (k_{f-1}+1, 0))& \text{if	}f-1\notin\Ws.
	\end{dcases}
\end{align*}

\noindent where $r\in\Ws$ is the largest integer with $\nu_p(a_2^{(i)})=1$.
Note that Proposition \ref{compute-hom-with-data-prop} will necessarily provide a reducible reduction from this data.

%%%%%%%%%%%%%%%%%%%%%%%%%%%%%%%%%%%%%%%%%%%%%%%%%%%%%%%%%%%%%%%%%%%%%%%%%%%%%%%%%%%%%%% 
\subsection{\texorpdfstring{Explicit Reductions for the Complete $f=2$ Case}{Explicit Reductions for the Complete f=2 Case}}\label{red-f=2-sec}

As opposed to the general $f>2$ case, the simplicity of the $f=2$ case allows us to explicitly compute reductions in all cases as well as describe the unramified characters resulting from the $k_F$-elements in the determinant.

Using the notation from Section \ref{complete-f=2-sec}, we congregate all the results of Theorem \ref{desc-alg-thm} and Propositions \ref{dec-assump-large-val-prop}, \ref{dec-assump-even-prop}, \ref{dec-assump-odd-prop}, \ref{dec-assump-I,I-prop}, and \ref{dec-assump-II,I-prop} in terms of the three irreducible Type combinations $(\I,\I)$, $(\I,\II)$ and $(\II,\I)$. 
	\begin{figure}[H]
	\centering
	\begin{tabular}{|c|l|l|}
		\hline
		Type & Conditions on valuations & Reduced Frobenius $(\overline{\Af},\overline{\Bf})$ \\
		\hline
		$(\I,\I)$ & 
		$\nu_p(a_2)\ge 1$ & 
		$\begin{pmatrix}
			0 & \overline{a}_1u^{k_0} \\ 1 & (\overline{a_2/p})\overline{\alpha}_1 u^p
		\end{pmatrix},\begin{pmatrix}
			0 & \overline{b}_1u^{k_1} \\ 1 & \overline{b}_2
		\end{pmatrix}$ \\
		
		\hline
		$(\I,\I)$ & 
		$\nu_p(a_2)<1$ and $\nu_p(b_2)\ge 1-\nu_p(a_2)$ & 
		$\begin{pmatrix}
			0 & -(\overline{a_1/\alpha_1})u^{k_0-p} \\ \overline{\alpha}_1 u^p & 1
		\end{pmatrix},\begin{pmatrix}
			1 & (\overline{a_2b_2/p}) \\ 0 & \overline{b}_1u^{k_1}
		\end{pmatrix}$\\
		
		\hline
		$(\I,\I)$ & 
		$\nu_p(a_2)<1$ and $\nu_p(b_2)< 1-\nu_p(a_2)$ & 
		$\begin{pmatrix}
			0 & -(\overline{a_1/\alpha_1})u^{k_0-p} \\ -\overline{\alpha}_1 u^p & 0
		\end{pmatrix},\begin{pmatrix}
			-\overline{b}_1u^{k_1} & 0 \\ 0 & 1
		\end{pmatrix}$\\
		
		\hline
		$(\II,\I)$ & 
		$\nu_p(a_2)\ge 1$ & 
		$\begin{pmatrix}
			\overline{a}_1u^{k_0} & 0 \\ (\overline{a_2/p})\overline{\alpha}_1 u^p & 1
		\end{pmatrix},\begin{pmatrix}
			0 & \overline{b}_1u^{k_1} \\ 1 & 0
		\end{pmatrix}$\\
		
		\hline
		$(\II,\I)$ & 
		$\nu_p(a_2)< 1$ and $\nu_p(b_2)\ge 1$ & 
		$\begin{pmatrix}
			0 & -(\overline{a_1/\alpha_1})u^{k_0-p} \\ -\overline{\alpha}_1u^p & 0
		\end{pmatrix},\begin{pmatrix}
			1 & (\overline{b_2/p}) \\ 0 & \overline{b}_1 u^{k_1}
		\end{pmatrix}$\\
		
		\hline
		$(\II,\I)$ & 
		$\nu_p(a_2)< 1$ and $\nu_p(b_2)< 1$ & 
		$\begin{pmatrix}
			0 & -(\overline{a_1/\alpha_1})u^{k_0-p} \\ \overline{\alpha}_1 u^p & 0
		\end{pmatrix},\begin{pmatrix}
			-\overline{b}_1 u^{k_1} & 0 \\ 0 & 1
		\end{pmatrix}$\\
		
		\hline
		$(\I,\II)$ & 
		$\nu_p(a_2)\ge 1$ & 
		$\begin{pmatrix}
			0 & \overline{a}_1u^{k_0} \\ 1 & (\overline{a_2/p})\overline{\alpha}_1 u^p
		\end{pmatrix},\begin{pmatrix}
			\overline{b}_1u^{k_1} & 0 \\ 0 & 1
		\end{pmatrix}$\\
		
		\hline
		$(\I,\II)$ & 
		$\nu_p(a_2)< 1$ & 
		$\begin{pmatrix}
			0 & -(\overline{a_1/\alpha_1})u^{k_0-p} \\ -\overline{\alpha}_1 u^p & 0
		\end{pmatrix},\begin{pmatrix}
			\overline{b}_1u^{k_1} & 0 \\ 0 & 1
		\end{pmatrix}$\\
		\hline
	\end{tabular}
	
	\caption{Frobenius in $f=2$ over $k_F[\![u]\!]$ depending on Type and valuations.}\label{red-frob-f=2-table}
\end{figure}
Let us recall that by Proposition \ref{irred-vals-prop}, when we are of Type $(\I,\I)$, the product $a_2b_2\in\varpi\oh_F$ and in the other two Type combinations, we have both $a_2,b_2\in\varpi\oh_F$. Let $\overline{\As}\coloneqq\As\pmod\varpi$ and $\overline{\Bs}\coloneqq\Bs\pmod\varpi$. Let us make the observation that all cases where $\nu_p(b_2)$ and $\nu_p(a_2)$ are less than one have the same shape. 
	
For the rest of this section, let us denote $\chi(c)$ to be the unramified character sending Frobenius to $c\in\overline{\F}_p$. 

\subsubsection{Explicit Reductions for Type $(\I,\I)$}
We begin by making the simplifying observation that the five cases in the Type $(\I,\I)$ can be simplified to three cases depending on $\nu_p(a_2b_2)$ when computing reductions. Indeed, we will suppose each subcase below and show that the cases are equivalent.

\subsubsection*{Subcase: $\nu_p(a_2b_2)>1$}
The bound $\nu_p(a_2b_2)>1$ is satisfied three times in Figure \ref{red-frob-f=2-table}.
\begin{enumerate}
	\item $\nu_p(a_2)>1$ and $\nu_p(b_2)\ge 0$:
	\begin{align*}
		\overline{\Af}&=\begin{pmatrix}
			0 & \overline{a}_1u^{k_0} \\ 1 & 0
		\end{pmatrix} & \overline{\Bf}&=\begin{pmatrix}
			0 & \overline{b}_1u^{k_1} \\ 1 & \overline{b}_2
		\end{pmatrix};
	\end{align*}
	
	\item $\nu_p(a_2)=1$ and $\nu_p(b_2)>0$:
	\begin{align*}
		\overline{\Af}&=\begin{pmatrix}
			0 & \overline{a}_1u^{k_0} \\ 1 & (\overline{a_2/p})\overline{\alpha}_1 u^p
		\end{pmatrix} & \overline{\Bf}&=\begin{pmatrix}
			0 & \overline{b}_1u^{k_1} \\ 1 & 0
		\end{pmatrix};
	\end{align*}
	
	\item $\nu_p(a_2)<1$ and $\nu_p(b_2)>1-\nu_p(a_2)$:
	\begin{align*}
		\overline{\Af}&= \begin{pmatrix}
			0 & -(\overline{a_1/\alpha_1})u^{k_0-p} \\ \overline{\alpha}_1 u^p & 1
		\end{pmatrix} & \overline{\Bf}^{(i)}&=\begin{pmatrix}
			1 & 0 \\ 0 & \overline{b}_1u^{k_1}
		\end{pmatrix}.
	\end{align*}
\end{enumerate}

\noindent In (a) we follow the the methods of Section \ref{verylarge-red-data-subsec}, in (b) we follow Section \ref{integral-red-data-subsec} and for (c) we utilize Section \ref{red-data-even-subseec}. The result will be the pair of reduced matrices in each case given by
\begin{enumerate}
	\item $\nu_p(a_2)>$ and $\nu_p(b_2)\ge 0$:
	\begin{align*}
		\overline{\Af}&=\begin{pmatrix}
			1 & 0 \\ 0 & \overline{a}_1u^{k_0-p}
		\end{pmatrix} & 
	\overline{\Bf}&=\begin{pmatrix}
		\overline{b}_1u^{k_1} & 0 \\ 0 & u
	\end{pmatrix};
	\end{align*}

	\item $\nu_p(a_2)=1$ and $\nu_p(b_2)>0$:
	\begin{align*}
		\overline{\Af}&=\begin{pmatrix}
			\overline{a}_1u^{k_0-p} & 0 \\ 0 & 1
		\end{pmatrix} & \overline{\Bf}&=\begin{pmatrix}
			u & 0 \\ 0 & \overline{b}_1u^{k_1}
		\end{pmatrix};
	\end{align*}
	
	\item $\nu_p(a_2)<1$ and $\nu_p(b_2)>1-\nu_p(a_2)$:
	\begin{align*}
	\overline{\Af}&=\begin{pmatrix}
		\overline{a}_1u^{k_0-p} & 0 \\ 0 & 1
	\end{pmatrix} & \overline{\Bf}&=\begin{pmatrix}
		u & 0 \\ 0 & \overline{b}_1u^{k_1}
	\end{pmatrix};
	\end{align*}
\end{enumerate}

\noindent It is then easy to see that all three cases above will result in a product 
\[\overline{\Af}\varphi(\overline{\Bf})=\begin{pmatrix}
	\overline{a}_1u^{k_0} & 0 \\ 0 & \overline{b}_1u^{pk_1}
\end{pmatrix}.\]
Hence, we can use Proposition \ref{compute-Hom-prop} to compute the reduction
\[\overline{V(A,B)}|_{G_\infty}=\chi(\overline{a}_1)\omega_2^{k_0}\oplus \chi(\overline{b}_1)\omega_2^{pk_1}.\]

\subsubsection*{Subcase: $\nu_p(a_2b_2)=1$}
The bound $\nu_p(a_2b_2)=1$ is satisfied twice in Figure \ref{red-frob-f=2-table}.
\begin{enumerate}
	\item $\nu_p(a_2)=1$ and $\nu_p(b_2)=0$:
	\begin{align*}
		\overline{\Af}&=\begin{pmatrix}
			0 & \overline{a}_1u^{k_0} \\ 1 & (\overline{a_2/p})\overline{\alpha}_1 u^p
		\end{pmatrix} & \overline{\Bf}&=\begin{pmatrix}
			0 & \overline{b}_1u^{k_1} \\ 1 & \overline{b}_2
		\end{pmatrix};
	\end{align*}
	
	\item $\nu_p(a_2)<1$ and $\nu_p(b_2)=1-\nu_p(a_2)$:
	\begin{align*}
		\overline{\Af}&=\begin{pmatrix}
			0 & -(\overline{a_1/\alpha_1})u^{k_0-p} \\ \overline{\alpha}_1 u^p & 1
		\end{pmatrix} & \overline{\Bf}&=\begin{pmatrix}
			1 & (\overline{a_2b_2/p}) \\ 0 & \overline{b}_1u^{k_1}
		\end{pmatrix}.
	\end{align*}
\end{enumerate}

\noindent In (a) we follow Section \ref{integral-red-data-subsec} and for (b) we follow Section \ref{red-data-even-subseec}. The result will be
\begin{enumerate}
	\item $\nu_p(a_2)=1$ and $\nu_p(b_2)=0$:
	\begin{align*}
		\overline{\Af}&=\begin{pmatrix}
			-\frac{\overline{a}_1}{(\overline{a_2/p})\overline{\alpha}_1}u^{k_0-p} & 0 \\ 0 & (\overline{a_2/p})\overline{\alpha}_1 
		\end{pmatrix} & \overline{\Bf}&=\begin{pmatrix}
			-(\overline{b_1/b_2})u^{k_1} & 0 \\ 0 & \overline{b}_2u
		\end{pmatrix};
	\end{align*}
	
	\item $\nu_p(a_2)<1$ and $\nu_p(b_2)=1-\nu_p(a_2)$:
	\begin{align*}
		\overline{\Af}&=\begin{pmatrix}
			\overline{\alpha}_1  & 0 \\ 0 & -(\overline{a}_1/\overline{\alpha}_1)u^{k_0-p}
		\end{pmatrix} & \overline{\Bf}&=\begin{pmatrix}
			(\overline{a_2b_2/p})u & 0 \\ 0 & -\frac{\overline{b}_1}{(\overline{a_2b_2/p})}u^{k_1}
		\end{pmatrix}.
	\end{align*}
\end{enumerate}

\noindent Hence, in both cases we will have a product
\[\overline{\Af}\varphi(\overline{\Bf})=\begin{pmatrix}
	(\overline{\frac{a_2b_2\alpha_1}{p}})u^p & 0 \\ 0 & (\overline{\frac{a_1b_1p}{a_2b_2\alpha_1}})u^{(k_0-p)+pk_1}
\end{pmatrix}.\]
We may now use Proposition \ref{compute-Hom-prop} to get a reducible reduction
\[\overline{V(A,B)}|_{G_\infty}=\chi_2\omega_2^{p}\oplus \chi_1\chi_2^{-1}\omega_2^{(k_0-p)+pk_1}\]
where $\chi_1=\chi(\overline{a_1b_1})$ and $\chi_2=\chi(\overline{a_2b_2\alpha_1/p})$.

\subsubsection*{Subcase: $\nu_p(a_2b_2)<1$}
Under the bounds $\nu_p(a_2b_2)<1$ we have reduced partial Frobenius from Figure \ref{red-frob-f=2-table} after scaling by $u$ via $X^{(f-1)}=\Diag(u,1)$,
\begin{align*}
	\overline{\Af}&=\begin{pmatrix}
		0 & -(\overline{a_1/\alpha_1})u^{k_0-p} \\ -\overline{\alpha}_1& 0
	\end{pmatrix} & \overline{\Bf}&=\begin{pmatrix}
		-\overline{b}_1u^{k_1+1} & 0 \\ 0 & 1
	\end{pmatrix}.
\end{align*}

\noindent Thankfully we can take the product immediately so that
\[\overline{\Af}\varphi(\overline{\Bf})=\begin{pmatrix}
	0 & -(\overline{a_1/\alpha_1})u^{k_0-p} \\ (\overline{b_1\alpha_1})u^{p(k_1+1)} & 0
\end{pmatrix}.\]
For the purposes of identifying the unramified character, consider the following base change over $\overline{\F}_p$ given by
\[\begin{pmatrix}
	1 & 0 \\ 0 & \frac{(\overline{-a_1b_1})^{1/2}}{\overline{b_1\alpha}_1}\end{pmatrix} * _\varphi\begin{pmatrix}
	0 & -(\overline{a_1/\alpha_1})u^{k_0-p} \\ (\overline{b_1\alpha_1})u^{p(k_1+1)} & 0
	\end{pmatrix}=(\overline{-a_1b_1})^{1/2}\begin{pmatrix}
	0 & u^{k_0-p} \\ u^{p(k_1+1)} & 0
	\end{pmatrix}.\]
This allows us to invoke Proposition \ref{compute-Hom-prop} to compute the irreducible reduction
\[\overline{V(A,B)}|_{G_\infty}=\chi((-\overline{a_1b_1})^{1/2})\ind_{G_{\Q_{p^{4}}}}^{G_{\Q_{p^{2}}}}\left(\omega_4^{(k_0-p)+p^2(k_1+1)}\right).\]

\subsubsection{Explicit Reductions for Type $(\II,\I)$}
Let us begin with the case of $\nu_p(a_2)>1$ so that we have reduced partial Frobenius from Figure \ref{red-frob-f=2-table} given by
\begin{align*}
	\overline{\Af}&=\begin{pmatrix}
		\overline{a}_1u^{k_0} & 0 \\ 0 & 1
	\end{pmatrix} & \overline{\Bf}&=\begin{pmatrix}
		0 & \overline{b}_1u^{k_1} \\ 1 & 0
	\end{pmatrix}.
\end{align*}

\noindent Proceeding as in Section \ref{verylarge-red-data-subsec}, we obtain
\begin{align*}
	\overline{\Af}&=\begin{pmatrix}
		\overline{a}_1u^{k_0-p} & 0 \\ 0 & 1
	\end{pmatrix} & \overline{\Bf}&=\begin{pmatrix}
		0 & \overline{b}_1u^{k_1+1} \\ 1 & 0
	\end{pmatrix}.
\end{align*}

\noindent The product will then be given by
\[\overline{\Af}\varphi(\overline{\Bf})=\begin{pmatrix}
	0 & (\overline{a_1b_1})u^{(k_0-p)+p(k_1+1)} \\ 1 & 0
\end{pmatrix}.\]
Scaling by $(\overline{a_1b_1})^{1/2}$ will make
\[\overline{\Af}\varphi(\overline{\Bf})=(\overline{a_1b_1})^{1/2}\begin{pmatrix}
	0 & u^{(k_0-p)+p(k_1+1)} \\ 1 & 0
\end{pmatrix}.\]
Hence, we may apply Proposition \ref{compute-Hom-prop} to get an irreducible reduction
\[\overline{V(A,B)}|_{G_\infty}=\chi((\overline{a_1b_1})^{1/2})\ind_{G_{\Q_{p^{4}}}}^{G_{\Q_{p^{2}}}}\left(\omega_4^{(k_0-p)+p(k_1+1)}\right)=\chi((\overline{a_1b_1})^{1/2})\ind_{G_{\Q_{p^{4}}}}^{G_{\Q_{p^{2}}}}\left(\omega_4^{k_0+pk_1}\right).\]

Next, we have $\nu_p(a_2)=1$ so that we get reduced partial Frobenius from Figure \ref{red-frob-f=2-table} given by
\begin{align*}
	\overline{\Af}&=\begin{pmatrix}
		\overline{a}_1u^{k_0} & 0 \\ (\overline{a_2/p})\overline{\alpha}_1 u^p & 1
	\end{pmatrix} & \overline{\Bf}&=\begin{pmatrix}
		0 & \overline{b}_1u^{k_1} \\ 1 & 0
	\end{pmatrix}.
\end{align*}

\noindent We proceed as in Section \ref{integral-red-data-subsec} to get 
\begin{align*}
	\overline{\Af}&=\begin{pmatrix}
		-(\overline{\frac{a_1 p}{a_2\alpha_1}})u^{k_0-p} & 0  \\ 0 & (\overline{a_2/p})\overline{\alpha}_1 
	\end{pmatrix} & \overline{\Bf}&=\begin{pmatrix}
		1 & 0 \\ 0 & \overline{b}_1u^{k_1+1}
	\end{pmatrix}.
\end{align*}

\noindent We can then easily see the product
\[\overline{\Af}\varphi(\overline{\Bf})=\begin{pmatrix}
	-(\overline{\frac{a_1 p}{a_2\alpha_1}})u^{k_0-p} & 0 \\ 0 & (\overline{\frac{b_1a_2\alpha_1}{p}})u^{p(k_1+1)}
\end{pmatrix}.\]
We can now use Proposition \ref{compute-Hom-prop} to get a reducible reduction
\[\overline{V(A,B)}|_{G_\infty}=\chi_1\chi_3^{-1}\omega_2^{k_0-p}\oplus \chi_2\chi_3\omega_2^{p(k_1+1)}\]
where $\chi_1=\chi(\overline{a}_1)$, $\chi_2=\chi(\overline{b}_1)$ and $\chi_3=\chi(\overline{a_2\alpha_1/p})$.

Now we consider the case of $\nu_p(a_2)<1$ but $\nu_p(b_2)>1$ so that by Figure \ref{red-frob-f=2-table} we have partial reduced Frobenius given by
\begin{align*}
	\overline{\Af}&=\begin{pmatrix}
		0 & -(\overline{a_1/\alpha_1})u^{k_0-p} \\ -\overline{\alpha}_1u^p  & 0
	\end{pmatrix} & \overline{\Bf}&=\begin{pmatrix}
		1 & 0 \\ 0 & \overline{b}_1u^{k_1}
	\end{pmatrix}.
\end{align*}

\noindent Following \ref{red-data-odd-subsec}, we have
\begin{align*}
	\overline{\Af}&=\begin{pmatrix}
		0 & -(\overline{a_1/\alpha_1})u^{k_0-p} \\ -\overline{\alpha}_1  & 0
	\end{pmatrix} & \overline{\Bf}&=\begin{pmatrix}
		u & 0 \\ 0 & \overline{b}_1u^{k_1}
	\end{pmatrix}.
\end{align*}

\noindent We may immediately compute the product
\[\overline{\Af}\varphi(\overline{\Bf})=\begin{pmatrix}
	0 & -(\overline{a_1b_1/\alpha_1})u^{(k_0-p)+pk_1} \\ -(\overline{\alpha}_1)u^{p} & 0
\end{pmatrix}.\]
Scaling by $(-(\overline{a_1b_1})^{1/2}/\overline{\alpha}_1)\in \overline{\F}_p$ will give 
\[\overline{\Af}\varphi(\overline{\Bf})=(\overline{a_1b_1})^{1/2}\begin{pmatrix}
	0 & u^{(k_0-p)+pk_1} \\ u^{p} & 0
\end{pmatrix}.\]
Hence we may use Proposition \ref{compute-Hom-prop} to get an irreducible reduction
\[\overline{V(A,B)}|_{G_\infty}=\chi((\overline{a_1b_1})^{1/2})\ind_{G_{\Q_{p^{4}}}}^{G_{\Q_{p^{2}}}}\left(\omega_4^{(k_0-p)+pk_1+p^2}\right).\]

The penultimate case is when $\nu_p(a_2)<1$ but $\nu_p(b_2)=1$ so that we have reduced partial Frobenius given by Figure \ref{red-frob-f=2-table}
\begin{align*}
	\overline{\Af}&=\begin{pmatrix}
		0 & -(\overline{a_1/\alpha_1})u^{k_0-p} \\ -\overline{\alpha}_1u^p & 0
	\end{pmatrix} & \overline{\Bf}&=\begin{pmatrix}
		1 & (\overline{b_2/p}) \\ 0 & \overline{b}_1 u^{k_1}
	\end{pmatrix}.
\end{align*}

\noindent Proceeding as in Section \ref{red-data-odd-subsec} will result in
\begin{align*}
	\overline{\Af}&=\begin{pmatrix}
		-\overline{\alpha}_1 & 0 \\ 0 & -(\overline{a_1/\alpha_1})u^{k_0-p}
	\end{pmatrix} & \overline{\Bf}&=\begin{pmatrix}
		(\overline{b_2/p})u & 0 \\ 0 & -(\overline{\frac{b_1p}{b_2}})u^{k_1}
	\end{pmatrix}.
\end{align*}

\noindent Hence, we have a product
\[\overline{\Af}\varphi(\overline{\Bf})=\begin{pmatrix}
	-(\overline{\frac{b_2\alpha_1}{p}})u^p & 0 \\ 0 & (\overline{\frac{a_1b_1p}{b_2\alpha_1}})u^{(k_0-p)+pk_1}
\end{pmatrix}.\]
We may now invoke Proposition \ref{compute-Hom-prop} to get a reducible reduction 
\[\overline{V(A,B)}|_{G_\infty}=\chi_2\omega_2^p\oplus \chi_1\chi_2^{-1}\omega_2^{(k_0-p)+pk_1}\]
where $\chi_1=(-\overline{a_1b_1})$ and $\chi_2=\chi(-\overline{b_2\alpha_1/p})$. Notice that up to restriction to inertia, this is exactly the same reduction that we got in the Type $(\I,\I)$ case with $\nu_p(a_2b_2)=1$.

Our final case is when both $\nu_p(b_2),\nu_p(a_2)<1$. Thankfully, we have already noticed that the reduced partial Frobenius from Figure \ref{red-frob-f=2-table} is exactly the same, up to a scaling operation, to the case of Type $(\I,\I)$ with $\nu_p(a_2b_2)<1$. As a result, we have an irreducible reduction
\[\overline{V(A,B)}|_{G_\infty}=\chi((\overline{a_1b_1})^{1/2})\ind_{G_{\Q_{p^{4}}}}^{G_{\Q_{p^{2}}}}\left(\omega_4^{(k_0-p)+p^2(k_1+1)}\right).\]

\subsubsection{Explicit Reductions for Type $(\I,\II)$}
For our final Type combination, we recall the shape of the reduced partial Frobenius from Figure \ref{red-frob-f=2-table} depend only on $\nu_p(a_2)$ as $\Ws=\emptyset$ in this case. Hence, we begin as before with the case of $\nu_p(a_2)>1$ so that we have partial reduced Frobenius given by
\begin{align*}
	\overline{\Af}&=\begin{pmatrix}
		0 & \overline{a}_1u^{k_0} \\ 1 & 0
	\end{pmatrix} & \overline{\Bf}&=\begin{pmatrix}
		\overline{b}_1u^{k_1} & 0 \\ 0 & 1
	\end{pmatrix}.
\end{align*}

\noindent Proceeding as in \ref{verylarge-red-data-subsec}, we will get
\begin{align*}
	\overline{\Af}&=\begin{pmatrix}
		0 & \overline{a}_1u^{k_0-p} \\ 1 & 0
	\end{pmatrix} & \overline{\Bf}&=\begin{pmatrix}
		\overline{b}_1u^{k_1} & 0 \\ 0 & u
	\end{pmatrix}.
\end{align*}

\noindent Hence, we may take the product to be
\[\overline{\Af}\varphi(\overline{\Bf})=\begin{pmatrix}
	0 & (\overline{a}_1)u^{(k_0-p)+p} \\ (\overline{b}_1)u^{pk_1} & 0
\end{pmatrix}.\]
Scaling the anti-diagonal of the above product by $((\overline{a_1b_1})^{1/2}/\overline{b}_1)\in\overline{\F}_p$ will make 
\[\overline{\Af}\varphi(\overline{\Bf})=(\overline{a_1b_1})^{1/2}\begin{pmatrix}
	0 & u^{(k_0-p)+p} \\ u^{pk_1} & 0
\end{pmatrix}.\]
We may now use Proposition \ref{compute-Hom-prop} to get an irreducible reduction
\[\overline{V(A,B)}|_{G_\infty}=\chi((\overline{a_1b_1})^{1/2})\ind_{G_{\Q_{p^{4}}}}^{G_{\Q_{p^{2}}}}\left(\omega_4^{(k_0-p)+p^2k_1+p}\right)=\chi((\overline{a_1b_1})^{1/2})\ind_{G_{\Q_{p^{4}}}}^{G_{\Q_{p^{2}}}}\left(\omega_4^{k_0+p^2k_1}\right).\]

We now consider the case of $\nu_p(a_2)=1$ so that by Figure \ref{red-frob-f=2-table} we have reduced partial Frobenius given by
\begin{align*}
	\overline{\Af}&=\begin{pmatrix}
		0 & \overline{a}_1u^{k_0} \\ 1 & (\overline{a_2/p})\overline{\alpha}_1 u^p
	\end{pmatrix} & \overline{\Bf}&=\begin{pmatrix}
		\overline{b}_1u^{k_1} & 0 \\ 0 & 1
	\end{pmatrix}.
\end{align*}

\noindent Proceeding as in Section \ref{integral-red-data-subsec}, we will receive
\begin{align*}
	\overline{\Af}&=\begin{pmatrix}
		-(\overline{\frac{a_1p}{a_2\alpha_1}})u^{k_0-p} & 0 \\ 0 & (\overline{a_2/p})\overline{\alpha}_1 
	\end{pmatrix} & \overline{\Bf}&=\begin{pmatrix}
		\overline{b}_1u^{k_1} & 0 \\ 0 & u
	\end{pmatrix}.
\end{align*}

\noindent We may now take the product 
\[\overline{\Af}\varphi(\overline{\Bf})=\begin{pmatrix}
	-(\overline{\frac{a_1b_1p}{a_2\alpha_1}})u^{(k_0-p)+pk_1} & 0 \\ 0 & (\overline{\frac{a_2\alpha_1}{p}}) u^p
\end{pmatrix}.\]
Hence, we may invoke Proposition \ref{red-frob-f=2-table} to get a reducible reduction
\[\overline{V(A,B)}|_{G_\infty}=\chi_2\omega_2^{p}\oplus \chi_1\chi_2^{-1}\omega_2^{(k_0-p)+pk_1}.\]
We may notice that the this reduction has already happened twice; The first time in Type $(\I,\I)$ with $\nu_p(a_2b_2)=1$ and the second in Type $(\II,\I)$ with $\nu_p(a_2)<1$ and $\nu_p(b_2)=1$.

Our final case is when $\nu_p(a_2)<1$. Thankfully, we have already noticed that the reduced partial Frobenius from Figure \ref{red-frob-f=2-table} is exactly the same, up to a scaling operation, to the case of Type $(\I,\I)$ with $\nu_p(a_2b_2)<1$ and Type $(\II,\I)$ with $\nu_p(a_2),\nu_p(b_2)<1$. As a result, we have an irreducible reduction
\[\overline{V(A,B)}|_{G_\infty}=\chi((\overline{a_1b_1})^{1/2})\ind_{G_{\Q_{p^{4}}}}^{G_{\Q_{p^{2}}}}\left(\omega_4^{(k_0-p)+p^2(k_1+1)}\right).\]

\subsection{Summary of Reductions}\label{summary-sec}
Above we have found, up to an unramified twist or restriction to inertia, the reduction data in the sense of Definition \ref{red-data-def} for Kisin modules $\overline{\Mf(\Af)}$ over $k_F[\![u]\!]$ with labeled heights in the range $p+2\le k_{0}\le 2p-4$ and $2\le k_i\le p-3$ for $1\le i\le f-1$. By the discussion found in Section \ref{comp-red-subsec}, this data suffices to compute the reduction of the associated crystalline representations $V(A)$. Moreover, such representations cover all possible irreducible, two-dimensional crystalline representations of $G_K$ with labeled Hodge-Tate weights in the same range and under certain bounds on the $p$-adic valuations of $a_2^{(i)}$. As a result, we have computed all such reductions and we summarize our results below in terms of the $p$-adic valuations of the $a_2^{(i)}$ parameters originating from the matrices $A=(A^{(i)})$ and the Type combination of said matrices.

\subsection*{Large Valuation $\nu_p(a_2^{(0)})>1$:}
Recall that $\Ss$ denotes the set of $i\in\Z/f\Z$ corresponding to Type $\I$ Frobenius and $\Ts$ is the set of $i\in\Z/f\Z$ corresponding to Type $\II$ Frobenius. 
\begin{itemize}
	\item $\nu_p(a_2^{(i)})>0$ for all $i\in\Ss$:
	\begin{align*}
		\mu_{0}&=\begin{dcases}
			(S, (0,k_{0}-p))& \text{if	} 0 \in\Ss\\
			(I, (k_{0}-p, 0)) & \text{if	}  0\in\Ts
		\end{dcases}\\
		\mu_i&=\begin{dcases}
			(S, (0,k_i))& \text{if	} 0<i< f-1 \text{ and } i \in\Ss\\
			(I, (k_i, 0)) & \text{if	} 0<i< f-1 \text{ and } i \in\Ts\\
		\end{dcases}\\
		\mu_{f-1}&=\begin{dcases}
			(S, (1,k_{f-1})) & \text{if	}f-1\in\Ss \text{ and } 0\in \Ss \\
			(S, (0,k_{f-1}+1)) & \text{if	}f-1\in\Ss \text{ and } 0\in \Ts \\
			(I, (k_{f-1},1)) & \text{if	}f-1\in\Ts \text{ and } 0\in \Ss \\
			(I, (k_{f-1}+1,0)) & \text{if	}f-1\in\Ts \text{ and } 0\in \Ts.
		\end{dcases}
	\end{align*}
	
	\noindent Note that this case is reducible if and only if $|\Ss|$ is even and it is irreducible otherwise.

	\item $\nu_p(a_2^{(i)})=0$ for some $i\in\Ss$:
	
	\noindent Recall that the subset $\Xs\subset\Z/f\Z$ is as defined in Definition \ref{X-sandwich-subset-def}.
	\begin{align*}
		\mu_{0}&=\begin{dcases}
			(I, (0, k_{0}-p)) & \text{if	}  0\in\Xs\\
			(I, (k_{0}-p, 0))& \text{if	}  0\notin\Xs
		\end{dcases}\\
		\mu_i&=\begin{dcases}
			(I, (0, k_i)) & \text{if	}0<i<f-1\text{ and }  i\in\Xs\\
			(I, (k_i,0))& \text{if	}  0<i<f-1\text{ and }i\notin\Xs
		\end{dcases}\\
		\mu_{f-1}&=\begin{dcases}
			(I, (0,k_{f-1}+1)) & \text{if }0,f-1\in\Xs \\
			(I, (1, k_{f-1})) & \text{if }f-1\in\Xs\text{ and }0\notin\Xs\\
			(I, (k_{f-1},1)) & \text{if }f-1\in\Xs\text{ and }0\in\Xs \\
			(I, (k_{f-1}0+1,0)) & \text{if }0,f-1\notin\Xs.
		\end{dcases}
	\end{align*}
	
	\noindent Note that this case is necessarily reducible.
\end{itemize}

\subsection*{Integer Valuation $\nu_p(a_2^{(0)})=1$:}
Recall that the subset $\Ys\subset\Z/f\Z$ is defined as in Definition \ref{Y-sandwich-subset-def}. We encourage the reader to notice the intersections of reductions between multiple Type combinations.
\begin{itemize}
	\item $0\in\Ss$ and $f-1\notin\Ys$ \textbf{OR} $0\in\Ts$ and $f-1\in\Ys$:
\begin{align*}
	\mu_{0}&=(I, (k_{0}-p, 0)).\\
	\mu_i&=\begin{dcases}
		(I, (k_i, 0)) & \text{if }1\le i<f-1 \text{ and }i\notin\Ys \\
		(I, (0, k_i)) & \text{if }1\le i<f-1 \text{ and }i\in \Ys
	\end{dcases}\\
	\mu_{f-1}&=\begin{dcases}
		(I, (k_{f-1}, 1)) & \text{if }f-1\notin\Ys \\
		(I, (0, k_{f-1}+1)) & \text{if }f-1\in \Ys.
	\end{dcases}
	\end{align*}

\noindent Note that this case is necessarily reducible.
	
	\item $0\in\Ts$ and $f-1\notin\Ys$ \textbf{OR} $0\in\Ss$ and $f-1\in\Ys$:
	\begin{align*}
		\mu_{0}&=(I, (k_{0}-p, 0))\\
		\mu_i&=\begin{dcases}
			(I, (k_i, 0)) & \text{if }1\le i<f-1 \text{ and }i\notin\Ys \\
			(I, (0, k_i)) & \text{if }1\le i<f-1 \text{ and }i\in \Ys
		\end{dcases}\\
		\mu_{f-1}&=\begin{dcases}
			(I, (k_{f-1}+1, 0)) & \text{if }f-1\notin\Ys \\
			(I, (1, k_{f-1})) & \text{if }f-1\in \Ys.
		\end{dcases}
	\end{align*}
	
	\noindent Note that this case is necessarily reducible.
\end{itemize}

\subsection*{Small Valuation $\nu_p(a_2^{(0)})<1$:}
Recall that the subset $\Ws\subset\Z/f\Z$ is defined in Definition \ref{sandwich-subset-def}. We encourage the reader to notice the intersections of reductions between multiple Type combinations.
\begin{itemize}
	\item $|\Ss|$ even and $\nu_p(a_2^{(i)})>1-\nu_p(a_2^{(0)})$ for all $i\in\Ws$ \textbf{OR} $|\Ss|$ odd and and $\nu_p(a_2^{(i)})>1$ for all $i\in\Ws$ with $\nu_p(a_2^{(0)})=0$:
	\begin{align*}
		\mu_{0}&=(I, (k_{0}-p,0))\\
		\mu_i&=\begin{dcases}
			(I, (0, k_i)) & \text{if	} 0<i< f-1 \text{ and }i\in\Ws \\
			(I, (k_i,0))&\text{if	} 0<i< f-1 \text{ and }i\notin\Ws 
		\end{dcases}\\
		\mu_{f-1}&=\begin{dcases}
			(I, (1, k_{f-1})) & \text{if	} f-1\in\Ws \\
			(I, (k_{f-1}+1,0)) &\text{if	} f-1\notin\Ws.
		\end{dcases}
	\end{align*}
	
	\noindent Note that this case is necessarily reducible.
	
	\item $|\Ss|$ even and $r\in\Ws$ is the largest integer with $\nu_p(a_2^{(r)})=1-\nu_p(a_2^{(0)})$ \textbf{OR} $|\Ss|$ odd and $r\in\Ws$ is the largest integer with $\nu_p(a_2^{(r)})=1$:
\begin{align*}
	\mu_{0}&=(I, (0,k_{0}-p))\\
	\mu_j&=\begin{dcases}
		(I, (k_j, 0))& \text{if }0<j<r \text{ and }j\in\Ws \\
		(I, (0,k_j))& \text{if }0<j<r \text{ and }j\notin\Ws
	\end{dcases}\\
	\mu_i&=\begin{dcases}
		(I, (0, k_i)) & \text{if }r\le i<f-1 \text{ and }i\in\Ws \\
		(I, (k_i, 0))& \text{if	}r\le i<f-1 \text{ and }i\notin\Ws
	\end{dcases}\\
	\mu_{f-1}&=\begin{dcases}
		(I, (1, k_{f-1})) & \text{if }f-1\in\Ws \\
		(I, (k_{f-1}+1, 0))& \text{if	}f-1\notin\Ws.
	\end{dcases}
\end{align*}

	\noindent Note that this case is necessarily reducible.
	
	\item $|\Ss|$ odd and $\nu_p(a_2^{(i)})>1$ for all $i\in\Ws$ with $\nu_p(a_2^{(0)})\neq 0$:
\begin{align*}
	\mu_{0}&=(S, (0 ,k_{0}-p))\\
	\mu_i&=\begin{dcases}
		(I, (0, k_i)) & \text{if } 0<i< f-1 \text{ and }i\in\Ws\\
		(I, (k_i, 0)) & \text{if } 0<i< f-1 \text{ and }i\notin\Ws
	\end{dcases}\\
	\mu_{f-1}&=\begin{dcases}
		(I, (1, k_{f-1})) & \text{if	} f-1\in\Ws\\
		(I, (k_{f-1}+1, 0)) & \text{if	} f-1\notin\Ws.
	\end{dcases}
\end{align*}

	\noindent Note that this case is necessarily irreducible.
\end{itemize}

\subsection*{Reductions for $f=2$:}

When $f=2$, we are able to explicitly compute the reductions along with the associated unramified characters. As before, the shape of these reductions turn out to depend on both the $p$-adic valuations of the parameters $a_2$ and $b_2$ along with the three distinct Type combinations. The table above summarizes these computations. Recall that $\chi(c)$ denotes the unramified character sending Frobenius to $c\in\overline{\F}_p$. All inductions are to be taken from $G_\Qpfour$ to $G_\Qptwo$.

\begin{figure}[H]
	\centering
	\renewcommand{\arraystretch}{1.5}
	\begin{tabular}{|c|l|l|l|}
		\hline
		Type & Conditions on valuations & Reduction $\overline{V(A,B)}|_{G_\infty}$ & Unramified Characters\\
		\hline
		$(\I,\I)$ & 
		$\nu_p(a_2b_2)> 1$ & $
		\chi_1\omega_2^{k_0}\oplus \chi_2\omega_2^{pk_1}$ & 
		$\begin{aligned}
			\chi_1 &=\chi(\overline{a}_1) \\ \chi_2 &=\chi(\overline{b}_1)
		\end{aligned}$ \\
		\hline
		
		$(\I,\I)$ & 
		$\nu_p(a_2b_2)= 1$ & 
		$\chi_1\chi_2^{-1}\omega_2^{(k_0-p)+pk_1}\oplus \chi_2\omega_2^p$ & 
		$\begin{aligned}
			\chi_1&=\chi(\overline{a_1b_1}) \\ \chi_2&=\chi(\overline{a_2b_2\alpha_1/p})
		\end{aligned}$\\
		\hline
		
		$(\I,\I)$ & 
		$\nu_p(a_2b_2)<1$ & 
		$\chi\ind\left(\omega_4^{(k_0-p)+p^2(k_1+1)}\right)$&  $\chi=\chi((-\overline{a_1b_1})^{1/2})$\\
		\hline
		
		$(\II,\I)$ & 
		$\nu_p(a_2)> 1$ &  
		$\chi\ind\left(\omega_4^{k_0+pk_1}\right)$& $\chi=\chi((\overline{a_1b_1})^{1/2})$\\
		\hline
		
		$(\II,\I)$ & 
		$\nu_p(a_2)= 1$ & 
		$\chi_1\chi_3^{-1}\omega_2^{k_0-p}\oplus \chi_2\chi_3\omega_2^{p(k_1+1)}$& $\begin{aligned}
			\chi_1&=\chi(\overline{a}_1) \\ \chi_2&=\chi(\overline{b}_1) \\ \chi_3&=\chi(\overline{a_2\alpha_1/p})
		\end{aligned}$\\
		\hline
		
		$(\II,\I)$ & 
		$\nu_p(a_2)< 1$ and 
		$\nu_p(b_2)> 1$ & 
		$\chi\ind\left(\omega_4^{(k_0-p)+p(k_1+p)}\right)$& $\chi=\chi((\overline{a_1b_1})^{1/2})$\\
		\hline
		
		$(\II,\I)$ & 
		$\nu_p(a_2)< 1$ and $\nu_p(b_2)= 1$ & 
		$\chi_1\chi_2^{-1}\omega_2^{(k_0-p)+pk_1}\oplus\chi_2\omega_2^p$& 
		$\begin{aligned}
			\chi_1&=\chi(-\overline{a_1b_1}) \\ \chi_2&=\chi(-\overline{b_2\alpha_1/p})
		\end{aligned}$\\
		\hline
		
		$(\II,\I)$ & 
		$\nu_p(a_2)< 1$ and $\nu_p(b_2)< 1$ & $\chi\ind\left(\omega_4^{(k_0-p)+p^2(k_1+1)}\right)$& $\chi=\chi((\overline{a_1b_1})^{1/2})$\\
		\hline
		
		$(\I,\II)$ & 
		$\nu_p(a_2)> 1$ & 
		$\chi\ind\left(\omega_4^{k_0+p^2k_1}\right)$& $\chi=\chi((\overline{a_1b_1})^{1/2})$\\
		\hline
		
		$(\I,\II)$ & 
		$\nu_p(a_2)= 1$ & 
		$\chi_1\chi_2^{-1}\omega_2^{(k_0-p)+pk_1}\oplus\chi_2\omega_2^{p}$ & 
		$\begin{aligned}
			\chi_1&=\chi(-\overline{a_1b_1}) \\ \chi_2&=\chi(\overline{a_2\alpha_1/p})
		\end{aligned}$\\
		\hline
		
		$(\I,\II)$ & 
		$\nu_p(a_2)< 1$ &$\chi\ind\left(\omega_4^{(k_0-p)+p^2(k_1+1)}\right)$& $\chi=\chi((\overline{a_1b_1})^{1/2})$ \\
		\hline
	\end{tabular}
	\caption{Reductions in $f=2$ depending on Type and $p$-adic valuations of $a_2$ and $b_2$.}\label{f=2-red-sum-fig}
\end{figure}

Let us explicitly point out the intersections of reductions up to unramified characters.
\begin{itemize}
	\item The reducible reduction $\omega_2^{(k_0-p)+pk_1}\oplus\omega_2^{p}$ occurs three times:
	\begin{itemize}
		\item Type $(\I,\I)$ with $\nu_p(a_2b_2)=1$ which if we recall from Section \ref{red-f=2-sec} happened when $\nu_p(a_2)=1$ and $\nu_p(b_2)=0$ \textit{or} $\nu_p(a_2)<1$ and $\nu_p(b_2)=1-\nu_p(a_2)$;
		\item Type $(\II,\I)$ with $\nu_p(a_2)<1$ and $\nu_p(b_2)=1$;
		\item Type $(\I,\II)$ with $\nu_p(a_2)=1$.
	\end{itemize}
		\item The irreducible reduction $\ind\omega_4^{(k_0-p)+p^2(k_1+1)}$ occurs three times:
	\begin{itemize}
		\item Type $(\I,\I)$ with $\nu_p(a_2b_2)<1$;
		\item Type $(\II,\I)$ with $\nu_p(a_2)<1$ and $\nu_p(b_2)<1$;
		\item Type $(\I,\II)$ with $\nu_p(a_2)<1$.
	\end{itemize}
\end{itemize}

%%%%%%%%%%%%%%%%%%%%%%%%%%%%%%%%%%%%%%%%%%%%%%%%%%%%%%%%%%%%%%%%%%%%%%%%%%%%%%%%%%%%%%
%%%%%%%%%%%%%%%%%%%%%%%%%%%%%%%%%%%%%%%%%%%%%%%%%%%%%%%%%%%%%%%%%%%%%%%%%%%%%%%%%%%%%%
%%%%%%%%%%%%%%%%%%%%%%%%%%%%%%%%%%%%%%%%%%%%%%%%%%%%%%%%%%%%%%%%%%%%%%%%%%%%%%%%%%%%%%

\printbibliography[title = {References}]

%%%%%%%%%%%%%%%%%%%%%%%%%%%%%%%%%%%%%%%%%%%%%%%%%%%%%%%%%%%%%%%%%%%%%%%%%%%%%%%%%%%%%%
%%%%%%%%%%%%%%%%%%%%%%%%%%%%%%%%%%%%%%%%%%%%%%%%%%%%%%%%%%%%%%%%%%%%%%%%%%%%%%%%%%%%%%
%%%%%%%%%%%%%%%%%%%%%%%%%%%%%%%%%%%%%%%%%%%%%%%%%%%%%%%%%%%%%%%%%%%%%%%%%%%%%%%%%%%%%%

\end{document}